# NEUTROSOPHIC RINGS

W. B. Vasantha Kandasamy
Florentin Smarandache

# NEUTROSOPHIC RINGS

**Dr. W. B. Vasantha Kandasamy**
Department of Mathematics
Indian Institute of Technology Madras
Chennai 600 036
India
vasanthakandasamy@gmail.com
http://www.vasantha.net
http://mat.iitm.ac.in/~wbv

**Dr. Florentin Smarandache**
Department of Mathematics
University of New Mexico
Gallup, NM, USA
smarand@unm.edu

2006

# CONTENTS

Chapter Four
**SUGGESTED PROBLEMS** 107

**REFERENCES** 135

**INDEX** 149

**ABOUT THE AUTHORS** 154

# PREFACE

In this book we define the new notion of neutrosophic rings. The motivation for this study is two-fold. Firstly, the classes of neutrosophic rings defined in this book are generalization of the two well-known classes of rings: group rings and semigroup rings. The study of these generalized neutrosophic rings will give more results for researchers interested in group rings and semigroup rings. Secondly, the notion of neutrosophic polynomial rings will cause a paradigm shift in the general polynomial rings. This study has to make several changes in case of neutrosophic polynomial rings. This would give solutions to polynomial equations for which the roots can be indeterminates. Further, the notion of neutrosophic matrix rings is defined in this book. Already these neutrosophic matrixes have been applied and used in the neutrosophic models like neutrosophic cognitive maps (NCMs), neutrosophic relational maps (NRMs) and so on.

This book has four chapters. Chapter one is introductory in nature, for it recalls some basic definitions essential to make the book a self-contained one. Chapter two, introduces for the first time the new notion of neutrosophic rings and some special neutrosophic rings like neutrosophic ring of matrix and neutrosophic polynomial rings. Chapter three gives some new classes of neutrosophic rings like group neutrosophic rings, neutrosophic group neutrosophic rings, semigroup neutrosophic rings, S-semigroup neutrosophic rings which can be realized as a type of extension of group rings or generalization of group rings. Study of these structures will throw light on the research on the algebraic structure of group rings. Chapter four is entirely devoted to the problems on this new topic, which is an added attraction to researchers.

A salient feature of this book is that it gives 246 problems in Chapter four. Some of the problems are direct and simple, some little difficult and some can be taken up as a research problem.

We express our sincere thanks to Kama for her help in the layout and Meena for cover-design of the book. The authors express their whole-hearted gratefulness to Dr.K.Kandasamy whose invaluable support and help, and patient proofreading contributed to a great extent to the making of this book.

W.B.VASANTHA KANDASAMY
FLORENTIN SMARANDACHE
1 August 2006

Chapter One

# INTRODUCTION

In this book we assume all fields to be real fields of characteristic 0, all vector spaces are taken as real spaces over reals and we denote the indeterminacy by 'I' as i will make a confusion, as it denotes the imaginary value, viz. $i^2 = -1$ that is $\sqrt{-1} = i$. The indeterminacy I is such that I . I = I $^2$ = I. In this book the authors mainly use the notion of neutrosophic semigroup to construct special type of neutrosophic rings viz. neutrosophic semigroup ring and neutrosophic semigroup neutrosophic ring.

In this chapter we just recall some of the basic neutrosophic structures in this book. This chapter has three sections. In section one, we just recall the notion of neutrosophic groups and their properties, section two recalls the basic concept of neutrosophic semigroups and section three gives the definition of neutrosophic field.

## 1.1 Neutrosophic Groups and their Properties

In this section we recall the notion of neutrosophic groups introduced in [142-3]; neutrosophic groups in general do not have group structure. We also define yet another notion called pseudo neutrosophic groups which have group structure. As neutrosophic groups do not have group structure the classical theorems viz. Sylow, Lagrange or Cauchy are not true in general which forces us to define notions like Lagrange neutrosophic groups, Sylow neutrosophic groups and Cauchy elements. Examples are given for the understanding of these new concepts.

We just give the basic definition alone as we use it only in the construction of neutrosophic group rings and neutrosophic group neutrosophic rings which are analogous structure of group rings. In fact neutrosophic group neutrosophic rings forms the most generalized neutrosophic rings.

**DEFINITION 1.1.1:** *Let (G, *) be any group, the neutrosophic group is generated by I and G under * denoted by N(G) = {⟨G ∪I⟩, *}.*

**Example 1.1.1:** Let $Z_7$ = {0, 1, 2, …, 6} be a group under addition modulo 7. N(G) = {⟨$Z_7$ ∪ I⟩, '+' modulo 7} is a neutrosophic group which is in fact a group. For N(G) = {a + bI / a, b ∈ $Z_7$} is a group under '+' modulo 7. Thus this neutrosophic group is also a group.

**Example 1.1.2:** Consider the set G = $Z_5$ \ {0}, G is a group under multiplication modulo 5. N(G) = {⟨G ∪ I⟩, under the binary operation, multiplication modulo 5}. N(G) is called the neutrosophic group generated by G ∪ I. Clearly N(G) is not a group, for $I^2$ = I and I is not the identity but only an indeterminate, but N(G) is defined as the neutrosophic group.

Thus based on this we have the following theorem:

**THEOREM 1.1.1:** *Let (G, *) be a group, N(G) = {⟨G ∪I⟩, *} be the neutrosophic group.*

1. *N(G) in general is not a group.*
2. *N(G) always contains a group.*

*Proof:* To prove N(G) in general is not a group it is sufficient we give an example; consider ⟨$Z_5$ \ {0} ∪ I⟩ = G = {1, 2, 4, 3, I,

2 I, 4 I, 3 I}; G is not a group under multiplication modulo 5. In fact {1, 2, 3, 4} is a group under multiplication modulo 5.

N(G) the neutrosophic group will always contain a group because we generate the neutrosophic group N(G) using the group G and I. So G $\subseteq_{\neq}$ N(G); hence N(G) will always contain a group.

Now we proceed onto define the notion of neutrosophic subgroup of a neutrosophic group.

**DEFINITION 1.1.2:** *Let N(G) = ⟨G ∪ I⟩ be a neutrosophic group generated by G and I. A proper subset P(G) is said to be a neutrosophic subgroup if P(G) is a neutrosophic group i.e. P(G) must contain a (sub) group.*

***Example 1.1.3:*** Let N($Z_2$) = ⟨$Z_2$ ∪ I⟩ be a neutrosophic group under addition. N($Z_2$) = {0, 1, I, 1 + I}. Now we see {0, I} is a group under + in fact a neutrosophic group {0, 1 + I} is a group under '+' but we call {0, I} or {0, 1 + I} only as pseudo neutrosophic groups for they do not have a proper subset which is a group. So {0, I} and {0, 1 + I} will be only called as pseudo neutrosophic groups (subgroups).

We can thus define a pseudo neutrosophic group as a neutrosophic group, which does not contain a proper subset which is a group. Pseudo neutrosophic subgroups can be found as a substructure of neutrosophic groups. Thus a pseudo neutrosophic group though has a group structure is not a neutrosophic group and a neutrosophic group cannot be a pseudo neutrosophic group. Both the concepts are different.

Now we see a neutrosophic group can have substructures which are pseudo neutrosophic groups which is evident from the following example.

***Example 1.1.4:*** Let N($Z_4$) = ⟨$Z_4$ ∪ I⟩ be a neutrosophic group under addition modulo 4. ⟨$Z_4$ ∪ I⟩ = {0, 1, 2, 3, I, 1 + I, 2I, 3I, 1 + 2I, 1 + 3I, 2 + I, 2 + 2I, 2 + 3I, 3 + I, 3 + 2I, 3 + 3I}. o(⟨$Z_4$ ∪ I⟩) = $4^2$.

Thus neutrosophic group has both neutrosophic subgroups and pseudo neutrosophic subgroups. For T = {0, 2, 2 + 2I, 2I} is a neutrosophic subgroup as {0 2} is a subgroup of $Z_4$ under addition modulo 4. P = {0, 2I} is a pseudo neutrosophic group under '+' modulo 4.

Now we are not sure that general properties, which are true of groups, are true in case of neutrosophic groups for neutrosophic groups are not in general groups. We see that in case of finite neutrosophic groups the order of both neutrosophic subgroups and pseudo neutrosophic subgroups do not divide the order of the neutrosophic group. Thus we give some problems in the chapter 4.

**THEOREM 1.1.2:** *Neutrosophic groups can have non-trivial idempotents.*

*Proof:* For I ∈ N(G) we have $I^2 = I$ .

***Note:*** We cannot claim from this that N(G) can have zero divisors because of the idempotent as our neutrosophic groups are algebraic structures with only one binary operation multiplication in this case $I^2 = I$ .

We illustrate this by the following examples:

***Example 1.1.5:*** Let N(G) = {1, 2, I, 2I} a neutrosophic group under multiplication modulo three. We see $(2I)^2 \equiv I$ (mod 3), $I^2$ = I. (2I) I = 2I, $2^2 \equiv 1$ (mod 3). So P = {1, I, 2I} is a pseudo neutrosophic subgroup. Also o(P) ∤ o(N (G)).

Thus we see order of a pseudo neutrosophic subgroup need not in general divide the order of the neutrosophic group.

We give yet another example, which will help us to see that Lagrange's theorem for finite groups in case of finite neutrosophic groups, is not true.

***Example 1.1.6:*** Let N(G) = {1, 2, 3, 4, I, 2I, 3I, 4I} be a neutrosophic group under multiplication modulo 5. Now consider P = {1, 4, I, 2I, 3I, 4I} ⊂ N(G). P is a neutrosophic subgroup. o(N(G)) = 8 but o(P) = 6, 6 ∤ 8. So clearly neutrosophic groups in general do not satisfy the Lagrange theorem for finite groups.

So we define or characterize those neutrosophic groups, which satisfy Lagrange theorem as follows:

**DEFINITION 1.1.3:** *Let N(G) be a neutrosophic group. The number of distinct elements in N(G) is called the order of N(G). If the number of elements in N(G) is finite we call N(G) a finite neutrosophic group; otherwise we call N(G) an infinite neutrosophic group, we denote the order of N(G) by o(N(G)) or |N(G)|.*

**DEFINITION 1.1.4:** *Let N(G) be a finite neutrosophic group. Let P be a proper subset of N(G), which under the operations of N(G) is a neutrosophic group. If o(P) / o(N(G)) then we call P to be a Lagrange neutrosophic subgroup. If in a finite neutrosophic group all its neutrosophic subgroups are Lagrange then we call N(G) to be a Lagrange neutrosophic group.*

*If N(G) has atleast one Lagrange neutrosophic subgroup then we call N(G) to be a weakly Lagrange neutrosophic group. If N(G) has no Lagrange neutrosophic subgroup then we call N(G) to be a Lagrange free neutrosophic group.*

We have already given examples of these. Now we proceed on to recall the notion called pseudo Lagrange neutrosophic group.

**DEFINITION 1.1.5:** *Let N(G) be a finite neutrosophic group. Suppose L is a pseudo neutrosophic subgroup of N(G) and if o(L) / o(N(G)) then we call L to be a pseudo Lagrange neutrosophic subgroup. If all pseudo neutrosophic subgroups of N(G) are pseudo Lagrange neutrosophic subgroups then we call N(G) to be a pseudo Lagrange neutrosophic group.*

*If N(G) has atleast one pseudo Lagrange neutrosophic subgroup then we call N(G) to be a weakly pseudo Lagrange neutrosophic group. If N(G) has no pseudo Lagrange*

*neutrosophic subgroup then we call N(G) to be pseudo Lagrange free neutrosophic group.*

Now we illustrate by some example some more properties of neutrosophic groups, which paves way for more definitions. We have heard about torsion elements and torsion free elements of a group.

We in this book recall the definition of neutrosophic element and neutrosophic free element of a neutrosophic group, which will be used in characterizing zero divisors in case of neutrosophic group rings, and neutrosophic group neutrosophic rings.

**DEFINITION 1.1.6:** *Let N(G) be a neutrosophic group. An element $x \in N(G)$ is said to be a neutrosophic element if there exists a positive integer n such that $x^n = I$, if for any y a neutrosophic element no such n exists then we call y to be a neutrosophic free element.*

We illustrate these by the following examples:

***Example 1.1.7:*** Let N(G) = {1, 2, 3, 4, 5, 6, I, 2I, 3I, 4I, 5I, 6I} be a neutrosophic group under multiplication modulo 7. We have $(3I)^6 = I$, $(4I)^3 = I$ $(6I)^2 = I$, $I^2 = I$, $(2I)^6 = I$, $(5I)^6 = I$. In this neutrosophic group all elements are either torsion elements or neutrosophic elements.

***Example 1.1.8:*** Let us now consider the set {1, 2, 3, 4, I, 2I, 3I, 4I, 1 + I, 2 + I, 3 + I, 4 + I, 1 + 2I, 1 + 3I, 1 + 4I, 2 + 2I, 2 + 3I, 2 + 4I, 3 + 2I, 3 + 3I, 3 + 4I, 4 + 2I, 4 + 3I, 4 + 4I}. This is a neutrosophic group under addition modulo 5. For {1, 2, 3, 4} = $Z_5 \setminus \{0\}$ is group under addition modulo 5.

**DEFINITION 1.1.7:** *Let N(G) be a neutrosophic group under multiplication. Let $x \in N(G)$ be a neutrosophic element such that $x^m = 1$ then x is called the pseudo neutrosophic torsion element of N(G).*

In the above example we have given several pseudo neutrosophic torsion elements of N(G).

Now we proceed on to define Cauchy neutrosophic elements of a neutrosophic group N(G).

**DEFINITION 1.1.8:** *Let N(G) be a finite neutrosophic group. Let $x \in N(G)$, if x is a torsion element say $x^m = 1$ and if m/o N(G)) we call x a Cauchy element of N(G); if x is a neutrosophic element and $x^t = I$ with t / o(N(G)), we call x a Cauchy neutrosophic element of N(G). If all torsion elements of N(G) are Cauchy we call N(G) as a Cauchy neutrosophic group. If every neutrosophic element is a neutrosophic Cauchy element then we call the neutrosophic group to be a Cauchy neutrosophic, neutrosophic group.*

We now illustrate these concepts by the following examples:

***Example 1.1.9:*** Let N(G) = {0, 1, 2, 3, 4, I, 2I, 3I, 4I} be a neutrosophic group under multiplication modulo 5. {1, 2, 3, 4} is a group under multiplication modulo 5. Now we see o (N (G)) = 9, $4^2 \equiv 1$ (mod 5) 2 ⫮ 9 similarly $(3I)^4 = I$ but 4 ⫮ 9. Thus none of these elements are Cauchy elements or Cauchy neutrosophic Cauchy elements of N(G).

Now we give yet another example.

***Example 1.1.10:*** Let N(G) be a neutrosophic group of finite order 4 where N(G) = {1, 2, I, 2I} group under multiplication modulo 3. Clearly every element in N(G) is either a Cauchy neutrosophic element or a Cauchy element.

Thus we give yet another definition.

**DEFINITION 1.1.9:** *Let N(G) be a neutrosophic group. If every element in N(G) is either a Cauchy neutrosophic element of N(G) or a Cauchy element of N(G) then we call N(G) a strong Cauchy neutrosophic group.*

The above example is an instance of a strong Cauchy neutrosophic group. Now we proceed on to define the notion of p-Sylow neutrosophic subgroup, Sylow neutrosophic group, weak Sylow neutrosophic group and Sylow free neutrosophic group.

**DEFINITION 1.1.10:** *Let N(G) be a finite neutrosophic group. If for a prime p, $p^\alpha$ / o(N(G)) and $p^{\alpha+1} \nmid o(N(G))$, N(G) has a neutrosophic subgroup P of order $p^\alpha$ then we call P a p-Sylow neutrosophic subgroup of N(G).*

*Now if for every prime p such that $p^\alpha$ / o(N(G)) and $p^{\alpha+1} \nmid o(N(G))$ we have an associated p-Sylow neutrosophic subgroup then we call N(G) a Sylow neutrosophic group.*

*If N(G) has atleast one p-Sylow neutrosophic subgroup then we call N(G) a weakly Sylow neutrosophic group. If N(G) has no p-Sylow neutrosophic subgroup then we call N(G) a Sylow free neutrosophic group.*

Now unlike in groups we have to speak about Sylow notion associated with pseudo neutrosophic groups.

**DEFINITION 1.1.11:** *Let N(G) be a finite neutrosophic group. Let P be a pseudo neutrosophic subgroup of N (G) such that $o(P) = p^\alpha$ where $p^\alpha$ / o (N(G)) and $p^{\alpha+1} \nmid o(N(G))$, p a prime, then we call P to be a p-Sylow pseudo neutrosophic subgroup of N(G).*

*If for a prime p we have a pseudo neutrosophic subgroup P such that $o(P) = p^\alpha$ where $p^\alpha$ / o(N(G)) and $p^{\alpha+1} \nmid o(N(G))$, then we call P to be p-Sylow pseudo neutrosophic subgroup of N(G). If for every prime p such that $p^\alpha$ / o(N(G)) and $p^{\alpha+1} \nmid o(N(G))$, we have a p-Sylow pseudo neutrosophic subgroup then we call N(G) a Sylow pseudo neutrosophic group.*

*If on the other hand N(G) has atleast one p-Sylow pseudo neutrosophic subgroup then we call N(G) a weak Sylow pseudo neutrosophic group. If N(G) has no p-Sylow pseudo neutrosophic subgroup then we call N(G) a free Sylow pseudo neutrosophic group.*

Now we proceed on to define neutrosophic normal subgroup.

**DEFINITION 1.1.12:** *Let N(G) be a neutrosophic group. Let P and K be any two neutrosophic subgroups of N(G). We say P and K are neutrosophic conjugate if we can find x, y ∈ N(G) with x P = K y.*

We illustrate this by the following example:

**Example 1.1.11:** Let N(G) = {0, 1, 2, 3, 4, 5, I, 2I, 3I, 4I, 5I, 1 + I, 2 + I, 3 + I, …, 5 + 5I} be a neutrosophic group under addition modulo 6. P = {0, 3, 3I, 3+3I} is a neutrosophic subgroup of N(G). K = {0, 2, 4, 2 + 2I, 4 + 4I, 2I, 4I} is a neutrosophic subgroup of N(G). For 2, 3 in N(G) we have 2P = 3K = {0}. So P and K are neutrosophic conjugate.

　　Thus in case of neutrosophic conjugate subgroups K and P we do not demand o(K) = o(P).

　　Now we proceed on to define neutrosophic normal subgroups.

**DEFINITION 1.1.13:** *Let N(G) be a neutrosophic group. We say a neutrosophic subgroup H of N(G) is normal if we can find x and y in N(G) such that H = xHy for all x, y ∈ N(G) (Note x = y or y = x⁻¹ can also occur).*

**Example 1.1.12:** Let N(G) be a neutrosophic group given by N(G) = {0, 1, 2, 3, 4, 5, 6, 7, I, 2I, 3I, 4I, 5I, 6I, 7I} be a neutrosophic group under multiplication modulo 8.

　　H = {1, 7, I, 7I} is a neutrosophic subgroup of N(G). For no x, y ∈ N(G), xHy = H so H is not normal in N(G).

**DEFINITION 1.1.14:** *A neutrosophic group N(G) which has no nontrivial neutrosophic normal subgroups is called a simple neutrosophic group.*

Now we define pseudo simple neutrosophic groups.

**DEFINITION 1.1.15:** *Let N(G) be a neutrosophic group. A proper pseudo neutrosophic subgroup P of N(G) is said to be normal if we have P = xPy for all x, y ∈ N(G). A neutrosophic* *group is said to be pseudo simple neutrosophic group if N(G) has no nontrivial pseudo normal subgroups.*

We do not know whether there exists any relation between pseudo simple neutrosophic groups and simple neutrosophic groups.

Now we proceed on to define the notion of right (left) coset for both the types of subgroups.

**DEFINITION 1.1.16:** *Let N(G) be a neutrosophic group. H be a neutrosophic subgroup of N(G) for n ∈ N(G), then H n = {hn / h ∈ H} is called a right coset of H in G.*

Similarly we can define left coset of the neutrosophic subgroup H in G.

It is important to note that as in case of groups we cannot speak of the properties of neutrosophic groups as we cannot find inverse for every x ∈ N (G).

So we make some modification before which we illustrate these concepts by the following examples:

***Example 1.1.13:*** Let N(G) = {1, 2, 3, 4, I, 2I, 3I, 4I} be a neutrosophic group under multiplication modulo 5. Let H = {1, 4, I, 4 I} be a neutrosophic subgroup of N(G). The right cosets of H are as follows:

|        |   |       |   |                  |
|--------|---|-------|---|------------------|
| H.2    | = | {2, 3, 2I, 3I}, |   |       |
| H.3    | = | {3, 2, 3I, 2I}, |   |       |
| H.1    | = | H4    | = | {1, 4, I, 4I},   |
| H. I   | = | {I 4I}  | = | H. 4I,          |
| H.2 I  | = | {2I, 3I} | = | H 3I = {3I, 2I}. |

Therefore the classes are

|      |   |      |   |                  |
|------|---|------|---|------------------|
| [2]  | = | [3]  | = | {2, 3, 2I, 3I}   |
| [1]  | = | [4]  | = | H = {1, 4, I 4I} |
| [I]  | = | [4I] | = | R {I, 4I}        |
| [2I] | = | [3I] | = | {3I, 2I}.        |

Now we are yet to know whether they will partition N(G) for we see here the cosets do not partition the neutrosophic group.

That is why we had problems with Lagrange theorem so only we defined the notion of Lagrange neutrosophic group.

We give yet another example before which we define the concept of commutative neutrosophic group.

**DEFINITION 1.1.17:** *Let N(G) be a neutrosophic group. We say N(G) is a commutative neutrosophic group if for every pair a, b ∈ N(G), a b = b a.*

We have seen several examples of commutative neutrosophic groups. So now we give an example of a non-commutative neutrosophic group.

***Example 1.1.14:*** Let

$$N(G) = \left\{ \begin{pmatrix} a & b \\ c & d \end{pmatrix} \mid a, b, c, d \in \{0, 1, 2, I, 2I\} \right\}.$$

N(G) under matrix multiplication modulo 3 is a neutrosophic group which is non commutative.
We now give yet another example of cosets in neutrosophic groups.

***Example 1.1.15:*** Let N(G) = {0, 1, 2, I, 2I, 1 + I, 1 + 2I, 2 + I, 2 + 2I} be a neutrosophic group under multiplication modulo 3.

Consider P = {1, 2, I, 2I}. ⊂ N(G); P is a neutrosophic subgroup.

| | | |
|---|---|---|
| P. 0 | = | {0} |
| P. 1 | = | {1, 2, I, 2I} |
| | = | P2 |
| P. I | = | {I 2I} |
| | = | P. 2I |
| P (1 + I) | = | {1 + I, 2 + 2I, 2I, I} |
| P (2 + I) | = | {2 + I, 1 + 2I, 0} |

P (1 + 2I)   =   {1 + 2I, 2 + I, 0}
            =   P (2 + I)
P (2 + 2I)   =   {2 + 2I, 1 + I, I, 2I}
            =   P (1 + I).

We see the coset does not partition the neutrosophic group.

Now using the concept of pseudo neutrosophic subgroup we define pseudo coset.

**DEFINITION 1.1.18:** *Let N(G) be a neutrosophic group. K be a pseudo neutrosophic subgroup of N(G). Then for a ∈ N(G), Ka = {ka | k ∈ K} is called the pseudo right coset of K in N(G).*

On similar lines we define the notion of pseudo left coset of a pseudo neutrosophic subgroup K of N (G). We illustrate this by the following example:

***Example 1.1.16:*** Let N (G) = {0, 1, I, 2, 2I, 1 + I, 1 + 2I, 2 + I, 2 + 2I} be a neutrosophic group under multiplication modulo 3. Take K = {1, 1 + I}, a pseudo neutrosophic group.
    Now we will study the cosets of K . K . 0 = {0}.

K 1          =   {1, 1 + I}
K (1 + I)    =   {1 + I}
K 2          =   {2, 2 + 2I}
K. I         =   {I, 2I}
K 2I         =   {2I, I}
            =   K. I.
K (1 + 2I)   =   {1 + 2I}
K (2 + I)    =   {2 + I}
K (2 + 2I)   =   {2 +2I, 2}
            =   K.2.

We see even the pseudo neutrosophic subgroups do not in general partition the neutrosophic group which is evident from the example.

Now we proceed on to define the concept of center of a neutrosophic group.

**DEFINITION 1.1.19:** *Let N(G) be a neutrosophic group, the center of N(G) denoted by C(N(G)) = {x ∈ N(G) | ax = xa for all a ∈ N(G)}.*

***Note:*** Clearly C(N(G)) ≠ φ for the identity of the neutrosophic group belongs to C(N(G)). Also if N(G) is a commutative neutrosophic group then C(N(G)) = N(G). As in case of groups we can define in case of neutrosophic groups also direct product of neutrosophic groups N(G).

**DEFINITION 1.1.20:** *Let N(G₁), N(G₂), ..., N(Gₙ) be n neutrosophic groups the direct product of the n-neutrosophic groups is denoted by N(G) = N(G₁) × ... × N(Gₙ) = {(g₁, g₂, ..., gₙ) | gᵢ ∈ N(Gᵢ); i = 1, 2,..., n}.*

N(G) is a neutrosophic group for the binary operation defined is component wise; for if $*_1, *_2, ..., *_n$ are the binary operations on N(G₁), ..., N(Gₙ) respectively then for X = (x₁, ..., xₙ) and Y = (y₁, y₂, ..., yₙ) in N(G), X * Y = (x₁, ..., xₙ) * (y₁, ..., yₙ) = (x₁*y₁, ..., xₙ*yₙ) = (t₁, ..., tₙ) ∈ N(G) thus closure axiom is satisfied. We see if e = (e₁, ..., eₙ) is the identity element where each eᵢ is the identity element of N(Gᵢ); 1 ≤ i ≤ n then X * e = e * X = X.

It is left as a matter of routine for the reader to check N(G) is a neutrosophic group. Thus we see that the concept of direct product of neutrosophic group helps us in obtaining more and more neutrosophic groups.

*Note*: It is important and interesting to note that if we take in N(Gᵢ), 1 ≤ i ≤ n. some N(Gᵢ) to be just groups still we continue to obtain neutrosophic groups.

We now give some examples as illustrations.

***Example 1.1.17:*** Let N(G₁) = {0, 1, I, 1 + I} and G₂ = {g | g³ = 1}. N(G) = N(G₁) × G₂ = {(0, g) (0, 1) (0, g²) (1, g) (1, 1) (1, g²) (I, g) (1 g²) (I, 1) (1 + I, 1) (1 + I, g) ( 1 + I, g²)} is a

neutrosophic group of order 12. Clearly $\{(1, 1), (1, g), (1, g^2)\}$ is the group in N(G).

Now several other properties, which we have left out, can be defined appropriately.

***Note:*** We can also define independently the notion of pseudo neutrosophic group as a neutrosophic group, which has no proper subset, which is a group, but the pseudo neutrosophic group itself is a group. We can give examples of them, the main difference between a pseudo neutrosophic group and a neutrosophic group is that a pseudo neutrosophic group is a group but a neutrosophic group is not a group in general but only contains a proper subset, which is a group.

Now we give an example of a pseudo neutrosophic group.

***Example 1.1.18:*** Consider the set N(G) = {1, 1 + I} under the operation multiplication modulo 3. {1, 1 + I} is a group called the pseudo neutrosophic group for this is evident from the table.

| *     | 1     | 1 + I |
|-------|-------|-------|
| 1     | 1     | 1 + I |
| 1 + I | 1 + I | 1     |

Clearly {1, 1 + I} is group but has no proper subset which is a group. Also this pseudo neutrosophic group can be realized as a cyclic group of order 2.

For more about these concepts please refer [142-3].

## 1.2 Neutrosophic Semigroups

In this section we recall the notion of neutrosophic semigroups. The notion of neutrosophic subsemigroups, neutrosophic ideals, neutrosophic Lagrange semigroups etc. are introduced just for the sake of completeness. We illustrate them with examples and give some of its properties.

**DEFINITION 1.2.1:** *Let S be a semigroup, the semigroup generated by S and I i.e. S ∪ I denoted by ⟨S ∪ I⟩ is defined to be a neutrosophic semigroup.*

It is interesting to note that all neutrosophic semigroups contain a proper subset, which is a semigroup.

***Example 1.2.1:*** Let $Z_{12}$ = {0, 1, 2, …, 11} be a semigroup under multiplication modulo 12. Let N(S) = ⟨$Z_{12}$ ∪ I⟩ be the neutrosophic semigroup. Clearly $Z_{12}$ ⊂ ⟨$Z_{12}$ ∪ I⟩ and $Z_{12}$ is a semigroup under multiplication modulo 12.

***Example 1.2.2:*** Let Z = {the set of positive and negative integers with zero}, Z is only a semigroup under multiplication. Let N(S) = {⟨Z ∪ I⟩} be the neutrosophic semigroup under multiplication. Clearly Z ⊂ N(S) is a semigroup.

Now we proceed on to define the notion of the order of a neutrosophic semigroup.

**DEFINITION 1.2.2:** *Let N(S) be a neutrosophic semigroup. The number of distinct elements in N(S) is called the order of N(S), denoted by o(N(S)).*

If the number of elements in the neutrosophic semigroup N(S) is finite we call the neutrosophic semigroup to be finite otherwise infinite. The neutrosophic semigroup given in example 1.2.1 is finite where as the neutrosophic semigroup given in example 1.2.2 is of infinite order.

Now we proceed on to define the notion of neutrosophic subsemigroup of a neutrosophic semigroup N(S).

**DEFINITION 1.2.3:** *Let N(S) be a neutrosophic semigroup. A proper subset P of N(S) is said to be a neutrosophic subsemigroup, if P is a neutrosophic semigroup under the operations of N (S). A neutrosophic semigroup N(S) is said to* *have a subsemigroup if N(S) has a proper subset, which is a semigroup under the operations of N(S).*

It is interesting to note a neutrosophic semigroup may or may not have a neutrosophic subsemigroup but it will always have a subsemigroup.

Now we proceed on to illustrate these by the following examples:

***Example 1.2.3:*** Let $Z^+ \cup \{0\}$ denote the set of positive integers together with zero. $\{Z^+ \cup \{0\}, +\}$ is a semigroup under the binary operation '+'. Now let N(S) = $\langle Z^+ \cup \{0\}^+ \cup \{I\}\rangle$. N(S) is a neutrosophic semigroup under '+'. Consider $\langle 2Z^+ \cup I\rangle$ = P, P is a neutrosophic subsemigroup of N(S). Take R = $\langle 3Z^+ \cup I\rangle$; R is also a neutrosophic subsemigroup of N(S).

Now we have the following interesting theorem:

**THEOREM 1.2.1:** *Let N(S) be a neutrosophic semigroup. Suppose $P_1$ and $P_2$ be any two neutrosophic subsemigroups of N(S) then $P_1 \cup P_2$ (i.e. the union) the union of two neutrosophic subsemigroups in general need not be a neutrosophic subsemigroup.*

*Proof:* We prove this by using the following example. Let $Z^+$ be the set of positive integers; $Z^+$ under '+' is a semigroup.

Let N(S) = $\langle Z^+ \cup I\rangle$ be the neutrosophic semigroup under '+'. Take $P_1 = \{\langle 2Z \cup I\rangle\}$ and $P_2 = \{\langle 5Z \cup I \rangle\}$ to be any two neutrosophic subsemigroups of N(S). Consider $P_1 \cup P_2$ we see $P_1 \cup P_2$ is only a subset of N (S) for $P_1 \cup P_2$ is not closed under the binary operation '+'. For take $2 + 4I \in P_1$ and $5 + 5I \in P_2$. Clearly $(2 + 5) + (4I + 5I) = 7 + 9I \notin P_1 \cup P_2$. Hence the claim.

We proved the theorem 1.2.1 mainly to show that we can give a nice algebraic structure to $P_1 \cup P_2$ viz. neutrosophic bisemigroup.

Now we proceed on to define the notion of neutrosophic monoid.

**DEFINITION 1.2.4:** *A neutrosophic semigroup N(S) which has an element e in N(S) such that e \* s = s \* e = s for all s ∈ N(S), is called as a neutrosophic monoid.*

It is interesting to note that in general all neutrosophic semigroups need not be neutrosophic monoids.

We illustrate this by an example.

**Example 1.2.4:** Let N(S) = ⟨$Z^+$ ∪ I⟩ be the neutrosophic semigroup under '+'. Clearly N(S) contains no e such that s + e = e + s = s for all s ∈ N(S). So N(S) is just a neutrosophic semigroup and not a neutrosophic monoid.

Now we give an example of a neutrosophic monoid.

**Example 1.2.5:** Let N(S) = ⟨$Z^+$ ∪ I⟩ be a neutrosophic semigroup generated under '×'. Clearly 1 in N(S) is such that 1 × s = s for all s ∈ N(S). So N(S) is a neutrosophic monoid.

It is still interesting to note the following:
1. From the examples 1.2.4 and 1.2.5 we have taken the same set ⟨$Z^+$ ∪ I⟩ with respect the binary operation '+', ⟨$Z^+$ ∪ I⟩ is only a neutrosophic semigroup but ⟨ $Z^+$ ∪ I⟩ under the binary operation × is a neutrosophic monoid.
2. In general all neutrosophic monoids need not have its neutrosophic subsemigroups to be neutrosophic submonoids.

First to this end we define the notion of neutrosophic submonoid.

**DEFINITION 1.2.5:** *Let N(S) be a neutrosophic monoid under the binary operation \*. Suppose e is the identity in N(S), that is s \* e = e \* s = s for all s ∈ N(S). We call a proper subset P of N(S) to be a neutrosophic submonoid if*

1. *P is a neutrosophic semigroup under '*'.*
2. *e ∈ P, i.e., P is a monoid under '*'.*

**Example 1.2.6:** Let $N(S) = \langle Z \cup I \rangle$ be a neutrosophic semigroup under '+'. N(S) is a monoid. $P = \langle 2Z^{+} \cup I \rangle$ is just a neutrosophic subsemigroup whereas $T = \langle 2Z \cup I \rangle$ is a neutrosophic submonoid. Thus a neutrosophic monoid can have both neutrosophic subsemigroups, which are different from the neutrosophic submonoids.

Now we proceed on to define the notion of neutrosophic ideals of a neutrosophic semigroup.

**DEFINITION 1.2.6**: *Let N(S) be a neutrosophic semigroup under a binary operation *. P be a proper subset of N(S). P is said to be a neutrosophic ideal of N(S) if the following conditions are satisfied.*
1. *P is a neutrosophic semigroup.*
2. *for all p ∈ P and for all s ∈ N(S) we have p * s and s * p are in P.*

**Note:** One can as in case of semigroups define the notion of neutrosophic right ideal and neutrosophic left ideal. A neutrosophic ideal is one, which is both a neutrosophic right ideal, and a neutrosophic left ideal. In general a neutrosophic right ideal need not be a neutrosophic left ideal.

Now we proceed on to give example to illustrate these notions.

**Example 1.2.7:** Let $N(S) = \langle Z \cup I \rangle$ be the neutrosophic semigroup under multiplication.
  Take P to be a proper subset of N(S) where $P = \langle 2Z \cup I \rangle$. Clearly P is a neutrosophic ideal of N(S). Since N(S) is a commutative neutrosophic semigroup we have P to be a neutrosophic ideal.

*Note:* A neutrosophic semigroup N(S) under the binary operation * is said to be a neutrosophic commutative semigroup if a * b = b * a for all a, b ∈ N(S).

**Example 1.2.8:** Let N(S) = $\left\{ \begin{pmatrix} a & b \\ c & d \end{pmatrix} \middle/ a, b, c, d, \in \langle Z \cup I \rangle \right\}$ be a neutrosophic semigroup under matrix multiplication. Take

$$P = \left\{ \begin{pmatrix} x & y \\ 0 & 0 \end{pmatrix} \middle/ x, y \in \langle Z \cup I \rangle \right\}.$$

Clearly

$$\begin{pmatrix} a & b \\ c & d \end{pmatrix} \begin{pmatrix} x & y \\ 0 & 0 \end{pmatrix} \notin P \,;$$

but

$$\begin{pmatrix} x & y \\ 0 & 0 \end{pmatrix} \begin{pmatrix} a & b \\ c & d \end{pmatrix} \in P.$$

Thus P is only a neutrosophic right ideal and not a neutrosophic left ideal of N(S).

Now we proceed on to define the notion of neutrosophic maximal ideal and neutrosophic minimal ideal of a neutrosophic semigroup N(S).

**DEFINITION 1.2.7:** *Let N(S) be a neutrosophic semigroup. A neutrosophic ideal P of N(S) is said to be maximal if P ⊂ J ⊂ N(S), J a neutrosophic ideal then either J = P or J = N(S). A neutrosophic ideal M of N(S) is said to be minimal if φ ≠ T ⊆ M ⊆ N(S) then T = M or T = φ.*

We cannot always define the notion of neutrosophic cyclic semigroup but we can always define the notion of neutrosophic cyclic ideal of a neutrosophic semigroup N(S).

**DEFINITION 1.2.8:** *Let N(S) be a neutrosophic semigroup. P be a neutrosophic ideal of N (S), P is said to be a neutrosophic*

*cyclic ideal or neutrosophic principal ideal if P can be generated by a single element.*

We proceed on to define the notion of neutrosophic symmetric semigroup.

**DEFINITION 1.2.9:** *Let S(N) be the neutrosophic semigroup. If S(N) contains a subsemigroup isomorphic to S(n) i.e. the semigroup of all mappings of the set (1, 2, 3, ..., n) to itself under the composition of mappings, for a suitable n then we call S (N) the neutrosophic symmetric semigroup.*

**Remark:** We cannot demand the subsemigroup to be neutrosophic, it is only a subsemigroup.

**DEFINITION 1.2.10:** *Let N(S) be a neutrosophic semigroup. N(S) is said to be a neutrosophic idempotent semigroup if every element in N (S) is an idempotent.*

***Example 1.2.9:*** Consider the neutrosophic semigroup under multiplication modulo 2, where N (S) = {0, 1, I, 1 + I}. We see every element is an idempotent so N (S) is a neutrosophic idempotent semigroup.

Next we proceed on to define the notion of weakly neutrosophic idempotent semigroup.

**DEFINITION 1.2.11:** *Let N(S) be a neutrosophic semigroup. If N(S) has a proper subset P where P is a neutrosophic subsemigroup in which every element is an idempotent then we call P a neutrosophic idempotent subsemigroup.*

*If N(S) has at least one neutrosophic idempotent subsemigroup then we call N(S) a weakly neutrosophic idempotent semigroup.*

We illustrate this by the following example:

***Example 1.2.10:*** Let N(S) = {0, 2, 1, I, 2I, 1 + I, 2 + 2 I, 1 + 2I, 2 + I} be the neutrosophic semigroup under multiplication modulo 3. Take P = {1, I, 1 + 2I, 0}; P is a neutrosophic idempotent subsemigroup. So N(S) is only a weakly neutrosophic idempotent semigroup. Clearly N(S) is not a neutrosophic idempotent semigroup as $(1 + I)^2 = 1$ which is not an idempotent of N(S).

**DEFINITION 1.2.12:** *Let N(S) be a neutrosophic semigroup (monoid). An element $x \in N(S)$ is called an element of finite order if $x^m = e$ where e is the identity element in N(S) i.e. (se = es = s for all $s \in S$ ) (m the smallest such integer).*

**DEFINITION 1.2.13:** *Let N (S) be a neutrosophic semigroup (monoid) with zero. An element $0 \neq x \in N(S)$ of a neutrosophic semigroup is said to be a zero divisor if there exist $0 \neq y \in N(S)$ with x . y = 0. An element $x \in N(S)$ is said to be invertible if there exist $y \in N(S)$ such that xy = yx = e ($e \in N(S)$), is such that se = es = s for all $s \in N(S)$).*

***Example 1.2.11:*** Let N (S) = {0, 1, 2, I, 2I, 1 + I, 2 + I, 1 + 2I, 2 + 2I} be a neutrosophic semigroup under multiplication modulo 3. Clearly $(1 + I) \in N (S)$ is invertible for (1 + I) (1 + I) = 1 (mod 3). (2 + 2I) is invertible for $(2 + 2I)^2 = 1$ (mod 3). N(S) also has zero divisors for (2 + I) I = 2I + I = 0(mod 3). Also (2 + I ) 2 I = 0 (mod 3) is a zero divisor.

Thus this neutrosophic semigroup has idempotents, units and zero divisors. For more about neutrosophic semigroups please refer [142-3].

### 1.3 Neutrosophic Fields

In this section we just recall the definition of neutrosophic fields for they are used in the construction of neutrosophic rings like neutrosophic group rings, neutrosophic semigroup rings.

**DEFINITION 1.3.1:** *Let K be the field of reals. We call the field generated by $K \cup I$ to be the neutrosophic field for it involves the indeterminacy factor in it. We define $I^2 = I$, I + I = 2I i.e., I*

+...+ $I = nI$, and if $k \in K$ then $k.I = kI$, $0I = 0$. We denote the neutrosophic field by $K(I)$ which is generated by $K \cup I$ that is $K(I) = \langle K \cup I \rangle$. $\langle K \cup I \rangle$ denotes the field generated by $K$ and $I$.

*Example 1.3.1:* Let R be the field of reals. The neutrosophic field of reals is generated by R and I denoted by $\langle R \cup I \rangle$ i.e. R(I) clearly $R \subset \langle R \cup I \rangle$.

*Example 1.3.2:* Let Q be the field of rationals. The neutrosophic field of rationals is generated by Q and I denoted by Q(I).

**DEFINITION 1.3.2:** *Let K(I) be a neutrosophic field we say K(I) is a prime neutrosophic field if K(I) has no proper subfield, which is a neutrosophic field.*

*Example 1.3.3:* Q(I) is a prime neutrosophic field where as R(I) is not a prime neutrosophic field for $Q(I) \subset R(I)$.

Likewise we can define neutrosophic subfield.

**DEFINITION 1.3.3:** *Let K(I) be a neutrosophic field, $P \subset K(I)$ is a neutrosophic subfield of P if P itself is a neutrosophic field. K(I) will also be called as the extension neutrosophic field of the neutrosophic field P.*

We can also define neutrosophic fields of prime characteristic p (p is a prime).

**DEFINITION 1.3.4:** *Let $Z_p = \{0,1, 2, ..., p-1\}$ be the prime field of characteristic p. $\langle Z_p \cup I \rangle$ is defined to be the neutrosophic field of characteristic p. Infact $\langle Z_p \cup I \rangle$ is generated by $Z_p$ and I and $\langle Z_p \cup I \rangle$ is a prime neutrosophic field of characteristic p.*

*Example 1.3.4:* $Z_7 = \{0, 1, 2, 3, ..., 6\}$ be the prime field of characteristic 7. $\langle Z_7 \cup I \rangle = \{0, 1, 2, ..., 6, I, 2I, ..., 6I, 1 + I, 1 + 2I, ..., 6 + 6I \}$ is the prime field of characteristic 7.

Chapter Two

# NEUTROSOPHIC RINGS AND THEIR PROPERTIES

In this book the authors for the first time define the new notion of neutrosophic rings. This chapter has two sections. In the first section we introduce the notion of neutrosophic rings and their substructures like ideals and subrings. Section two introduces several types of neutrosophic rings.

## 2.1 Neutrosophic Rings and their Substructures

In this section we define the new notion called neutrosophic rings. Just as complex rings or complex fields include in them the notion of imaginary element i with $i^2 = -1$ or the complex number i, in the neutrosophic rings we include the indeterminate element I where $I^2 = I$. In most of the real world problems or in situations we see in general we cannot always predict the occurrence or non occurrence of an event, the occurrence or the non occurrence may be an indeterminate so we introduce the neutrosophic rings which can very easily cater to such situations. In many a cases we felt the concept of indeterminacy was more concrete than the notion of 'imaginary' so we have ventured to define neutrosophic algebraic concepts.

   Throughout this section $Z^+$ will denote the set of positive integers, $Z^-$ the set of negative integers, Z the set of positive and negative integers with zero. $Z^+ \cup \{0\}$ will show we have adjoined 0 with $Z^+$. Like wise the rationals Q, the reals R, the complex number C can have one part of it denoted by $Q^+$, $Q^-$, $R^+$, $R^-$ and so on.

Also $Z_n$ will denote the set of modulo integers i.e. {0, 1, 2, …, n − 1} i.e. n ≡ 0 (mod n).

Now we proceed on to define the neutrosophic ring.

**DEFINITION 2.1.1:** *Let R be any ring. The neutrosophic ring ⟨R ∪ I⟩ is also a ring generated by R and I under the operations of R.*

We first illustrate this by some examples.

***Example 2.1.1:*** Let Z be the ring of integers; ⟨Z ∪ I⟩ = {a + bI / a, b ∈ Z}. ⟨Z ∪ I⟩ is a ring called the neutrosophic ring of integers. Also Z $\underset{\neq}{\subseteq}$ ⟨Z ∪ I⟩.

***Example 2.1.2:*** Let Q be the ring of rationals. ⟨Q ∪ I⟩ = {r + sI | r, s ∈ Q} is the ring called the neutrosophic ring of rationals.

It is important to note that though Q is the field, ⟨Q ∪ I⟩ neutrosophic ring of rationals is not a field it is only a ring for $I^2$ = I and I has no inverse, yet we call it the neutrosophic field of rationals.

***Example 2.1.3:*** Let $\mathfrak{R}$ be the ring of reals, ⟨$\mathfrak{R}$ ∪ I⟩ is the neutrosophic ring called the neutrosophic ring of reals. Here also ⟨$\mathfrak{R}$ ∪ I⟩ is only a ring and not a field yet we call it the neutrosophic field of reals.

***Example 2.1.4:*** Let C be the complex field, ⟨C ∪ I⟩ is the neutrosophic ring of complex numbers. ⟨C ∪ I⟩ is not a field. ⟨C ∪ I⟩ is called the complex neutrosophic ring, but we call it the neutrosophic field of complex numbers.

***Example 2.1.5:*** Let $Z_n = \{0, 1, 2, \ldots, n-1\}$ be the ring of integers modulo n. $\langle Z_n \cup I \rangle$ is the neutrosophic ring of modulo integers n.

***Example 2.1.6:*** Let $Z_2 = \{0, 1\}$ be the ring of integers modulo 2. $\langle Z_2 \cup I \rangle = \{0, 1, I, 1+I\}$ is the neutrosophic ring. Clearly $\langle Z_2 \cup I \rangle$ is not the field of prime characteristic two. It is only a neutrosophic field of characteristic 2.

This notion will be dealt in this book in due course of time.

**DEFINITION 2.1.2:** *Let $\langle R \cup I \rangle$ be a neutrosophic ring. We say $\langle R \cup I \rangle$ is a neutrosophic ring of characteristic zero if nx = 0 (n an integer) for all $x \in \langle R \cup I \rangle$ is possible only if n = 0, then we call the neutrosophic ring to be a neutrosophic ring of characteristic zero.*

*If in the neutrosophic ring $\langle R \cup I \rangle$, nx = 0 (n an integer) for all $x \in \langle R \cup I \rangle$, then we say the neutrosophic ring $\langle R \cup I \rangle$ is of characteristic n and n is called the characteristic of the neutrosophic ring $\langle R \cup I \rangle$.*

It is pertinent to mention n can be a prime or a composite number. We will first illustrate these by the following examples:

***Example 2.1.7:*** Let $\langle Q \cup I \rangle$ be the neutrosophic ring of rationals. $\langle Q \cup I \rangle$ is the neutrosophic ring of characteristic zero.

***Example 2.1.8:*** Consider $\langle C \cup I \rangle$ the neutrosophic ring of complex numbers. $\langle C \cup I \rangle$ is the neutrosophic ring of characteristic zero.

***Note:*** We can use the term neutrosophic ring for the neutrosophic field also as every field is a ring.

***Example 2.1.9:*** Consider the neutrosophic ring $\langle Z_4 \cup I \rangle = \{0, 1, 2, 3, I, 2I, 3I, 1+I, 2+I, 3+I, 2I+1, 2I+2, 2I+3, 3I+1, 3I+2, 3I+3\}$. $\langle Z_4 \cup I \rangle$ is a neutrosophic ring of characteristic 4. $(2I+2)(2I+2) = 4I^2 + 4 + 8I = 4$ $([I+1+2I]) \equiv 0 \pmod 4$.

We call all neutrosophic rings, which are of finite characteristic as finite characteristic neutrosophic rings.

*Example 2.1.10:* Let $\langle Z_5 \cup I \rangle$ = {0, 1, 2, 3, 4, I, 2I, 3I, 4I, 1 + I, 2 + I, 3 + I, 4 + I, 2 + 2I, 3 + 2I, 4 + 2I, 3 + 3I, 3 + 4I, 4 + 2I, 4 + 4I, 3I + 2, 3I + 1, 2I + 1, 4I + 2, 4I + 1}, $\langle Z_5 \cup I \rangle$ is a neutrosophic ring of characteristic five; or we say $\langle Z_5 \cup I \rangle$ is a neutrosophic ring of prime characteristic. All elements in $\langle Z_5 \cup I \rangle$ are such that 5 x $\equiv$ 0 (mod 5) for all x $\in \langle Z_5 \cup I \rangle$.

Having defined the characteristic of a neutrosophic ring we proceed on to prove every neutrosophic ring is a ring and every neutrosophic ring contains a proper subset, which is just a ring.

**THEOREM 2.1.1:** *Let $\langle R \cup I \rangle$ be a neutrosophic ring. $\langle R \cup I \rangle$ is a ring.*

*Proof:* Let $\langle R \cup I \rangle$ be a neutrosophic ring. R is a ring in which '+' and '•' are the binary operations. To show $\langle R \cup I \rangle$ is a ring we have to prove $\langle R \cup I \rangle$ is an abelian group under '+' and a semigroup under '•' $\langle R \cup I \rangle$ = {a + bI / a, b $\in$ R}. Let a + bI, c + dI $\in \langle R \cup I \rangle$. To show closure it is enough if we show (a + bI) + (c + dI) is in $\langle R \cup I \rangle$.

Consider (a + bI) + (c + dI) = (a + c) + (b + d) I; as a + c and b + d $\in$ R we see (a + bI) + (c + dI) is in $\langle R \cup I \rangle$. Thus closure axiom under the operation + is satisfied. It is easily verified (a + bI) + (c + dI) = (c + dI) + (a + bI) as a + c = c + a and b + d = d + b in R. Thus the binary operation '+' is commutative on $\langle R \cup I \rangle$.

Since R is ring '0' is its additive identity. Take (a + bI + 0) = a + 0 + (b + 0) I = a + bI. So '0' is the additive identity of $\langle R \cup I \rangle$. For every a + bI $\in \langle R \cup I \rangle$ we have a unique – a + (– b)I in $\langle R \cup I \rangle$ which is such that (a + bI) + (–a + (–b) I) = 0. This follows from the fact for every a $\neq$ 0 in R we have a unique –a in R with a + (–a) = 0. Thus $\langle R \cup I \rangle$ under '+' is an abelian group.

To show ⟨R ∪ I⟩ is ring it is sufficient if we show ⟨R ∪ I⟩ is closed with respect to '•' for other properties follow easily. Take a + bI, c + dI ∈ ⟨R ∪ I⟩ consider (a + bI) • (c + dI) = ac + bcI + adI + bdI$^2$ = ac + (bc + ad + bd) I = x + yI ∈ ⟨R ∪ I⟩ where x = ac and y = bc + ad + bd as x and y are in R. Hence the claim. Thus ⟨R ∪ I⟩ is a ring under '+' and '•'.

*Note:* We have proved this mainly to show in general a neutrosophic ring is a ring; but a neutrosophic group may not have a group structure. This is a vast difference between these two algebraic structures.

**DEFINITION 2.1.3:** *Let ⟨R ∪ I⟩ be a neutrosophic ring. We say ⟨R ∪ I⟩ is a commutative neutrosophic ring if for all x, y ∈ ⟨R ∪ I⟩ we must have x y = y x. If even for a single pair xy ≠ y x we call the neutrosophic ring to be a non commutative neutrosophic ring.*

We illustrate them by the following example:

**Example 2.1.11:** Let ⟨Q ∪ I⟩ be a neutrosophic ring. ⟨Q ∪ I⟩ is a commutative ring of characteristic zero.

**Example 2.1.12:** Let ⟨R ∪ I⟩ = $\left\{ \begin{bmatrix} a & b \\ c & d \end{bmatrix} \middle/ \ a, b, c, d \in \langle Q \cup I \rangle \right\}$

be the set of all 2 × 2 matrices. ⟨R ∪ I⟩ is a non commutative ring of characteristic zero under matrix addition and matrix multiplication. We say ⟨R ∪ I⟩ is a neutrosophic ring with unit, if 1 ∈ ⟨R ∪ I⟩ is such that 1 x = x. 1 = x for all ∈ ⟨R ∪ I⟩. All the neutrosophic rings given as examples are neutrosophic rings with unit.

Now we go for examples of non-commutative neutrosophic rings of finite order.

**Example 2.1.13:** Let ⟨R ∪ I⟩ = $\left\{ \begin{pmatrix} a & b \\ c & d \end{pmatrix} \middle| a, b, c, d, \in \{0, 1, 2, \right.$

I, 2I, 1 + I, 2 + I, 1 + 2I, 2 + 2I}. ⟨R ∪ I⟩ is a neutrosophic ring

under usual matrix addition and multiplication and ⟨R ∪ I⟩ is a non commutative neutrosophic ring of finite order.

We will define subrings and ideals of a neutrosophic ring.

**DEFINITION 2.1.4:** *Let ⟨R ∪ I⟩ be a neutrosophic ring. A proper subset P of ⟨R ∪ I⟩ is said to be a neutrosophic subring if P itself is a neutrosophic ring under the operations of ⟨R ∪ I⟩. It is essential that P = ⟨S ∪ nI⟩, n a positive integer where S is a subring of R. i.e. {P is generated by the subring S together with n I. (n ∈ Z⁺)}.*

*Note:* Even if P is a ring and cannot be represented as ⟨S ∪ nI⟩ where S is a subring of R then we do not call P a neutrosophic subring of ⟨R ∪ I⟩.

**Example 2.1.14:** Let ⟨$Z_{12}$ ∪ I⟩ be a finite neutrosophic ring of characteristic 12. P = {0, 6, I, 6I, 2I, 3I, 4I, 5I, 7I, 8I, …, 11I, 6 + I, 6 + 2I, …, 6 + 11I} is a neutrosophic subring of ⟨$Z_{12}$ ∪ I⟩ as S = {0, 6} is a subring of $Z_{12}$.

So P = ⟨S ∪ I⟩ is a neutrosophic subring of ⟨$Z_{12}$ ∪ I⟩. Take $P_1$ = {0, 2, 4, 6, 8, 10, 2I, 4I, 6I, 8I, 10I, 2 + 2I, 2 + 4I, …, 2 + 10I, 10 + 2I, …, 10 + 10I} is a neutrosophic subring of ⟨$Z_{12}$ ∪ I⟩. Take $P_2$ = {0, 6 + 6I} this is a subring under the operations of ⟨$Z_{12}$ ∪ I⟩ but is not a neutrosophic subring.

How to characterize such neutrosophic rings which are not generated by a ring? To this we define a notion called pseudo neutrosophic rings.

**DEFINITION 2.1.5:** *Let T be a non empty set with two binary operations '+' and '×'. T is said to be a pseudo neutrosophic ring if the following conditions are true:*

1. *T contains elements of the form x + yI (x, y are reals y ≠ 0 for atleast one value).*
2. *(T, +) is an abelian group.*
3. *(T, ×) is a semigroup.*

4.    $t (t_1 + t_2) = tt_1 + tt_2$ and $(t_1 + t_2) t = tt_1 + t_2t$ for all $t, t_1, t_2 \in T$.

*We do not see any visible relations between a neutrosophic ring and a pseudo neutrosophic ring.*

**Example 2.1.15:** T = {0, 2I, 4I, 6I, 8I, 10I} is a pseudo neutrosophic ring under addition modulo 12 and usual multiplication modulo 12. Clearly T is not a neutrosophic ring.

**Example 2.1.16:** Take T = {2ZI ∪ 0} = {0, ±2I, ±4I, ±6 I, ±8I, …, ∞}, T is a pseudo neutrosophic ring which is not a neutrosophic ring.

One can construct several such examples.

Next we proceed on to define the notion of ideals in a neutrosophic rings.

**DEFINITION 2.1.6:** *Let ⟨R ∪I⟩ be any neutrosophic ring, a non empty subset P of ⟨R ∪I⟩ is defined to be a neutrosophic ideal of ⟨R ∪I⟩ if the following conditions are satisfied;*

1.   *P is a neutrosophic subring of ⟨R ∪I⟩.*
2.   *For every p ∈ P and r ∈ ⟨R ∪I⟩, rp and pr ∈ P.*

**Note:** The notion of neutrosophic right ideal and neutrosophic left ideal can also be defined provided ⟨R ∪ I⟩ is a non commutative neutrosophic ring.

**Example 2.1.17:** Let ⟨Z_{12} ∪ I⟩ be a neutrosophic ring under '+' and × modulo 12. Let {0, 6, 2I, 4I, 6I, 8I, 10I, 6 + 2I, …, 6 + 10I} be the neutrosophic subring of ⟨Z_{12} ∪ I⟩. It is easily verified ⟨Z_{12} ∪I⟩ is a neutrosophic ideal of ⟨Z_{12} ∪I⟩.

Now we can define the notion of neutrosophic pseudo ideal of the neutrosophic ring ⟨R ∪ I⟩.

**DEFINITION 2.1.7:** *Let $\langle R \cup I \rangle$ be a neutrosophic ring. P a pseudo neutrosophic subring of R. If for all $p \in P$ and $r \in \langle R \cup I \rangle$ we have pr, r p $\in$ P. Then we call P to be a pseudo neutrosophic ideal of $\langle R \cup I \rangle$.*

We illustrate this by the following example 2.1.18.

**Example 2.1.18:** Let $\langle Z_{16} \cup I \rangle$ be a neutrosophic ring of characteristic 16. Take P = {0, 2I, 4I, 6I, 8I, 10I, 12I, 14I} in $\langle Z_{16} \cup I \rangle$. P is pseudo neutrosophic subring. In fact P is a pseudo neutrosophic ideal of $\langle Z_{16} \cup I \rangle$. o$\langle Z_{16} \cup I \rangle = 16^2$.

Now we make the following observations.

*Observation 1:* Let $\langle Z_n \cup I \rangle$ be a neutrosophic ring $n < \infty$, o$(\langle Z_n \cup I \rangle)$ = $n^2$.

*Observation 2:* Does the order of the neutrosophic ideal or the pseudo neutrosophic ideal divide the order of the neutrosophic ring $\langle Z_n \cup I \rangle$; $n < \infty$.

We just illustrate by some more examples before we propose some open problems.

**Example 2.1.19:** Let $\langle Z_5 \cup I \rangle$ = {0, 1, 2, 3, 4, I, 2I, 3I, 4I, 1 + I, 2 + I,…, 4 + 3I, 4 + 4I} be a finite neutrosophic ring. Clearly o$(\langle Z_5 \cup I \rangle)$ = 25. Take {0, I, 2I, 3I, 4I} = P; P is a pseudo neutrosophic ideal of $\langle Z_5 \cup I \rangle$. $\langle Z_5 \cup I \rangle$ has no neutrosophic ideals and has only one pseudo neutrosophic ideal of order 5.

**THEOREM 2.1.2:** *Let $\langle Z_p \cup I \rangle$ be a neutrosophic ring, where p is a prime. o $(\langle Z_p \cup I \rangle)$ = $p^2$.*

1. *$\langle Z_p \cup I \rangle$ has no neutrosophic ideals and*
2. *$\langle Z_p \cup I \rangle$ has only one pseudo neutrosophic ideal of order p.*

*Proof:* Given $\langle Z_p \cup I \rangle$ is a neutrosophic ring of order $p^2$, p a prime. To prove $\langle Z_p \cup I \rangle$ has no neutrosophic ideals. On the contrary, suppose $\langle Z_p \cup I \rangle$ has a neutrosophic ideal say P; then P = $\langle P_1 \cup I \rangle$, $P_1$ a proper subring of $Z_p$. Since p is a prime, $Z_p$ has no subrings other than {0} and $Z_p$. Hence P = $\langle P_1 \cup I \rangle$ is impossible for any $P_1$ as a subring. Thus $\langle Z_p \cup I \rangle$ has no neutrosophic ideals.

To show $\langle Z_p \cup I \rangle$ has only one pseudo neutrosophic ideal.

Take K = {0, I, 2I, 3I, …, (p − 1)I}; K is a pseudo neutrosophic ideal of $\langle Z_p \cup I \rangle$. o(K) = p. Since p is a prime every element of K \ {0} generates K so K is maximal and minimal. Hence the claim.

Now we proceed on to give examples of our observations.

***Example 2.1.20:*** Let $\langle Z_4 \cup I \rangle$ be a neutrosophic ring of characteristic 4. $\langle Z_4 \cup I \rangle$ = {0, 1, 2, 3, I, 2I, 3I, 1 + I, 1 + 2I, 1 + 3I, 2 + I, 2 + 2I, 2 + 3I, 3 + I, 3 + 2I, 3 + 3I}. ($\langle Z_4 \cup I \rangle$ has both neutrosophic ideals and pseudo neutrosophic ideals.

Take P = {0, 2, 2I, 2 + 2I} ⊂ $\langle Z_4 \cup I \rangle$, P is a neutrosophic ideal of $\langle Z_4 \cup I \rangle$.

Consider L = {0, I, 2I, 3I} ⊂ $\langle Z_4 \cup I \rangle$, L is a pseudo neutrosophic ideal of $\langle Z_4 \cup I \rangle$. Thus $\langle Z_4 \cup I \rangle$ has both neutrosophic ideals and pseudo neutrosophic ideals.

***Example 2.1.21:*** Let $\langle Z_6 \cup I \rangle$ = {0, 1, 2, 3, 4, 5, I, 2I, 3I, 4I, 5I, 1 + I, 1 + 2I, 1 + 3I, …, 5 + I, 5 + 2I, 5 + 3I, 5 + 4I, 5 + 5I} be a neutrosophic ring o ($\langle Z_6 \cup I \rangle$) = $6^2$.

Take P = {0, 3, 3I, 3 + 3I}, P is a neutrosophic ideal of $\langle Z_6 \cup I \rangle$.

Let S = {0, I, 2I, 3I, 4I, 5I}, S is a pseudo neutrosophic ideal of $\langle Z_6 \cup I \rangle$. T = {0, 2, 4, 2I, 4I, 4 + 2I, 2 + 2I, 2 + 4I, 4I + 4} ⊂ $\langle Z_6 \cup I \rangle$; T is a neutrosophic ideal of $\langle Z_6 \cup I \rangle$. o(T) = 9 and o(P) = 3. L = {0, 3, I, 2I, 3I, 4I, 5I 3 + I, 3 + 2I, 3 + 3I, 3 + 4I, 3 + 5I} is a neutrosophic ideal of $\langle Z_6 \cup I \rangle$, o(L) = 12.

Now we define the quotient neutrosophic ring and pseudo quotient neutrosophic ring.

**DEFINITION 2.1.8** *Let $\langle R \cup I \rangle$ be a neutrosophic ring, P be a neutrosophic ideal of $\langle R \cup I \rangle$. The quotient neutrosophic ring*

$$\frac{\langle R \cup I \rangle}{P} = \{\, a + P \mid a \in \langle R \cup I \rangle \,\}$$

*provided $\dfrac{\langle R \cup I \rangle}{P}$ is a neutrosophic ring. It is easily verified that $\dfrac{\langle R \cup I \rangle}{P}$ is a neutrosophic ring with P acting as the additive identity.*

We illustrate this by the following example:

***Example 2.1.22:*** Let $\langle Z \cup I \rangle$ be a neutrosophic ring. P = $\langle 3Z \cup I \rangle$, the neutrosophic ideal of $\langle Z \cup I \rangle$.

The quotient neutrosophic ring

$$\frac{\langle Z \cup I \rangle}{P} = \{a + P \mid a \in \langle Z \cup I \rangle\}$$

= {P, 1 + P, 2 + P, I + P, 2I + P, (1 + I) + P, (2 + I) + P, (1 + 2I) + P, (2 + 2I) + P} is a neutrosophic ring.

We give another example.

***Example 2.1.23:*** Let $\langle Z_{12} \cup I \rangle$ = {0, 1, 2, 3, …, 11, I, 2I, …, 11I, 1 + I, 2 + I, …, 11 + I, …, 11 + 11I} be a neutrosophic ring.

Take P = {0, 6, I, 2I, …, 11I, 1 + 6, …, 11I + 6} be the neutrosophic ideal of $\langle Z_{12} \cup I \rangle$.

$$\frac{\langle Z_{12} \cup I \rangle}{P} = \{P, 1 + P, 2 + P, …, 5 + P\}, \text{ the quotient}$$

is a ring but is not a neutrosophic ring.

We make a special definition in this case.

**DEFINITION 2.1.9:** *Let $\langle R \cup I \rangle$ be a neutrosophic ring, P be a neutrosophic ideal of $\langle R \cup I \rangle$.*

*If the quotient $\dfrac{\langle R \cup I \rangle}{P}$ is just a ring and not a neutrosophic ring then we call $\dfrac{\langle R \cup I \rangle}{P}$ to be a false neutrosophic quotient ring.*

The neutrosophic ring given in example 2.1.23 is only a false neutrosophic quotient ring.

**Example 2.1.24:** Let $\langle Z_{12} \cup I \rangle$ be a neutrosophic ring, J = {0, 2I, 4I, 6I, 8I, 10I} be the pseudo neutrosophic ideal of $\langle Z_{12} \cup I \rangle$.

The quotient ring $\dfrac{\langle Z_{12} \cup I \rangle}{J}$ = {J, 1 + J, 2 + J, ..., 11 + J, I J, (1 + I) + J, ..., (11 + I) + J} is a neutrosophic ring which we call differently.

**DEFINITION 2.1.10:** *Let $\langle R \cup I \rangle$ be a neutrosophic ring. Suppose J be a pseudo neutrosophic ideal of $\langle R \cup I \rangle$.*

*If the quotient ring $\dfrac{\langle R \cup I \rangle}{J}$, is a neutrosophic ring then we call it as the pseudo quotient neutrosophic ring. If $\dfrac{\langle R \cup I \rangle}{J}$ is just a ring we call it a false pseudo quotient neutrosophic ring.*

We just give an example of the same.

**Example 2.1.25:** Let $\langle Z_6 \cup I \rangle$ be a neutrosophic ring. Let P = {0, I, 2I, 3I, 4I, 5I} be the pseudo neutrosophic ring.

The quotient ring

$$\frac{\langle Z_6 \cup I \rangle}{P} = \{P, 1 + P, 2 + P, 3 + P, 4 + P, 5 + P\}$$

is a false pseudo neutrosophic ring.

Now we just define one small notion before we proceed to define some more quotient.

**DEFINITION 2.1.11:** *Let $\langle R \cup I \rangle$ be a neutrosophic ring, P a proper subset of $\langle R \cup I \rangle$ which is just a ring. Then we call P just a subring.*

It is important to note that we cannot have ideals in $\langle R \cup I \rangle$ which are not neutrosophic. But we have quotient ring which may not be a neutrosophic ring.

We just illustrate by an example.

**Example 2.1.26:** Let $\langle Z_7 \cup I \rangle$ be a neutrosophic ring.
    Take P = {0, I, 2I, 3I, …, 6I};

$$\frac{\langle Z_7 \cup I \rangle}{P} = \{P, 1 + P, 2 + P, 3\,P, …, 6 + P\}$$

is a false neutrosophic ring isomorphic to $Z_7$.

In view of this we have the following theorem:

**THEOREM 2.1.3:** *Let $\langle Z_n \cup I \rangle$ be a neutrosophic ring. $\langle Z_n \cup I \rangle$ has a pseudo ideal P such that $\frac{\langle Z_n \cup I \rangle}{P} \cong Z_n$.*

*Proof:* Given $\langle Z_n \cup I \rangle$ is a neutrosophic ring. Choose P = {0, I, 2I, …, (n − 2) I, (n − 1)I}, Clearly P is a pseudo neutrosophic ideal of $\langle Z_n \cup I \rangle$ and o(P) = n..
    Now the pseudo quotient neutrosophic ring

$$\frac{\langle Z_n \cup I \rangle}{P} = \{P, 1 + P, 2 + P, …, (n − 2) + P, (n − 1) + P\} \cong Z_n.$$

Hence the claim. Further $\frac{\langle Z_n \cup I \rangle}{P}$ is only a false neutrosophic quotient ring.

**Example 2.1.27:** Let $\langle Z_{10} \cup I \rangle$ be a neutrosophic ring. Take P = {0, 5, 5I, 5 + 5I} $\subset \langle Z_{10} \cup I \rangle$; P is an ideal of $\langle Z_{10} \cup I \rangle$.

$\dfrac{\langle Z_{10} \cup I \rangle}{P}$ = {P, 1 + P, 2 + P, 3 + P, 4 + P, I + P, 2I + P, 3I + P, 4I + P, (1 + 4I) + P, + P (1 + I) + P (1 + 3I) + P, (1 + 4I) + P, 2 + I + P, 2 + 2I + P, 2 + 3I + P, 2 + 4I + P, 3 + I + P, 3 + 2I + P, 3 + 3I + P, 3 + 4I + P, 4 + I + P, 4 + 2I + P, 4 + 3I + P, 4 + 4I + P}, is the neutrosophic quotient ring, o $\left[ \dfrac{\langle Z \cup I \rangle}{P} \right]$ = 25.

Now we proceed on to define homomorphism of two neutrosophic rings.

**DEFINITION 2.1.12:** *Let $\langle R \cup I \rangle$ and $\langle S \cup I \rangle$ be any two neutrosophic rings. A map $\phi$ from $\langle R \cup I \rangle$ to $\langle S \cup I \rangle$ is said to be a neutrosophic ring homomorphism if*

1. *$\phi$ is a ring homomorphism.*
2. *$\phi(I) = I$ for I the neutrosophic element of $\langle R \cup I \rangle$ i.e. the indeterminate "I' remains as it is.*

*We can have several examples of such homomorphisms. (Clearly Ker $\phi = \{x \in \langle R \cup I \rangle \mid \phi(x) = 0\}$ is always a subring which is never a neutrosophic subring. Ker $\phi$ can never be an ideal in a neutrosophic ring.*

Now we proceed on to define the notion of some special type of neutrosophic rings.

## 2.2 Special Type of Neutrosophic Rings

In this section we introduce several types of neutrosophic rings and study them. Some of them are commutative some non-commutative, some finite, some finite characteristic but infinite and some are zero characteristic and infinite. We introduce neutrosophic polynomial ring, neutrosophic matrix ring, neutrosophic direct product rings, neutrosophic integral domains, neutrosophic unique factorization domains, neutrosophic division rings, neutrosophic integral quaternions, neutrosophic ring of real quaternions.

**DEFINITION 2.2.1:** *Let $\langle R \cup I \rangle$ be a commutative ring with unity 1. (i.e. 1 r = r 1 = r for all $r \in \langle R \cup I \rangle$). x be variate $\langle R \cup I \rangle$ [x], the neutrosophic ring of polynomials in the variable x consists of the set of symbols.*

$$a_0 + a_1 x + \dots + a_n x^n$$

*where n can be any non negative integer and where the coefficients $a_0, a_1, \dots, a_n$ are all in $\langle R \cup I \rangle$. We call $a_0 + a_1 x + \dots + a_n x^n$ a neutrosophic polynomial in the variate x. In order to make ($\langle R \cup I \rangle$ [x]) the neutrosophic ring we must be able to recognize when two elements in it are equal, we must be able to add and multiply elements of ($\langle R \cup I \rangle$ [x]). As in case of polynomial rings we in case of neutrosophic polynomial rings also we say when two polynomials are equal. Let*

$$p\ (x) = a_0 + a_1 x + \dots + a_m x^m \text{ and}$$
$$q\ (x) = b_0 + b_1 x + \dots + b_n x^n$$

*be in ($\langle R \cup I \rangle$ [x]), then p(x) = q(x) if and only if for every integer $i \geq 0$, $a_i = b_i$, $a_0, \dots, a_m, b_0, b_1, \dots, b_n$ are from $\langle R \cup I \rangle$ and m = n.*

*Addition of two neutrosophic polynomials is done as in case of usual polynomials; i.e.*

$$p(x) + q(x) = C_0 + C_1 + \dots + C_t x^t$$

*where for each i, $C_i = a_i + b_i$. In other words add two neutrosophic polynomials by adding their coefficients and collecting terms.*

*Suppose we want to add $7 + I + (2 - 5I)\ x$ and $8 + 3I - (3 + 12I)\ x + (5 - I)\ x^2$. We consider $(7 + I) + (2 - 5I)\ x + (0 + 0I)\ x^2 + (8 + 3I) - (3 + 12I)\ x + (5 - I)\ x^2 = 15 + 4I + (-1 + 7I)\ x + (5 - I)x^2$.*

Now how to define product i.e. multiplication of any two neutrosophic polynomials in ($\langle R \cup I \rangle$ [x]). If

$$p(x) = a_0 + a_1 x + \ldots + a_m x^m \text{ and}$$
$$q(x) = b_0 + b_1 x + \ldots + b_n x^n,$$

where $p(x)$, $q(x) \in (\langle R \cup I \rangle [x])$ then
$$p(x) \, q(x) = C_0 + C_1 x + \ldots + C_k x^k, \; C_t = a_t b_0 + \ldots + a_0 b_t$$
here we use $I^2 = I$. Let

| | | |
|---|---|---|
| $p(x)$ | $=$ | $(2 - I) + (7 + 5I) x - (8I) x^2$ |
| $q(x)$ | $=$ | $(7 + I) + (2 + 5I) x^2 + (3 + I) x^3$ |
| $a_0$ | $=$ | $(2 - I),$ |
| $a_1$ | $=$ | $(7 + 5I),$ |
| $a_2$ | $=$ | $-8I,$ |
| $a_3 = a_4$ | $= \ldots = 0$ and | |
| $b_0$ | $=$ | $7 + I$ |
| $b_1$ | $=$ | $0,$ |
| $b_2$ | $=$ | $(2 + 5I)$ |
| $b_3$ | $=$ | $3 + I$ |
| $b_4 = b_5$ | $= \ldots = 0.$ | |

$$
\begin{aligned}
C_0 &= a_0 b_0 \\
&= (2 - I)(7 + I) \\
&= 14 - 7I + 2I - I^2 \\
&= 14 - 5I - I \; (\because) \; I^2 = I \\
&= 14 - 6I.
\end{aligned}
$$

$$
\begin{aligned}
C_1 &= a_1 b_0 + a_0 b_1 \\
&= (7 + 5I)(7 + I) + (2 - I).0 \\
&= 49 + 5I^2 + 35I + 7I + 0 \\
&= 49 + 47I.
\end{aligned}
$$

$$
\begin{aligned}
C_2 &= a_2 b_0 + a_1 b_1 + a_0 b_2 \\
&= -8I(7 + I) + (7 + 5I)0 + (2 - I)(2 + 5I) \\
&= 4 - 68I.
\end{aligned}
$$

$$
\begin{aligned}
C_3 &= a_3 b_0 + a_2 b_1 + a_1 b_2 + a_0 b_3 \\
&= 0 - 8I(0) + (7 + 5I)(2 + 5I) + (2 - I)(3 + I) \\
&= 0 - 0 + 14 + 10I + 25I^2 + 35I + 6 - I^2 - 3I + 2I \\
&= 20 - 69I.
\end{aligned}
$$

$$C_4 \quad = \quad a_4b_0 + a_3b_1 + a_2b_2 + a_1\,b_3 + a_0b_4$$
$$= \quad 0b_0 + 0b_1 + (-8I)(2+5I) + (7+5I)(3+I) + (2-I)\,0$$
$$= \quad -16I - 40I^2 + 21 + 5I^2 + 15I + 7I + 0$$
$$= \quad 21 - 29I.$$

$$C_5 \quad = \quad a_5b_0 + a_4b_1 + a_3b_2 + a_2b_3 + a_1\,b_4 + a_0\,b_5$$
$$= \quad 0b_0 + 0b_0 + 0b_2 + (-8I)\,(3+I) + a_10 + a_0.0.$$
$$= \quad -24I^2 - 8I\ (\because\ I^2 = I).$$

$$C_5 \quad = \quad -32I.$$

$$C_6 \quad = \quad a_6b_0 + a_5b_1 + a_4b_2 + a_3b_3 + a_2b_4 + a_1b_5 + a_0\,b_6$$
$$= \quad 0b_0 + 0b_1 + 0b_2 + 0b_3 + (-8I)\,0 + (7+5I).0 + a_0\,b_6.0.$$
$$= \quad 0.$$

$$p(x)\,q(x) \quad = \quad (14 - 6I) + (49 + 47I)x + (4 - 61I)x^2 + (20 - 68I)\,x^3 + (21 - 29I)\,x^4 - 32I\,x^5.$$

$(\langle R \cup I\rangle[x])$ is a commutative neutrosophic ring with unit element.

It is left as an exercise for the reader to prove $(\langle R \cup I\rangle\,[x])$ is a commutative ring.

**THEOREM 2.2.1:** *Every neutrosophic polynomial ring contains a polynomial ring.*

*Proof:* Given $(\langle R \cup I\rangle\,[x])$ is the polynomial neutrosophic ring where R is a commutative ring with unit; $\langle R \cup I\rangle$ is the neutrosophic ring generated by R and I; which is again a commutative ring with unity. Now $R[x] \subseteq (\langle R \cup I\rangle)[x]$, hence the claim.

Now we wish to make a few remarks as definitions.

**DEFINITION 2.2.2:** *Let $(\langle R \cup I\rangle\,[x])$ be a neutrosophic polynomial ring. A polynomial $p(x) = a_o + a_1\,x + \ldots + a_n\,x$ is called a strong neutrosophic polynomial if each $a_i$ is of the form $a + bI$ where $b \neq 0$, $a, b \in R$ the polynomial $q\,(x) = b_0 + b_1x +$*

*... + $b_m x^m$ is called the mixed neutrosophic polynomial if each $b_i$ ∈ $\langle R \cup I \rangle$, i.e. some $b_i$ are from R some $b_j$'s are of the form c + dI (d ≠ 0). A polynomial t (x) = $t_0 + t_1 x + ... + t_r x^r$ is said to be a polynomial if each $t_i$ ∈ R.*

Now we give some examples of neutrosophic polynomial rings before we proceed onto define substructures in them.

***Example 2.2.1:*** Let $\langle Z \cup I \rangle$ be the neutrosophic ring of integers, $[\langle Z \cup I \rangle]$ [x] is the neutrosophic polynomial ring with integer coefficients i.e. in this neutrosophic polynomial ring we have coefficients are of the form a + bI; a and b are positive or negative integers i.e. if p(x) = $a_0 + a_1 x + ... + a_n x^n$ $a_i$ ∈ $\langle Z \cup I \rangle$ and $a_i = c_i + d_i I$, $c_i$, $d_i$ ∈ Z.

***Example 2.2.2:*** Let Q be the prime field of rationals. $\langle Q \cup I \rangle$ is the neutrosophic ring of rationals. $[\langle Q \cup I \rangle]$ [x] is the neutrosophic polynomial ring.

***Example 2.2.3:*** Let R be the field of reals. $\langle R \cup I \rangle$ be the neutrosophic ring of reals, $[\langle R \cup I \rangle]$ [x] is the neutrosophic polynomial ring over the neutrosophic ring of reals.

We just give the example of complex neutrosophic polynomial ring.

***Example 2.2.4:*** Let D be the field of complex numbers. $\langle D \cup I \rangle$ be the neutrosophic ring of complex numbers. $[\langle D \cup I \rangle]$ [x] is the neutrosophic polynomial ring with coefficients from the neutrosophic complex field.

We see all the neutrosophic polynomial rings, which we have given as examples are rings of characteristic zero and they contain infinitely many elements. Now we proceed on to give examples of finite characteristic neutrosophic polynomial rings.

***Example 2.2.5:*** Let $Z_n$ be the ring of integers modulo n. $\langle Z_n \cup I \rangle$ be the neutrosophic ring of modulo integers, $[\langle Z_n \cup I \rangle]$ [x] is the polynomial neutrosophic ring of characteristic n, (n < ∞) and n any positive integer.

***Example 2.2.6:*** Let $Z_3$ = {0, 1, 2} be the ring of integers modulo 3. $\langle Z_3 \cup I \rangle$ = {0, 1, 2, I, 2I, 1 + I, 2 + I, 1 + 2I, 2 + 2I}. Clearly [$\langle Z_3 \cup I \rangle$] [x] is the neutrosophic polynomial ring with coefficients from $\langle Z_3 \cup I \rangle$. Characteristic of [$\langle Z_3 \cup I \rangle$] [x] is 3.

Now a natural question would be what is the neutrosophic polynomial ring structure. Does it have zero divisors or not? To this end first we make certain things clear. First Z is the integral domain but $\langle Z \cup I \rangle$ is not an integral domain. For consider (7 – 7I). 3I = 0; 3I ≠ 0 and 7 – 7I ≠ 0 yet 3 I (7 – 7I) = 0. Thus we see the basic nature of the ring is not carried out faithfully by the neutrosophic rings.

**THEOREM 2.2.2:** *Let F be field of characteristic zero. $\langle F \cup I \rangle$ be the neutrosophic ring. $\langle F \cup I \rangle$ has non trivial divisors of zero.*

*Proof:* Given F is a field of characteristic zero. So the set of integers Z is contained in F; further F ⊂ $\langle F \cup I \rangle$. Take x = (8 – 8I) and y = 5I, x y = (8 – 8I). 5I = 40I – 40I² = 0, (as I² = I).

Thus a neutrosophic ring $\langle R \cup I \rangle$, even when the ring R is a field or an integral domain has zero divisors.

***Example 2.2.7:*** Let $\langle Z_5 \cup I \rangle$ be the neutrosophic ring; [$\langle Z_5 \cup I \rangle$][x] is the neutrosophic polynomial ring.
Take

| | | |
|---|---|---|
| p(x) | = | (2 + 3I) $x^3$ , |
| q(x) | = | 4Ix. |
| p(x) q(x) | = | 0, |

though $Z_5$ is the prime field of characteristic 5. [$\langle Z_5 \cup I \rangle$] [x] the neutrosophic polynomial ring has zero divisors.

Based on this we have the following interesting result:

**THEOREM 2.2.3:** *A neutrosophic polynomial ring [⟨R ∪ I⟩] [x] is not an integral domain even if R is an integral domain or a field.*

The proof is left as an exercise for the reader. As in case of any polynomial in case of neutrosophic polynomial $p(x) = a_0 + a_1 x + \ldots + a_n x^n$ we say degree of $p(x) = n$ provided $a_n \neq 0$ and n is the highest degree of x; $a_0, a_1, \ldots, a_n \in \langle R \cup I \rangle$.

***Result:*** Let f(x), g(x) be any two neutrosophic polynomials in [⟨R ∪ I⟩] [x]. R is a field or an integral domain. If f(x), g(x) are two non zero polynomials of [⟨R ∪ I⟩] [x] then deg (f(x). g(x)) ≤ deg f(x) + deg g(x).

**DEFINITION 2.2.3:** *Let [⟨R ∪ I⟩] [x] be a neutrosophic polynomial ring.*
 *A neutrosophic polynomial p(x) is said to be neutrosophic reducible if p(x) = r(x) . s(x) where both r(x) and s(x) are neutrosophic polynomials; if p(x) = r(x) s(x) but only one of r(x) or s(x) is a neutrosophic polynomial then we say p(x) is only semi neutrosophic reducible. If p(x) = r(x) s(x) where r(x) = I or 1 and s(x) = p(x) or I s(x) = p(x) then we call p(x) a irreducible neutrosophic polynomial.*

**DEFINITION 2.2.4:** *Let [⟨R ∪ I⟩] [x] be a neutrosophic polynomial ring. Let I be an ideal of [⟨R ∪ I⟩][x], if I is generated by an irreducible neutrosophic polynomial, then we call I the principal neutrosophic ideal of [⟨R ∪ I⟩] [x].*

***Example 2.2.8:*** Let [⟨Z₂ ∪ I⟩] [x] = {a + bx | a, b ∈ {0, 1, I, 1 + I}}, I + (1 + I)x = p(x) is an irreducible polynomial of [⟨Z₂ ∪ I⟩][x]. The ideal generated by p(x) is a neutrosophic principal ideal of [⟨Z₂ ∪ I⟩] [x].
 The ideal generated by p (x) is a neutrosophic principal ideal of [⟨Z₂ ∪ I⟩] [x].

$$\frac{\langle Z_2 \cup I \rangle \, [x]}{N = \langle I + (1+I) \, x \rangle} \cong \{N, \, 1 + N\}.$$

Now we define neutrosophic prime ideal of a neutrosophic polynomial ring.

**DEFINITION 2.2.5:** *Let [⟨R ∪ I⟩] [x] be a neutrosophic polynomial ring. Let N be a neutrosophic polynomial ideal of [⟨R ∪ I⟩] [x]. We say N is a prime neutrosophic polynomial ideal of [⟨R ∪ I⟩] [x] if p(x) q(x) = r(x) ∈ N then p(x) or q(x) is in N.*

Obtain some examples of prime neutrosophic polynomial ideal of [⟨R ∪ I⟩] [x].

Is every principal neutrosophic polynomial ideal prime? Justify your answer. The notion of maximal and minimal ideal of a neutrosophic polynomial ring can be defined. Clearly a neutrosophic polynomial ring is not a neutrosophic integral domain.

Is [⟨R ∪ I⟩] [x] a neutrosophic unique factorization domain given R is a UFD? Justify your claim! Is division algorithm for neutrosophic polynomial rings true? Study and analyze them! Obtain any other interesting properties about neutrosophic polynomial rings.

We call two neutrosophic polynomials p(x), q(x) of the neutrosophic polynomial ring [⟨R ∪ I⟩] [x] to be relatively neutrosophic prime, if they do not have a neutrosophic polynomial r(x) such that r(x) | p(x) and r(x) | q(x) i.e. (p(x) , q(x)) = r(x) is not possible; Otherwise they have a factor or a common divisor which is a neutrosophic polynomial. It is to be noted that p(x) and q(x) can have a usual polynomial, which is a common factor, yet the neutrosophic polynomials are said to be neutrosophically relatively prime. We say a pair of two neutrosophic polynomials are said to be strongly relatively neutrosophic prime if they do not have factor which is a neutrosophic polynomial or a usual polynomial i.e. (p(x), q(x)) = 1 or I.

Several properties enjoyed by usual polynomial rings can be easily extended to the case of neutrosophic polynomial rings with appropriate modifications. It is to be noted in many a case neutrosophic polynomial rings behave differently.

Neutrosophic polynomial rings, which we have defined, happen to be commutative rings with unity. Now we proceed on to define neutrosophic matrix rings which are non commutative.

**DEFINITION 2.2.6:** *Let $\langle R \cup I \rangle$ be any neutrosophic ring. The collection of all $n \times n$ matrices with entries from $\langle R \cup I \rangle$ is called the neutrosophic matrix ring; i.e. $M_{n \times n} = \{M = (a_{ij}) \mid a_{ij} \in \langle R \cup I \rangle\}$. The operations are the usual matrix addition and matrix multiplication.*

It is left as an exercise for the reader to verify that $M_{n \times n} = \{M = (a_{ij}) \mid a_{ij} \in \langle R \cup I \rangle\}$ is a neutrosophic matrix ring. If R is of characteristic zero than $M_{n \times n}$ is a neutrosophic ring of characteristic zero. Clearly $M_{n \times n}$ has infinite number of elements and infact $M_{n \times n}$ is a non commutative neutrosophic ring.

***Example 2.2.9:*** Let $\langle Z \cup I \rangle$ be a neutrosophic ring.

$$M_{2 \times 2} = \left\{ \begin{pmatrix} a & b \\ c & d \end{pmatrix} \mid a,\, b,\, c\ d \in \langle Z \cup I \rangle \right\}$$

is the neutrosophic matrix ring of characteristic zero.

***Example 2.2.10:*** Let $\langle Z_2 \cup I \rangle = \{0,\ 1,\ I\ 1 + I\}$ be the neutrosophic ring of characteristic two.

$$M_{2 \times 2} = \left\{ \begin{pmatrix} a & b \\ c & d \end{pmatrix} \mid a,\, b,\, c,\, d \in \langle Z_2 \cup I \rangle \right\}$$

is the neutrosophic matrix ring of characteristic two and infact $|M_{2 \times 2}| < \infty$ and $M_{2 \times 2}$ is a non commutative finite neutrosophic ring. Find the number of elements in $M_{2 \times 2}$.

All properties associated with any general neutrosophic ring can be easily extended in the case of neutrosophic matrix rings. We define a neutrosophic integral domain.

**DEFINITION 2.2.7:** *Let ⟨R ∪ I⟩ be a neutrosophic ring. If ⟨R ∪ I⟩ is a commutative ring with no divisors of zero then we call ⟨R ∪ I⟩ to be an neutrosophic integral domain.*

**Example 2.2.11:** Let ⟨Z₂ ∪ I⟩ be the neutrosophic ring. ⟨Z₂ ∪ I⟩ is not a neutrosophic integral domain.

Does their exist neutrosophic integral domains?

Now we proceed on to define pseudo neutrosophic integral domain.

**DEFINITION 2.2.8:** *Let T be a commutative pseudo neutrosophic ring. If T has no zero divisors then we call T a pseudo neutrosophic integral domain.*

**Example 2.2.12:** Let P = {0, I, 2I, 3I, 4I, 5I, 6I}; P is a pseudo neutrosophic ring under addition and multiplication modulo 7, which is a pseudo neutrosophic integral domain; The operations on P being addition and multiplication modulo 7.

**Example 2.2.13:** Let S = {0, I, 2I, 3I, 4I, 5I}; S is a pseudo neutrosophic ring which is not a pseudo neutrosophic integral domain. S is a pseudo neutrosophic ring under addition and multiplication modulo 6. For 2I . 3I ≡ (0) mod 6.

**Example 2.2.14:** Let T = {0, 3, I, 2I, 3I, 4I, 5I, 3 + I, 3 + 2I, 3 + 3I, 3 + 4I, 3 + 5I}; T is a neutrosophic ring and T is not a neutrosophic integral domain but it contains a pseudo neutrosophic integral domain.

Now we proceed on to define neutrosophic division rings.

**DEFINITION 2.2.9:** *Let ⟨R ∪ I⟩ be any neutrosophic ring without zero divisors. If ⟨R ∪ I⟩ is non commutative and has no zero divisors then we all ⟨R ∪ I⟩ a neutrosophic division ring.*

***Example 2.2.15:*** Let $\langle Z_2 \cup I \rangle = \{0, 1, I, 1 + I\}$ be a neutrosophic ring. Let

$$D_I = \left\{ \begin{pmatrix} a & b \\ c & d \end{pmatrix} \mid a, b, c\ d \in \langle Z_2 \cup I \rangle \right\}$$

be the collection of all $2 \times 2$ matrices with entries from $\langle Z_2 \cup I \rangle$, $D_I$ is non commutative for take

$$x = \begin{pmatrix} I & 0 \\ 1 & I \end{pmatrix} \text{ and } y = \begin{pmatrix} 1 & I \\ I & 0 \end{pmatrix}$$

$$x.\ y = \begin{pmatrix} I & 0 \\ 1 & I \end{pmatrix} \begin{pmatrix} 1 & I \\ I & 0 \end{pmatrix} = \begin{pmatrix} I & I \\ 1+I & I \end{pmatrix}$$

$$y.\ x = \begin{pmatrix} 1 & I \\ I & 0 \end{pmatrix} \begin{pmatrix} I & 0 \\ 1 & I \end{pmatrix} = \begin{pmatrix} 0 & I \\ I & 0 \end{pmatrix}, xy \neq yx.$$

Hence $D_I$ is a non commutative neutrosophic ring. Clearly $D_I$ is not a division ring for take

$$a = \begin{pmatrix} 0 & 1 \\ 0 & 0 \end{pmatrix} \text{ and } b = \begin{pmatrix} 1 & 0 \\ 0 & 0 \end{pmatrix} \text{ in } D_I$$

$$a.\ b = \begin{pmatrix} 0 & 1 \\ 0 & 0 \end{pmatrix} \begin{pmatrix} 1 & 0 \\ 0 & 0 \end{pmatrix} = \begin{pmatrix} 0 & 0 \\ 0 & 0 \end{pmatrix}.$$

Consider

$$b.a = \begin{pmatrix} 1 & 0 \\ 0 & 0 \end{pmatrix} \begin{pmatrix} 0 & 1 \\ 0 & 0 \end{pmatrix} = \begin{pmatrix} 0 & 1 \\ 0 & 0 \end{pmatrix} \neq \text{ zero matrix.}$$

Thus we see a.b = (0) but b.a ≠ (0). So a.b is only a one sided zero divisor of $D_I$. Based on this observation we would be defining the notion of left or right zero divisors in a neutrosophic ring. We would also define neutrosophic zero divisors and semi neutrosophic zero divisors of a neutrosophic ring.

Likewise we would also define idempotents and neutrosophic idempotents in a neutrosophic ring. Consider

$$x = \begin{pmatrix} 1 & 0 \\ 0 & 0 \end{pmatrix} \in D_I. \quad x^2 = \begin{pmatrix} 1 & 0 \\ 0 & 0 \end{pmatrix} \begin{pmatrix} 1 & 0 \\ 0 & 0 \end{pmatrix} = \begin{pmatrix} 1 & 0 \\ 0 & 0 \end{pmatrix} = x.$$

We call x just an idempotent of the neutrosophic ring. Take

$$y = \begin{pmatrix} 1+I & 0 \\ 0 & 0 \end{pmatrix} \in D_I,$$

$$y^2 = \begin{pmatrix} 1+I & 0 \\ 0 & 0 \end{pmatrix} \begin{pmatrix} 1+I & 0 \\ 0 & 0 \end{pmatrix} = \begin{pmatrix} 1+I & 0 \\ 0 & 0 \end{pmatrix} = y,$$

y is called the neutrosophic idempotent of the neutrosophic ring. Does $D_I$ have a proper subset, which is a division ring or a neutrosophic division ring?

**DEFINITION 2.2.10:** *Let D be a division ring. The ring generated by $\langle D \cup I \rangle$ is called the neutrosophic division ring. If $\langle D \cup I \rangle$ has a proper subset P and P is a neutrosophic ring and is a division ring then we call $\langle D \cup I \rangle$ a weak neutrosophic division ring.*

***Note:*** We need this definition because even if D is a division ring we may not have $\langle D \cup I \rangle$ to be a division ring. Is it possible to find a neutrosophic division ring $\langle D \cup I \rangle$ where D is a division ring?

Consider the neutrosophic ring in this example.

***Example 2.2.16:*** Let $\langle Z_2 \cup I \rangle = \{0, 1, I, 1 + I\}$ is the neutrosophic ring. P = $\{0, 1\}$, T = $\{0, I\}$, S = $\{0, 1 + I\}$ are division rings in fact integral domains but none of them are neutrosophic division rings or neutrosophic integral domains.

***Example 2.2.17:*** Let $\langle Z_3 \cup I \rangle$ = {0, 1, 2, I, 2I, 1 + I, 1 + 2I, 2 + I, 2 + 2I} be the neutrosophic ring, P = {0, I, 2I} is an integral domain which is not a neutrosophic integral domain or usual integral domain.

We make from the above example one more observation viz. (2 + 2I) (2 + 2I) = 1 Also $(1 + I)^2 = 1$, $(1 + 2I)^2 = 1 + 2I$, 2I (1 + I) = I, 2I (1 + 2I) = 0, I (2 + I) = 0, I (2 + 2I) = I. We make provisions for all these by defining new notions.

Now we proceed on to define the notion of neutrosophic ring of real quaternions, which will give another example of a non commutative neutrosophic ring apart from the neutrosophic ring of matrices.

**DEFINITION 2.2.11:** *Let $Q = \{\alpha_0 + \alpha_1 i + \alpha_2 j + \alpha_3 k$ where $\alpha_0$, $\alpha_1$, $\alpha_2$, $\alpha_3$, are real numbers, $i^2 = k^2 = j^2 = ijk = -1$, $ij = -ji = k$, $jk = -kj = i$, $ki = -ik = j\}$. Q is a non commutative ring called the ring of real quaternions infact a division ring. The neutrosophic ring generated by Q and I is called the neutrosophic ring of real quaternions. It is also a ring and not a neutrosophic division ring.*

We know $\langle R \cup I \rangle$, $\langle Q \cup I \rangle$ and $\langle Z \cup I \rangle$ are neutrosophic rings. Clearly $\langle Z \cup I \rangle$ is a neutrosophic ring of integers. Suppose ($a_0$, $a_1$, …, $a_n$) $\in$ Z we say ($a_0$, $a_1$, …, $a_n$) is primitive if the greatest common divisor of $a_0$, $a_1$, …, $a_n$ is 1.

Let $a_0$, $a_1$, …, $a_n$ $\in$ $\langle Z \cup I \rangle$ where $a_0$, $a_1$, …, $a_n$ can be integers or neutrosophic integers or elements of the form x + yI (x, y $\in$ Z).

We say the set is primitive if the greatest common divisor is 1. Consider {2 + I, 3 + 2I, 5 − I, I, 8, 31 − 2I} $\in$ $\langle Z \cup I \rangle$. This set is a neutrosophic set of numbers and their greatest common divisor is 1.

Consider {3, 3I, 3 − 21I, 21 + 6I, 27 − 9I, 9I + 18, 18} $\in$ $\langle Z \cup I \rangle$. The greatest common divisor is 3.

Now using this fact we define the notion of primitive polynomial.

**DEFINITION 2.2.12:** *Let $p(x) = a_0 + a_1x + ... + a_nx \in \langle R \cup I \rangle$ [x] where $\{a_0, a_1, a_2, ..., a_n\} \in \langle Z \cup I \rangle$. We say $p(x)$ is a primitive neutrosophic polynomial if the greatest common divisor of $\{a_0, a_1, ..., a_n\}$ is 1.*

Is the fact if f(x) and g(x) are primitive neutrosophic polynomials then f(x) g(x) is a primitive polynomial true, in case of neutrosophic polynomial rings?

***Example 2.2.18:*** We just illustrate with some examples, Let $\langle Z \cup I \rangle$[x] be a neutrosophic polynomial ring. Take

$$p(x) = 3I + (2 + 5I) x + 5x^2 \text{ and}$$
$$q(x) = 2 + 3Ix + (5 - 2I) x^3.$$

We see both p(x) and q(x) are primitive polynomials.

Consider the product

$$p(x) \, q(x) = 6I + (9I + 4 + 10I) x + (6I + 15I + 10) x^2 + (15I - 6I + 15I) x^3. = 6I + (4 + 19I) x + (10 + 21I) x^2 + 24I x^3.$$

Clearly p(x) q(x) is also a primitive neutrosophic polynomial of $\langle Z \cup I \rangle$ [x].

We define the content of the neutrosophic polynomial $p(x) = a_0 + a_1x + ... + a_nx^n \in \langle Z \cup I \rangle$ [x] to be the greatest common divisor (g.c.d) of $a_0, a_1, ..., a_n$. The g.c.d can be a real integer or a neutrosophic number.
For we illustrate this by the following example:
Take
$$p(x) = (2 + I) + 7 (2 + I) x^2 + (4 - I) x^3.$$
The content of p(x) is a neutrosophic number given by 2 + I; as
$$p(x) = (2 + I) + 7 (2 + I) x^2 + (2 + I) (2 - I) x^3.$$
Take
$$q(x) = (8 - 16I) + 4I x^2 + 40x^3 - 16 Ix^4.$$
The content of q(x) = 4.

Can Gauss Lemma for primitive polynomials be proved in case of primitive neutrosophic polynomials in $\langle Z \cup I \rangle$ [x]?

We call a neutrosophic polynomial integer monic if all coefficients are from $\langle Z \cup I \rangle$ and its highest coefficient is 1. It is pertinent to mention here that a neutrosophic integer monic polynomial is also primitive.

Now as in case of polynomial rings we in case of neutrosophic polynomial rings also define polynomial rings in n-variables.

By a neutrosophic polynomial ring in the variable $x_1$ over $\langle R \cup I \rangle$, a neutrosophic commutative ring with 1 we mean a set of formal symbols $a_0 + a_1x + \ldots + a_mx^m$ where $a_0, a_1, \ldots, a_m$ are in $\langle R \cup I \rangle$, we have proved $\langle R \cup I \rangle [x_1]$ is also a commutative ring with unit.

Let $R_1 = \langle R \cup I \rangle[x_1]$. Consider a variable $x_2$, $R_1 [x_2]$ is a neutrosophic polynomial ring over $R_1$ in the variable $x_2$, proceeding in this way $R_n = R_{n-1} [x_n]$ and $R_n$ is a neutrosophic polynomial ring over $R_{n-1}$ in the variable $x_n$. Thus $R_n$ is a ring of neutrosophic polynomials in the variables $x_1, x_2, \ldots, x_n$ over $\langle R \cup I \rangle$. Its element is of the form

$$\sum a_{i_1, i_2, \ldots, i_n} \ x_1^{i_1}, \ x_2^{i_2} \ldots x_n^{i_n} \text{ with}$$

$a_{i_1, i_2, \ldots, i_n} \in \langle R \cup I \rangle$, $i_1, \ldots, i_n$ are integers greater than or equal to 0. Multiplication is defined by

$$\left( x_1^{i_1}, x_2^{i_2}, \ldots, x_n^{i_n} \right) \left( x_1^{j_1}, x_2^{j_2}, \ldots, x_n^{j_n} \right) = x_1^{i_1+j_1}, \ x_2^{i_2+j_2}, \ldots, x_n^{i_n+j_n}$$

$R_n$ is called the neutrosophic ring of polynomials in n-variables $x_1, \ldots, x_n$ over the ring R.

***Example 2.2.19:*** Let $(\langle Z \cup I \rangle) [x_1, x_2] = R_2$ is the neutrosophic ring of polynomials in the 2 variables $x_1$ and $x_2$.

$p (x_1, x_2) = (2 + I) + (3 - 4I) \ x_1^2 x_2 + 5I \ x_1^4 \ x_2^3 + (2 - 6I) \ x_1^7 \ x_2^5$.

Now we can also replace Z in the example 2.2.19 by Q or R.

***Example 2.2.20:*** Let $\langle Z_3 \cup I \rangle [x_1, x_2]$ is the neutrosophic polynomial ring in two variables $x_1$ and $x_2$. $P (x_1, x_2) = (2 + I) + Ix_1x_2 + (2 + 2I) \ x_2^2 + (I + 1) \ x_1^3 \ x_2^4$. These types of neutrosophic polynomial rings in several variables can be studied over finite characteristic neutrosophic rings as well as infinite neutrosophic rings and also of zero characteristic.

Such study may lead to some new means of solving equations, which has solutions involving indeterminacy.

***Example 2.2.21:*** Let $\langle Z_2 \cup I \rangle$ [$x_1$, $x_2$, $x_3$] be a neutrosophic polynomial ring in 3 variables $x_1$, $x_2$, $x_3$ over

$$\langle Z_2 \cup I \rangle = \{0, 1, I, 1 + I\}.$$

Take the polynomial

$$p(x_1, x_2, x_3) = x_1 x_2 + x_1 + (1 + I) x_2 + 1 + I.$$

Is this polynomial factorizable or not;

$$p(x_1, x_2, x_3) = x_1 (x_2 + 1) + (1 + I) (x_2 + 1)$$
$$= [x_1 + (1 + I)] [x_2 + 1].$$

But this sort of factorization is not always possible. The interested reader can construct means and methods to factorize neutrosophic polynomials in single and several variables.

**DEFINITION 2.2.13:** *Let $\langle R \cup I \rangle$ be a neutrosophic ring. An element x in $\langle R \cup I \rangle$ where $x = a + bI$, $a \neq b$ or $-b$ is said to be a neutrosophic zero divisor if we can find $y = c + dI$, $c \neq d$ or $-d$ in $\langle R \cup I \rangle$ such that $x y = y x = 0$.*

*If $xy = 0$ and $y x \neq 0$ we call x a right zero divisor, y is its right zero divisor.*

On similar lines we can define left neutrosophic zero divisor. Thus we see y is a neutrosophic zero divisor of x if and only if y is both a neutrosophic right zero divisor as well as neutrosophic left zero divisor of x.

We illustrate both neutrosophic zero divisors and right (left) neutrosophic zero divisors. First it is to be noted the concept of right (left) neutrosophic zero divisors occur only when the neutrosophic ring is a non commutative one. It is still important to note that all zero divisors of the form $(a - aI) I = aI - aI^2 = 0$ ($\because$ $I^2 = I$) will be termed as trivial neutrosophic zero divisors of $\langle R \cup I \rangle$. For they arise from the basic fact $I^2 = I$.

***Example 2.2.22:*** Let $\langle Z_4 \cup I \rangle$ = {0, 1, 2, 3, I, 2I, 3I, 1 + I, 2 + I, 3 + I, 1 + 2I, 2 + 2I, 3 + 3I, 1 + 3I, 2 + 3I, 3 + 3I} be a neutrosophic ring of characteristic four.

Clearly (2 + 2I) I = 0 is a trivial neutrosophic zero divisor. $(2 + 2I)^2 = 0$ is a non trivial neutrosophic zero divisor of $\langle Z_4 \cup I \rangle$. (2 + 3I) (2 + 2I) = 0 is a non trivial neutrosophic divisor of zero; 2I. 2 = 0 is also a non trivial neutrosophic divisor of zero, which we name differently. (2 + 2I) . 2 = 0 is also a non-trial neutrosophic divisor of zero, which will be defined differently.

**DEFINITION 2.2.14:** *Let $\langle R \cup I \rangle$ be a neutrosophic ring $x = a + bI$, be a neutrosophic element of R ($a \neq 0$ and $b \neq 0$). If $y \in R$ is such that $x\,y = y\,x = 0$ then $y$ is called the semi neutrosophic divisor of zero.*

***Note:*** An element $x \in \langle R \cup I \rangle$ can have both neutrosophic zero divisors and semi neutrosophic zero divisors. For we see in example 2.2.22. (2 + 2I) has both neutrosophic divisors of zero as well as semi neutrosophic divisors of zero.

***Example 2.2.23:***

$$M_{2 \times 2} = \left\{ (a_{ij}) = \begin{pmatrix} a & b \\ c & d \end{pmatrix} \text{where } a, b, c, d \in \langle Z_2 \cup I \rangle = \{0, 1, I, 1 + I\} \right\}.$$

$M_{2 \times 2}$ is a neutrosophic ring which is non commutative. Take

$$x = \begin{pmatrix} 1+I & 0 \\ 0 & 0 \end{pmatrix} \text{ and } \begin{pmatrix} 0 & 1+I \\ 0 & 0 \end{pmatrix} = y\,.$$

$$x.y = \begin{pmatrix} 1+I & 0 \\ 0 & 0 \end{pmatrix} \begin{pmatrix} 0 & 1+I \\ 0 & 0 \end{pmatrix}$$

$$= \begin{pmatrix} 0 & 1+I \\ 0 & 0 \end{pmatrix} \neq \begin{pmatrix} 0 & 0 \\ 0 & 0 \end{pmatrix}.$$

$$y\,x = \begin{pmatrix} 0 & 1+I \\ 0 & 0 \end{pmatrix}\begin{pmatrix} 1+I & 0 \\ 0 & 0 \end{pmatrix} = \begin{pmatrix} 0 & 0 \\ 0 & 0 \end{pmatrix}.$$

This x is only a one sided neutrosophic divisor of zero, consider

$$a = \begin{pmatrix} 1 & 0 \\ 0 & 0 \end{pmatrix} \text{ and } b = \begin{pmatrix} 0 & 1+I \\ 0 & 0 \end{pmatrix} \text{ in } M_{2\times 2}.$$

$$b\,a = \begin{pmatrix} 0 & 1+I \\ 0 & 0 \end{pmatrix}\begin{pmatrix} 1 & 0 \\ 0 & 0 \end{pmatrix} = \begin{pmatrix} 0 & 0 \\ 0 & 0 \end{pmatrix}.$$

$$a\,b = \begin{pmatrix} 1 & 0 \\ 0 & 0 \end{pmatrix}\begin{pmatrix} 0 & 1+I \\ 0 & 0 \end{pmatrix} = \begin{pmatrix} 0 & 1+I \\ 0 & 0 \end{pmatrix} \neq \begin{pmatrix} 0 & 0 \\ 0 & 0 \end{pmatrix}.$$

Thus $a = \begin{pmatrix} 1 & 0 \\ 0 & 0 \end{pmatrix}$ is a semi neutrosophic divisor of zero. For both x and a. The same is true of

$$y = \begin{pmatrix} 0 & 1+I \\ 0 & 0 \end{pmatrix},$$

acts to produce a neutrosophic zero divisor as well as semi neutrosophic zero divisor. Also this has several other one sided divisors of zero. For instance take

$$c = \begin{pmatrix} 0 & 0 \\ 1+I & 1+I \end{pmatrix} \text{ and } d = \begin{pmatrix} 0 & 1+I \\ 0 & 0 \end{pmatrix}$$

$$c\,d = \begin{pmatrix} 0 & 0 \\ 1+I & 1+I \end{pmatrix}\begin{pmatrix} 1+I & 0 \\ 0 & 0 \end{pmatrix} = \begin{pmatrix} 0 & 0 \\ 1+I & 0 \end{pmatrix} \neq \begin{pmatrix} 0 & 0 \\ 0 & 0 \end{pmatrix}.$$

Chapter Three

# NEUTROSOPHIC GROUP RINGS AND THEIR GENERALIZATIONS

This chapter introduces the notion of neutrosophic group rings; their special generalizations and neutrosophic semigroup rings and their generalizations. This chapter has three sections. Section one introduces the notion of neutrosophic group rings K⟨G ∪ I⟩, i.e., the neutrosophic group ⟨G ∪ I⟩ over the ring K. Section two studies the special properties of neutrosophic group rings. Section three introduces several concepts like semigroup neutrosophic rings, S-semigroup neutrosophic rings neutrosophic semigroup neutrosophic rings some of which special generalizations of group rings and analyses them.

## 3.1 Neutrosophic Group Rings

In this section we for the first time introduce the notion of neutrosophic group rings. Neutrosophic group rings are defined analogous to groups over rings i.e., group rings. Here groups are replaced by neutrosophic groups in the definition of group rings, these algebraic structures can be realized as a special type of generalization of group rings. Throughout this section we take neutrosophic group over rings unless otherwise a special mention is made.

**DEFINITION 3.1.1:** *Let* ⟨*G* ∪ *I*⟩ *be any neutrosophic group. R any ring with 1 which is commutative or a field. We define the neutrosophic group ring R (⟨G* ∪ *I*⟩*) of the neutrosophic group* ⟨*G* ∪ *I*⟩ *over the ring R as follows:*

1. *R(⟨G* ∪ *I*⟩) consists of all finite formal sums of the form*

$$\alpha = \sum_{i=1}^{n} r_i g_i, \ n < \infty \ ; \ r_i \in R \ and \ g_i \in \langle G \cup I \rangle \ (\alpha \in R(\langle G \cup I \rangle)).$$

2. *Two elements* $\alpha = \sum_{i=1}^{n} r_i g_i$ *and* $\beta = \sum_{i=1}^{m} s_i g_i$ *in R⟨G* ∪ *I⟩ are equal if and only if $r_i = s_i$ and m = n. Thus having defined distinctness of elements we proceed onto define sum and product in R* ⟨*G* ∪ *I*⟩*.*

3. *Let*

$$\alpha = \sum_{i=1}^{n} \alpha_i g_i \, , \ \beta = \sum_{i=1}^{m} \beta_i g_i \ \in R \ \langle G \cup I \rangle;$$

$$\alpha + \beta = \sum_{i=1}^{n} \ (\alpha_i + \beta_i) g_i \in R \ \langle G \cup I \rangle;$$

*as $\alpha_i$, $\beta_i \in R$ so $\alpha_i + \beta_i \in R$ and $g_i \in \langle G \cup I \rangle$.*

4. $0 = \sum_{i=1}^{n} 0 g_i$ *serves as the zero of R(⟨G* ∪ *I*⟩*).*

5. *Let*

$$\alpha = \sum_{i=1}^{n} \alpha_i g_i \ \in R \ \langle G \cup I \rangle$$

*then*

$$-\alpha = \sum_{i=1}^{n} (-\alpha_i) g_i$$

*is such that*

$$\begin{aligned} \alpha + (-\alpha) \ &= \ 0 \\ &= \ \Sigma \, (\alpha_i + (-\alpha_i)) g_i \\ &= \ \Sigma \, 0 g_i \, . \end{aligned}$$

*Thus we see R ⟨G ∪ I⟩ is an abelian group under the binary operation '+'.*

6. *We define product of two element α, β in R⟨G ∪ I⟩, as follows :*
   *Let*

   $$\alpha = \sum_{i=1}^{n} \alpha_i g_i \ \ and \ \ \beta = \sum_{j=1}^{m} \beta_j h_j$$

   $$\alpha \cdot \beta = \sum_{\substack{i \leq i \leq n \\ 1 \leq j \leq m}} \alpha_i \cdot \beta_j g_i h_j$$

   $$= \sum_{K} Y_K \, t_K$$

   *(where $Y_K = \sum \alpha_i \beta_j$ with $g_i h_j = t_K$, $t_K \in ⟨G \cup I⟩$ and $Y_K \in R$). Clearly $\alpha \cdot \beta \in R⟨G \cup I⟩$.*

7. *It is an easy task to verify that α (β + γ) = αβ + αγ and (β + γ) α = βα + γα for all α, β, γ ∈ R(⟨G ∪ I⟩). That is the distributive law is inherited by the algebraic structure.*

   *Clearly. R ⟨G ∪ I⟩ is a ring under the binary operation + and '·'.*

**Remark:**

1. The neutrosophic group ring R⟨G ∪ I⟩ is commutative if and only if ⟨G ∪ I⟩ is a commutative neutrosophic group.
2. The neutrosophic group ring R⟨G ∪ I⟩ has finite number of elements or is of finite order if and only if both the ring and the neutrosophic group ⟨G ∪ I⟩ are finite. Even if one of them is infinite the neutrosophic group ring R⟨G ∪ I⟩ is of infinite order.
3. Since 1 ∈ R we can always have ⟨G ∪ I⟩ to be a proper subset of R⟨G ∪ I⟩ i.e.1.g → g.

Now we will see illustrations of such structures before we proceed onto prove some of its related properties.

**Example 3.1.1:** Let R = Z, be the ring of integers. Take $\langle G \cup I \rangle$ = {1, g, g², I, gI, g²I / g³ = 1 and I² = I}.

$$Z\langle G \cup I \rangle = \left\{ \sum_{i=1}^{n} r_i g_i, \quad n \leq 6 \;\middle/\; g_i \in \langle G \cup I \rangle \text{ and } r_i \in Z \right\}$$

$Z\langle G \cup I \rangle$ is a ring called the neutrosophic group ring. The order of $Z\langle G \cup I \rangle$ is infinite.

Clearly Z is an integral domain but $Z\langle G \cup I \rangle$ is not a domain for it has zero divisors. Further $Z\langle G \cup I \rangle$ is a commutative neutrosophic group ring of infinite order.

Now we proceed onto give some examples of finite neutrosophic group rings.

**Example 3.1.2:** Let

$\langle G \cup I \rangle$ = {1, g, g², g³, g⁴, I, gI, g²I, g³I, g⁴I, / g⁵ = 1, I² = I}

be a neutrosophic group of order 10. Let $Z_2$ = {0, 1} be the ring of characteristic 2. $Z_2 \langle G \cup I \rangle$ is a neutrosophic ring of finite order, which is a commutative ring.

Clearly $\langle G \cup I \rangle \subseteq Z_2 \langle G \cup I \rangle$ and $Z_2 \subseteq Z_2 (\langle G \cup I \rangle)$.

Now we proceed to give yet another example.

**Example 3.1.3:** Let $Z_2$ = {0,1} the ring of characteristic two.

$\langle G \cup I \rangle$ = {1,g, I, g I / g² = 1 and I² = I}.

$Z_2 \langle G \cup I \rangle$ = {0, 1, g, I, gI, 1 + g, 1 + g + I, 1 + gI + I, 1 + I, g + gI + I, I + gI, g + I, g + gI, 1 + gI, 1 + g + gI, 1 + g + I + gI}.

$Z_2 \langle G \cup I \rangle$ is a finite commutative neutrosophic group ring having zero divisors for $(1 + g)^2 = 1 + 2g + g^2 = 0$, $(I + gI)^2 = 0$. This neutrosophic group ring has nontrivial idempotents, viz. $(1 + I)^2 = 1 + 2I + I^2 = 1 + I$ we call $I^2 = I$ to be a trivial idempotent of the neutrosophic group ring.

Now we proceed on to illustrate a non commutative neutrosophic group ring $R \langle G \cup I \rangle$.

***Example 3.1.4:*** Let $Z_6$ be the ring of integers module 6 i.e.,

$Z_6 = \{0, 1, 2, 3, 4, 5\}$ and $\langle G \cup I \rangle$
$= \{1, p_1, p_2, p_3, p_4, p_5, I, p_1I, p_2I, ..., p_5I\}$ where

$$1 = \begin{pmatrix} 1 & 2 & 3 \\ 1 & 2 & 2 \end{pmatrix}, \ p_1 = \begin{pmatrix} 1 & 2 & 3 \\ 1 & 3 & 2 \end{pmatrix},$$

$$p_2 = \begin{pmatrix} 1 & 2 & 3 \\ 3 & 2 & 1 \end{pmatrix}, \ p_3 = \begin{pmatrix} 1 & 2 & 3 \\ 2 & 1 & 3 \end{pmatrix},$$

$$p_4 = \begin{pmatrix} 1 & 2 & 3 \\ 2 & 3 & 1 \end{pmatrix} \text{ and } \ p_5 = \begin{pmatrix} 1 & 2 & 3 \\ 3 & 1 & 2 \end{pmatrix}.$$

Clearly $Z_6[\langle G \cup I \rangle]$ is a non commutative neutrosophic group ring of finite order having zero divisors and idempotents.

The following theorems are worth mentioning:

**THEOREM 3.1.1:** *Every neutrosophic group ring $R \langle G \cup I \rangle$ contains atleast one proper subset, which is a group ring.*

*Proof:* Given $R(\langle G \cup I \rangle)$ is a neutrosophic group ring. Clearly, $RG \subseteq R(\langle G \cup I \rangle)$, so $R(\langle G \cup I \rangle)$ contains a group ring.

Now we give the condition under which the neutrosophic group ring $R(\langle G \cup I \rangle)$ has one and only one proper subset, which is a group ring.

**THEOREM 3.1.2:** *Let R(⟨G ∪ I⟩) be a neutrosophic group ring of the neutrosophic group ⟨G ∪ I⟩ over the ring R. R(⟨G ∪ I⟩) has one and only one proper subset which is a group ring if the following condition are satisfied.*

1.  *R is a ring having no proper subrings or R is a prime field.*
2.  *G is a group having no proper subgroups.*

*Equivalently if conditions (1) and (2) are satisfied the neutrosophic group ring has one and only one proper subset which is a group ring.*

*Proof:* First we make a mention of what we mean by a proper subset being a group ring in a neutrosophic group ring R(⟨G ∪ I⟩).

A proper subset φ ≠ A of R(⟨G ∪ I⟩) is said to be group ring if A can be written as A = KH where K is a commutative ring with unity and a proper subset of R or a subfield and a proper subset of R or K = R. H is a proper subgroup of G or H = G, (A ⊊ R[⟨G ∪ I⟩]).

If G has no proper subgroups and R has no proper subrings with unit then RG is the only group ring of the neutrosophic group ring R(⟨G ∪ I⟩). Thus we see conditions (1) and (2) of the theorem are trivially satisfied.

If on the other hand conditions (1) and (2) are satisfied to show the neutrosophic group ring can have one and only proper subset which is a group ring.

We are given conditions (1) and (2) of the theorem are satisfied i.e. G has no proper subgroups other than itself and identity and R has no proper subset, which is a ring with unity. So the only group ring in the neutrosophic group ring R(⟨G ∪ I⟩) is RG. Hence the claim.

First we show we have a nontrivial class of neutrosophic group rings, which satisfies this condition. We call this class of neutrosophic group rings as a special class of neutrosophic rings or the group rings is known as the special neutrosophic group ring.

**THEOREM 3.1.3:** *Let $Z_p$ be the prime field of characteristic p or Z the ring of integers and G a cyclic group of prime order p. Then the neutrosophic group ring $Z_p$ $(G \cup I)$ $(Z(\langle G \cup I \rangle))$ has only one proper subset which is a group ring given by $Z_pG$ or ZG.*

*Proof:* Follows directly from the earlier theorem and the fact $Z_p$ is a prime field of characteristic p (Z has no proper subrings with unity) and G a cyclic group of prime order.

***Note:*** If we take Q to be the prime field of characteristic 0 we see $Z \underset{\neq}{\subset} Q$, Z is a ring with unit and ZG is the group ring also QG is a group ring. That is why we have not taken Q, only Z.

Before we proceed onto define other concepts we first illustrate this notion by the following examples:

***Example 3.1.5:*** Let Z be the ring of integers $\langle G \cup I \rangle = \{1, g, g^2, \ldots, g^6, IgI, g^2I, \ldots, g^5I / g^7 = 1 \text{ and } I^2 = I\}$.

$Z\langle G \cup I \rangle$ is a special neutrosophic group ring for ZG is the only group ring in $Z(\langle G \cup I \rangle)$.

***Example 3.1.6:*** Let $Z_5 = \{0,1, 2,3, 4\}$ be the prime field of characteristic 5 and $\langle G \cup I \rangle = \{1, g, g^2, I, gI, g^2I / g^3 = 1, I^2 = I\}$.

$Z_5 \langle G \cup I \rangle$ is a special neutrosophic group ring.

Thus the concept of special neutrosophic group rings divides the class of neutrosophic rings into two distinct (i.e. non overlapping classes) classes.

***Example 3.1.7:*** Let Q be the prime field of characteristic 0. $\langle G \cup I \rangle = \langle 1, g, g^2, g^3, g^4, I, gI, g^2I, g^3I, g^4I / g^5 = 1 \text{ and } I^2 = I \rangle$. $Q(\langle G \cup I \rangle)$ is not a special neutrosophic group ring for it contains QG and ZG as proper group rings.

In case of neutrosophic group rings we can define five types of subrings.

1. Subneutrosophic group rings.
2. Neutrosophic subrings.
3. Pseudo neutrosophic subrings.
4. Subgroup rings.
5. Subrings.

Now we proceed on to define these five concepts and illustrate with examples and also obtain conditions for them to contain such types of subrings.

**DEFINITION 3.1.2:** *Let R(⟨G ∪ I⟩) be a neutrosophic group ring of the neutrosophic group ⟨G ∪ I⟩ over the ring R. We call a proper subset P of R(⟨G ∪ I⟩) to be a subneutrosophic group ring if P = R(⟨H ∪ I⟩) or S(⟨G ∪ I⟩) or T(⟨H ∪ I⟩), in P = R(⟨H ∪ I⟩) i.e. R is the ring and ⟨H ∪ I⟩ is a proper neutrosophic subgroup of ⟨G ∪ I⟩ or in S(⟨G ∪ I⟩). S is a proper subring with unit of the ring R and ⟨G ∪ I⟩ is the neutrosophic group and if P = T(⟨H ∪ I⟩); T is a subring of R with unit and ⟨H ∪ I⟩ is a neutrosophic subgroup, H is a proper subgroup of G.*

We first illustrate this by the following example:

***Example 3.1.8:*** Let Q(⟨G ∪ I⟩) be a neutrosophic group ring where Q = the field of rationals and
⟨G ∪ I⟩ = ⟨1, g, g$^2$, g$^3$, g$^4$, g$^5$, I, gI, …, g$^5$I / g$^6$ = 1 and I$^2$ = I⟩.

Take P = Q(⟨H ∪ I⟩) where
⟨H ∪ I⟩ = {1, g$^2$, g$^4$, I, g$^2$I, g$^4$I}.
P is a subneutrosophic group ring.

Take T = Z⟨G ∪ I⟩, Z is a ring with unit, Z⟨G ∪ I⟩ is also a subneutrosophic group ring.

Consider Z(⟨P ∪ I⟩) where
⟨P ∪ I⟩ = {1, g$^3$, I, g$^3$I / g$^6$ = 1 and I$^2$ =I}.

Z($\langle P \cup I \rangle$) is also a subneutrosophic group ring. Thus this neutrosophic group ring has all the three types of subneutrosophic group rings.

Next we proceed on to define neutrosophic subrings.

**DEFINITION 3.1.3:** *Let R($\langle G \cup I \rangle$) be a neutrosophic group ring. A proper subset P of R($\langle G \cup I \rangle$) is said to be a neutrosophic subring if P = $\langle S \cup I \rangle$ where S is a subring of RG or R.*

We first illustrate this by the following example:

***Example 3.1.9:*** Let R($\langle G \cup I \rangle$) be a neutrosophic group ring of the neutrosophic group $\langle G \cup I \rangle$ over the ring R.

Clearly $\langle R \cup I \rangle$ is a neutrosophic subring of R($\langle G \cup I \rangle$).

Now we have the following interesting theorem:

**THEOREM 3.1.4:** *Let R($\langle G \cup I \rangle$) be the neutrosophic group ring of the neutrosophic group $\langle G \cup I \rangle$ over the ring R. Then R ($\langle G \cup I \rangle$) always has a nontrivial neutrosophic subring.*

*Proof:* Consider the neutrosophic ring generated by $\langle R \cup I \rangle$ clearly $\langle R \cup I \rangle \subset$ R($\langle G \cup I \rangle$). Hence the claim.

Now we proceed onto define pseudo neutrosophic subrings.

**DEFINITION 3.1.4:** *Let R($\langle G \cup I \rangle$) be a neutrosophic group ring of the neutrosophic group $\langle G \cup I \rangle$ over the ring R. If R($\langle G \cup I \rangle$) has a proper subset T which is just a pseudo neutrosophic ring, then we say the neutrosophic group ring has pseudo neutrosophic subring. i.e. T $\neq \langle S \cup I \rangle$ where S is a subring or subgroup ring of R or RG respectively.*

We just illustrate this by the following example:

***Example 3.1.10:*** Let $Z_6(\langle G \cup I \rangle)$ be a neutrosophic group ring of the neutrosophic group $\langle G \cup I \rangle$ over $Z_6$, the ring of integers modulo 6. Clearly $P = \{0, 3I\}$ is a proper subset of $Z_6(\langle G \cup I \rangle)$ which is a pseudo neutrosophic subring of $Z_6 (\langle G \cup I \rangle)$.

Next we proceed onto define the notion of subgroup rings.

**DEFINITION 3.1.5:** *Let $R(\langle G \cup I \rangle)$ be the neutrosophic group ring of the neutrosophic group $\langle G \cup I \rangle$ over the ring R. We say a proper subset P of $R(\langle G \cup I \rangle)$ to be a subgroup ring of $R(\langle G \cup I \rangle)$ if P = SH where S is a subring of R and H is just a subgroup of G and SH is the group ring of the subgroup H over the ring S.*

*Clearly all neutrosophic group rings have proper subgroup ring. For if $R(\langle G \cup I \rangle)$ is a neutrosophic group ring then RG is the subgroup ring of $R\langle G \cup I \rangle$.*

Now we proceed onto define the notion of subrings of the neutrosophic group ring. Clearly we assume these subrings do not have the structure of subgroup rings. They are only subrings.

**DEFINITION 3.1.6:** *Let $R(\langle G \cup I \rangle)$ be the neutrosophic group ring of the neutrosophic group $\langle G \cup I \rangle$ over the ring R. A proper subset P which is a subring but P should not have a group ring structure, is defined to be the subring of the neutrosophic group ring $R(\langle G \cup I \rangle)$.*

We illustrate this by the following example:

***Example 3.1.11:*** Let $Z_2\langle G \cup I \rangle$ be the neutrosophic group ring where $Z_2 = \{0, 1\}$ and $\langle G \cup I \rangle = \{1, g, g^2, g^3, I, gI, g^2I, g^3I / g^4 = 1$ and $I^2 = I\}$.

Clearly $P = \{0, 1 + g^2\}$ is a subring of $Z_2 (\langle G \cup I \rangle)$.

Here it has become important to mention that when we define ideals in a neutrosophic group ring we cannot have any of these types we can have only one type. Thus we see in case of neutrosophic group rings we have the vast difference between the ideals and subrings.

Now we proceed onto define neutrosophic ideals of a neutrosophic group ring R(⟨G ∪ I⟩).

**DEFINITION 3.1.7:** *Let R(⟨G ∪ I⟩) be the neutrosophic group ring of the neutrosophic group ⟨G ∪ I⟩ over the ring R. We say a proper subset J of R(⟨G ∪ I⟩) is said to be a neutrosophic ideal of R(⟨G ∪ I⟩)if (1) J is a neutrosophic subring or a sub neutrosophic group ring of R(⟨G ∪I⟩).*

*For every j ∈ J and α ∈ R(⟨G ∪I⟩) αj and jα belong to J.*

***Note:*** One can define the notion of right and left neutrosophic ideals of the neutrosophic group ring R(⟨G ∪ I⟩).

In case the neutrosophic group ring is commutative the notion of left and right neutrosophic ideals coincide.

Now we have yet another notion called pseudo neutrosophic ideal of the neutrosophic group ring R(⟨G ∪ I⟩).

**DEFINITION 3.1.8:** *Let R(⟨G ∪ I⟩) be the neutrosophic group ring of the neutrosophic group ⟨G ∪I⟩ over the ring R. Let P be a pseudo neutrosophic subring of the neutrosophic group ring R(⟨G ∪ I⟩) we say P is a pseudo neutrosophic ideal of the neutrosophic group ring R (⟨G ∪I⟩) if the following condition is satisfied.*

*For all p ∈ P and α ∈ R(⟨G ∪I⟩) pα and α p ∈ P.*

The reader is expected to construct more examples however we give an example of each.

***Example 3.1.12:*** Let R(⟨G ∪ I⟩) be the neutrosophic group ring of the neutrosophic group ⟨G ∪ I⟩ over the ring R, where R = {$Z_2$ = {0,1}} and

{G ∪ I} ={1, g, $g^2$, $g^3$, I gI, $g^2$I $g^3$I / $g^4$ = 1 and $I^2$ =I}.

Clearly P = ⟨0, $(1 + g + g^2 + g^3)$, $(I + gI + g^2I + g^3I)$, $(1 + g + g^2 + g^3 + I + gI + g^2I + g^3I)$⟩ is a pseudo neutrosophic ideal of the neutrosophic group ring R⟨G ∪ I⟩.

Here we give yet another example of a neutrosophic ideal.

***Example 3.1.13:*** Let Z⟨G ∪ I⟩ be a neutrosophic group ring of the neutrosophic group ⟨G ∪ I⟩ over the ring Z. ZP = 2Z (⟨G ∪ I⟩) is a neutrosophic subring is the neutrosophic ideal of the neutrosophic group ring.

***Note:*** Clearly 2Z ⟨G ∪ I⟩ is not a neutrosophic subring which is a neutrosophic group ring. Also a neutrosophic subring need not always contain the identity of the neutrosophic ring.

One can define as a matter of routine the notion of maximal neutrosophic ideal, minimal neutrosophic ideal, principal neutrosophic ideal and prime neutrosophic ideal. This is left as an exercise for the reader. Several problems in this direction are given in the last chapter.

Now we proceed onto define the notion of quasi neutrosophic ideals and loyal pseudo neutrosophic ideals of a neutrosophic group ring.

**DEFINITION 3.1.9:** *Let R(⟨G ∪ I⟩) be a neutrosophic group ring of the neutrosophic group ⟨G ∪ I⟩ over the ring R. Let S be any subring (neutrosophic or otherwise). Suppose their exists another subring P in R (⟨G ∪ I⟩) such that S is an ideal over P i.e., sp and ps ∈ S for all p ∈ P and s ∈ S, then we call S a quasi neutrosophic ideal of the neutrosophic group ring R(⟨G ∪ I⟩) relative to P.*

*If S happens to be only a right or left ideal only then we call S a quasi neutrosophic right (left) ideal relative to P.*

The following are interesting, so we just make a mention of them.

**DEFINITION 3.1.10:** *For a given or choosen S in the definition 3.1.9 one may have several $P_i$ such that S is a quasi neutrosophic ideal over $P_i$.*

*If in a neutrosophic group ring for a given S we have only one P such that S is a quasi neutrosophic ideal relative only to P and for no other P then we call such S to be a loyal quasi neutrosophic ideal relative to P.*

*If every subring S of a neutrosophic group ring happens to be a loyal quasi neutrosophic ideal relative to a unique P we call the neutrosophic group ring $R(\langle G \cup I \rangle)$ itself a loyal neutrosophic group ring.*

**DEFINITION 3.1.11:** *Let $R(\langle G \cup I \rangle)$ be a neutrosophic group ring of the neutrosophic group $\langle G \cup I \rangle$ over the ring R. If for S a subring P is another subring ($P \neq S$) such that S is a quasi neutrosophic ideal relative P. If in particular P happens to be quasi neutrosophic ideal relative to S then we call (P, S) a bonded quasi neutrosophic ideal of the neutrosophic group ring.*

***Note:*** We throw open several problems which are both interesting and needs more research.

1. Should a bonded quasi neutrosophic ideals be loyal?
2. Give examples of loyal quasi neutrosophic ideals which are not bonded!

As in case of usual ideals in group rings in case of bonded quasi neutrosophic ideals also we can define bonded quasi neutrosophic right ideals or bonded quasi neutrosophic left ideals if S in the definition happens to be either left (right) ideal. Similarly one can also define loyal quasi neutrosophic right or left ideals.

We have got a very interesting relation which is left as an exercise for the reader to prove.

Prove if S is a subring structure of a neutrosophic group ring and if S is quasi neutrosophic ideal over $P_1$ and $P_2$. Does it always imply $P_1 \subset P_2$ or $P_2 \subset P_1$? Or can $P_1$ and $P_2$ be disjoint?

Now we proceed onto define the notion of pseudo bonded quasi neutrosophic ideals and pseudo loyal quasi neutrosophic ideals of a neutrosophic group ring R(⟨G ∪ I⟩).

**DEFINITION 3.1.12:** *Let R(⟨G ∪ I⟩) be a neutrosophic group ring of the neutrosophic group over the ring R. Let S be a pseudo neutrosophic subring, if S is a quasi ideal over a subring P then we call S a pseudo neutrosophic ideal relative to the subring P.*

*P in general need not be pseudo neutrosophic, if in particular P also happens to be a pseudo neutrosophic subring of R(⟨G ∪ I⟩) then we call S a strong pseudo neutrosophic quasi ideal of R(⟨G ∪ I⟩).*

*It can be easily verified that in general all strong pseudo quasi neutrosophic ideals are pseudo quasi neutrosophic ideals and not conversely.*

Give an example of a pseudo quasi neutrosophic ideal which is not a strong pseudo quasi neutrosophic ideal.

Now we proceed onto define the notion of pseudo loyal quasi neutrosophic ideal and pseudo bonded quasi neutrosophic ideal.

**DEFINITION 3.1.13:** *Let R(⟨G ∪ I⟩) be a neutrosophic group ring of the neutrosophic group ⟨G ∪ I⟩ over the ring R.. Let S be a pseudo neutrosophic subring.*

*We say S is a quasi pseudo neutrosophic loyal ideal relative to P if S is a quasi pseudo neutrosophic ideal only relative to P and is not a quasi pseudo neutrosophic ideal relative to any other subring of R.*

*If in addition P also happens to be a pseudo neutrosophic subring then we call S a pseudo strong loyal quasi neutrosophic ideal.*

Clearly we see all pseudo strong loyal quasi neutrosophic ideals are pseudo loyal quasi neutrosophic ideals and not conversely.

Now one can define the notion of bonded pseudo quasi neutrosophic ideals and strongly bonded pseudo quasi neutrosophic ideals of a neutrosophic group ring $R(\langle G \cup I \rangle)$.

As in case of other ideals one can define the notion of left (and or) right pseudo bonded quasi neutrosophic ideals.

## 3.2 Some special properties of Neutrosophic Group Rings

Now we proceed onto define some special properties of neutrosophic group ring and element wise analysis of neutrosophic group rings and obtain some results analogous to group rings.

**DEFINITION 3.2.1**: *Let $R(\langle G \cup I \rangle)$ be a neutrosophic group ring of a neutrosophic group $\langle G \cup I \rangle$ over the ring R. We say an element $\alpha \in R\langle G \cup I \rangle$ is a neutrosophic idempotent if support of $\alpha$ has neutrosophic group elements and $\alpha^2 = \alpha$. If $\alpha \in R\langle G \cup I \rangle$ has no neutrosophic group elements then we just call $\alpha$ an idempotent if $\alpha^2 = \alpha$.*

We illustrate both these concepts by the following example.

**Example 3.2.1:** Let $R\langle G \cup I \rangle$ be the neutrosophic group ring where $R = Z_2$ and

$$\langle G \cup I \rangle = \{1, g, g^2, g^3, I\ gI, g^2I, g^3I\ /g^4 = 1 \text{ and } I^2 = I\}$$

$\alpha = 1 + g + g^2 + g^3 \in Z_2 \langle G \cup I \rangle$ is such that $\alpha^2 = \alpha$; $\alpha$ is an idempotent of the neutrosophic group ring. Take $(1 + I) = \beta \in Z_2 \langle G \cup I \rangle$; $(1 + I)^2 = 1 + I$, so $\beta$ is a neutrosophic idempotent of $Z_2\langle G \cup I \rangle$.

Now we proceed onto define units and neutrosophic units in a neutrosophic group ring.

**DEFINITION 3.2.2:** *Let $R(\langle G \cup I \rangle)$ be a neutrosophic group ring of a neutrosophic group $\langle G \cup I \rangle$ over a ring R. An element $\alpha$ is* *said to be unit in R($\langle G \cup I \rangle$) if their exists a β in R($\langle G \cup I \rangle$) such that αβ = 1. An element α' in the neutrosophic group ring R($\langle G \cup I \rangle$) is said to be a neutrosophic unit of R($\langle G \cup I \rangle$) if there exists a β' in R($\langle G \cup I \rangle$) such that α'β' = I.*

**Note**: As $I^2$ = I, I will be called as a trivial neutrosophic unit of the neutrosophic group ring. We also can call for α; β is the inverse in the definition.

Similarly for α' in R$\langle G \cup I \rangle$ we have β' to be the neutrosophic inverse of α'. It is further important to note α need not be an ordinary element. It can be a neutrosophic element.

To this end we give an example.

**Example 3.2.2:** Let $Z_4$ = {0, 1, 2, 3} be the ring of integers modulo 4 and G = $\langle g^4 = 1 \rangle$ and $\langle G \cup I \rangle$ = {1, g, $g^2$, $g^3$, I, gI, $g^2$I, $g^3$I/ $g^4$ = 1 and $I^2$ = I}. Let $Z_4$ $\langle G \cup I \rangle$ be the neutrosophic group ring. Clearly α = (1 + 2gI) is such that $α^2$ = 1 but α is a neutrosophic element but α is not a neutrosophic unit of $Z_4$ $\langle G \cup I \rangle$.

Next we proceed on to define the notion of zero divisors and neutrosophic zero divisors in a neutrosophic group ring R($\langle G \cup I \rangle$).

**DEFINITION 3.2.3:** *Let R($\langle G \cup I \rangle$) be the neutrosophic group ring of the neutrosophic group $\langle G \cup I \rangle$ over the ring R. An element of 0 ≠ α ∈ R($\langle G \cup I \rangle$) which is just a non neutrosophic element is called a zero divisor in R($\langle G \cup I \rangle$) if their exists a non zero β in R($\langle G \cup I \rangle$) with αβ = 0.*

*An element 0 ≠ α ∈ R($\langle G \cup I \rangle$) is a neutrosophic zero divisor if α is a neutrosophic element in R($\langle G \cup I \rangle$) and their exists a β ≠ 0 ∈ R($\langle G \cup I \rangle$) with αβ = 0.*

**Example 3.2.3:** Let R($\langle G \cup I \rangle$) be a neutrosophic group ring where R = $Z_2$ the ring of integers modulo 2 and $\langle G \cup I \rangle$ = {1, g, $g^2$, $g^3$, I, gI, $g^2$I, $g^3$I/$g^4$ = 1 and $I^2$ = I}. R $\langle G \cup I \rangle$ = $Z_2$ $\langle G \cup I \rangle$ is the neutrosophic group ring we see $Z_2$ $\langle G \cup I \rangle$ contains both zero divisors and neutrosophic zero divisors. For take $\alpha = (I + g)$ and $\beta = 1 + g + g^2 + g^3$ in $Z_2 \langle G \cup I \rangle$, $\alpha\beta = 0$ is a zero divisor. Take $\alpha = (I + g^2)$ and $\beta = (I + gI + g^2I + g^3I)$. Cleary $\alpha\beta = 0$. Thus $\alpha$ is a neutrosophic zero divisor in $Z_2 (\langle G \cup I \rangle)$.

We have the following interesting result about zero divisors and neutrosophic zero divisors in case of neutrosophic group rings.

Just we define neutrosophic torsion element of a group.

**DEFINITION 3.2.4:** *Let $\langle G \cup I \rangle$ be a neutrosophic group. An element g of $\langle G \cup I \rangle$ is said to be neutrosophic torsion element if $(g)^n = I$ for some finite integer n. If for no n we have $g^n = I$ where g is a neutrosophic element of $\langle G \cup I \rangle$ then we call g a neutrosophic torsion free element of $\langle G \cup I \rangle$.*

**Example 3.2.4:** Let $Q^+$ denote the set of positive rationals, $\langle Q^+ \cup I \rangle$ be the neutrosophic group under multiplication. $\langle Q^+ \cup I \rangle$ is a torsion free group as well as neutrosophic torsion free group.

**Example 3.2.5:** Let $\langle G \cup I \rangle = \{1, g, g^2, g^3, g^4, I, gI, g^2I, g^3I, g^4I$ / $g^5 = 1$ and $I^2 = I\}$ every element in $\langle G \cup I \rangle$ is either torsion or neutrosophic torsion.

Now we have the following result:

**LEMMA 3.2.1:** Let $\langle G \cup I \rangle$ be a neutrosophic group which is not torsion free and not neutrosophic torsion free then the neutrosophic group ring $K(\langle G \cup I \rangle)$ has proper zero divisors and neutrosophic zero divisors where K is a field.

*Proof:* By assumption $\langle G \cup I \rangle$ contains finite neutrosophic subgroup $\langle H \cup I \rangle \neq \langle 1 \cup I \rangle$. Set

$$\alpha = \sum_{x \in H \cup I} x \in (\langle G \cup I \rangle).$$

If $h \in H$ then $hH = H$ so $h\alpha = \alpha$.

Thus $\alpha^2 = n\alpha$ where n is the order of $\langle H \cup I \rangle$. So $\alpha(\alpha - n) = 0$, we see that $K\langle G \cup I \rangle$ has proper neutrosophic divisors of zero.

We illustrate the above lemma by the following example.

***Example 3.2.6:*** Let K be a field of characteristic zero and $\langle G \cup I \rangle = \{1, g, g^2, g^3, g^4, g^5, I, gI, g^2I, g^3I, g^4I, g^5I / g^6 = 1$ and $I^2 = I\}$. Let $K(\langle G \cup I \rangle)$ be a neutrosophic group ring of the neutrosophic group $(\langle G \cup I \rangle)$ over the field K. Clearly $\langle H \cup I \rangle = \{1, g^2, g^4, I, g^2I, g^4I\}$ and $|\langle H \cup I \rangle| = 6$.

Now if $\alpha = 1 + g^2 + g^4 + I + g^2I + g^4I$, then $\alpha(\alpha - 6) = 0$ which is a neutrosophic zero divisor of $K(\langle G \cup I \rangle)$. Let $\beta = 1 + g^2 + g^4$, $\beta [\beta - 3] = 0$ is a zero divisor of $K(\langle G \cup I \rangle)$.

Now we prove an interesting result about units and neutrosophic units of a neutrosophic group ring $R(\langle G \cup I \rangle)$.

**THEOREM 3.2.1:** *Let $\langle G \cup I \rangle$ be a neutrosophic group which is not torsion free and not neutrosophic torsion free and K any field with $|K| > 3$. Then the neutrosophic group ring KG has nontrivial units.*

***Proof:*** Given $\langle G \cup I \rangle$ is not a neutrosophic torsion free neutrosophic group. So G is not a torsion free group. K any field such that $|K| > 3$. To show KG has nontrivial units and nontrivial neutrosophic units. By the very assumption $\langle G \cup I \rangle$ contains neutrosophic subgroups. Let $\langle H \cup I \rangle$ be a neutrosophic subgroup. Set

$$\alpha = \sum_{x \in \langle H \cup I \rangle} x \in K (\langle G \cup I \rangle).$$

Let supp $\alpha = |\langle H \cup I \rangle| = n$ (say). Given $|K| > 3$ so we can choose $a \in K$ ($a \neq 0$ and 1 and a.n). Then $1 - a\alpha$ is a nontrivial unit in $K[\langle G \cup I \rangle]$ with inverse

$$b = \frac{a}{(an - 1)}$$

choose

$$b' = \frac{Ia}{\langle an - I \rangle}$$

then $I - b'\alpha$ is the neutrosophic inverse of $1 - a\alpha$.

Now we illustrate this theorem by the following example:

**Example 3.2.7:** Let $K = Q$, the field of rationals. Take $\langle G \cup I \rangle = \{1, g, g^2, g^3, g^4, g^5, I, gI, g^2I, g^3I, g^4I, g^5I \mid g^6 = 1$ and $I^2 = I\}$ as the neutrosophic group. $Q [\langle G \cup I \rangle]$ is the neutrosophic group ring of the neutrosophic group $\langle G \cup I \rangle$ over Q.

Take $\langle H \cup I \rangle = \{1, g^3, I, g^3I\}$. Let $\alpha = 1 + g^3 + I + g^3I \in Q [\langle G \cup I \rangle]$.

$$x = 1 - 5\alpha$$

is a nontrivial unit in $K [\langle G \cup I \rangle]$ with inverse.

$$y = \{1 - \frac{5}{19}(1 + g^3 + I + g^3I)\}$$

$$
\begin{aligned}
xy &= (1 - 5\alpha)\left[1 - \frac{5}{19}\alpha\right] \\[2mm]
&= 1 - 5\alpha - \frac{5\alpha}{19} + \frac{25 \times 4\alpha}{19} \\[2mm]
&= \frac{19 - 95\alpha - 5\alpha + 100\alpha}{19} \\[2mm]
&= 1.
\end{aligned}
$$

Thus the claim.

Take

$$
\begin{aligned}
x &= 1 - 5\alpha \text{ and} \\
y &= \left[I - \frac{5\alpha I}{19}\right] \\[2mm]
xy &= I,
\end{aligned}
$$

which shows x, y is a neutrosophic unit of $K [\langle G \cup I \rangle]$.

In case of neutrosophic group rings we have several interesting properties about them. There are also open problems analogous to the ones given in group rings.

Now just for the sake of completeness, we define the notion of two unique product neutrosophic group. Before we proceed on to define this we define a neutrosophic ordered group.

**DEFINITION 3.2.5:** *Let G be an ordered group (i.e., G admits a strict linear ordering < such that x < y implies xz < yz and zx < zy for all z ∈ G).*

*We say the neutrosophic group generated by this ordered group G and I i.e., ⟨G ∪ I⟩ is said to be neutrosophic ordered if xI < yI then (zI) xI < (zI) yI and (xI) zI < (yI) (zI) for all zI ∈ {⟨G ∪ I⟩} \ G.*

***Note:*** If x < y then we cannot in general say xI < yI it may so happen yI < xI. Thus the neutrosophic order in general need not preserve the order. If G is ordered under '<' the neutrosophic part of ⟨G ∪ I⟩ have the preservations of order under <; i.e., if x < y, x, y G then xI < yI may occur or may not occur.

Now we have the following interesting notion:

If G is an ordered group under '<' and if the neutrosophic part of the group ⟨G ∪ I⟩ is also ordered under the '<' then we call the neutrosophic group to be strongly ordered. In general all strongly ordered neutrosophic groups are ordered neutrosophic group and not conversely.

Also it is pertinent to note we do not compare an ordinary group element with a neutrosophic group element i.e., x and xI cannot be compared or ordered, also x and yI cannot be compared or ordered even in case of strongly ordered neutrosophic group.

Having defined the notion of ordered neutrosophic group now we proceed on to define the notion of t-u-p neutrosophic group and u-p neutrosophic group.

**DEFINITION 3.2.6:** *A group is said to be a u-p group (unique product group) if given any two nonempty finite subsets A and B of G then there exists atleast one element x ∈ G which has a unique representation in the form x = ab with a ∈ A and b ∈ B.*

*We say the neutrosophic group ⟨G ∪ I⟩ is said to be a two unique product neutrosophic group if given any two non empty finite subsets A and B of G, then there exists atleast one element x in G which has a unique representation in the form x = ab with a ∈ A and b ∈ B and if for any finite subsets A and B of ⟨G ∪ I⟩ there exist at least one element x in ⟨G ∪ I⟩ (x can be neutrosophic) which has a unique representation in the form x = ab with a ∈ A and b ∈ B.*

On similar lines we can define neutrosophic unique product group. Several interesting problems can be raised in this direction. The reader is expected to define strong unique product neutrosophic group and illustrate them with examples.

Now we prove some properties of neutrosophic modules in neutrosophic group rings. To this end we define the notion of neutrosophic modules.

**DEFINITION 3.2.7:** *A neutrosophic R-module V is a neutrosophic abelian group with law of composition written as "+" together with a scalar multiplication $R \times V \to V$ written as $r, v \mapsto rv$ where $r \in R$ (R is a commutative ring) and $v \in V$; satisfies these axioms.*

$$
\begin{aligned}
1.\, v &= v \\
r\,.I &= rI \\
1.\,I &= I \\
(rs)v &= r\,(sv) \\
(r+s)v &= rv + sv \\
r\,(v+v') &= rv + rv'
\end{aligned}
$$

*for all r, s ∈ R and v, v' ∈ V.*

***Note:*** Just in the definition of R-module we have replaced the abelian group by a neutrosophic abelian group; all other factors remain the same.

We first illustrate this by the following example:

***Example 3.2.8:*** Let R $[\langle G \cup I \rangle]$ be a neutrosophic group ring. R a commutative ring with unit, R $[\langle G \cup I \rangle]$ is a neutrosophic R-module as R $[\langle G \cup I \rangle]$ is an additive abelian group under usual addition '+'. i.e., R $\times$ R$[\langle G \cup I \rangle] \rightarrow$ R$[\langle G \cup I \rangle]$ satisfies all the axioms of the neutrosophic R-module.

Now we just make a mention of the definition of the neutrosophic submodule. Let V be an abelian neutrosophic group. R a commutative ring, V be a neutrosophic R-module. A proper-subset P of V is said to be a neutrosophic submodule of the R-module V if P is a nonempty set which is closed under addition and scalar multiplication. If the subset P of V is not a neutrosophic subgroup of V then P is not a neutrosophic submodule but only a module of V. It may still happen that the proper subset P of V is a pseudo neutrosophic subgroup in which case we call P the pseudo neutrosophic submodule.

Now we proceed on to illustrate them by an example.

***Example 3.2.9:*** Let Z be the ring of integers and $\langle G \cup I \rangle = \{1, g, g^2, g^3, g^4, g^5, I \, gI, g^2I, g^3I, g^4I, g^5I \, / \, g^6 = 1 \text{ and } I^2 = I\}$. Z $[\langle G \cup I \rangle]$ is a neutrosophic group ring. In fact Z $[\langle G \cup I \rangle]$ is also an additive abelian neutrosophic group and Z$[\langle G \cup I \rangle]$ is a neutrosophic Z module.

Z[G] is a subgroup of the neutrosophic group Z$[\langle G \cup I \rangle]$ and Z[G] is a submodule of Z$[\langle G \cup I \rangle]$. Z$[\langle H \cup I \rangle]$ where $\langle H \cup I \rangle = \{1, g^2, g^4, I \, g^2I, g^4I\}$ is a neutrosophic subgroup of Z$[\langle G \cup I \rangle]$ and Z$[\langle H \cup I \rangle]$ is a neutrosophic Z-submodule of Z$[\langle G \cup I \rangle]$.

Take P = $\{0, n \, (1 + g + g^2 + g^3 + g^4 + g^5 + I + gI + g^2I + g^3I + g^4I + g^5I) \, / \, n \in Z\}$. P is a pseudo neutrosophic Z submodule of Z$[\langle G \cup I \rangle]$.

It is interesting to note that all neutrosophic subgroups need not be neutrosophic submodules. Reader can construct examples for all these situations.

**Note:** To develop the notion of Jordan Holder theorem in case of neutrosophic K [⟨G ∪ I⟩] modules we need to develop the concept of full neutrosophic ring matrices over the neutrosophic field ⟨K ∪ I⟩. Thus at this stage we have only limited sources about the properties of neutrosophic rings.

Now we proceed on to define the notion of homomorphism of neutrosophic group rings.

**DEFINITION 3.2.8:** *Let R (⟨G ∪ I⟩) and S (⟨H ∪ I⟩) be any two neutrosophic group rings. A map φ from R (⟨G ∪ I⟩) to S (⟨H ∪ I⟩) is said to be a neutrosophic group ring homomorphism; if the following conditions are satisfied by φ :*

> 1. *φ (0) = 0*
> 2. *φ (I) = I*
> 3. *φ (1) = 1*
> 4. *φ (a + b) = φ (a) + φ (b) and*
> 5. *φ (ab) = φ (a). φ (b)*

*for all a, b, ∈ R (⟨G ∪ I⟩).*

The notion of neutrosophic isomorphism, neutrosophic automorphism are defined as in case of rings. One of the important notions to be noted is that kernel of any neutrosophic homomorphism is never a neutrosophic subring as φ (I) = I.

The kernel of any neutrosophic homomorphism of any neutrosophic group ring can only be either 0 or a subring and never a neutrosophic subring. Thus the kernel can never be a neutrosophic ideal. This is the marked difference between usual homomorphism of rings and homomorphism of neutrosophic rings. Thus the notion of quotient neutrosophic ring using kernel (φ a neutrosophic ring homomorphism) never arises in case of neutrosophic ring homomorphism.

Now as in case of group rings and rings we can define the notion of neutrosophic prime rings. In the first place it is important to note that a neutrosophic ring can never have ideals which are not neutrosophic i.e., all ideals in a neutrosophic group ring are neutrosophic ideals so we can only define a neutrosophic prime ring.

**DEFINITION 3.2.9:** *Let $R (\langle G \cup I \rangle)$ be a neutrosophic group ring. A neutrosophic ideal $P$ of $R (\langle G \cup I \rangle)$ is said to be neutrosophic prime if $R / P$ is a neutrosophic prime ring. (We say a neutrosophic group ring is neutrosophic prime if the neutrosophic group ring $R (\langle G \cup I \rangle)$ have neutrosophic ideals $A$, $B$ in $R (\langle G \cup I \rangle)$ with $AB = 0$ implies $A = 0$ or $B = 0$).*

*This $P$ is a neutrosophic prime ideal of $R (\langle G \cup I \rangle)$ if and only if for all neutrosophic ideals $A$, $B \subseteq R$ we have $AB \subset P$ implies $A \subseteq P$ or $B \subseteq P$. We say $R (\langle G \cup I \rangle)$ is neutrosophic prime if the intersection of all prime ideals of $R(\langle G \cup I \rangle)$ is generated by $\langle 0, I \rangle$ or just $\{0, I\}$.*

**Example 3.2.10:** Let $Z_2 [\langle G \cup I \rangle]$ be the neutrosophic group ring where $(\langle G \cup I \rangle) = \{1, g, I, g I / g^2 = 1$ and $I^2 = I\}$. $Z_2 (\langle G \cup I \rangle)$ is not a neutrosophic semiprime ring for P = $\{(1 + g), (I + gI), 1 + g + gI + I, 0\}$. Clearly $P^2 = 0$ and P is a neutrosophic ideal of $Z_2(\langle G \cup I \rangle)$. So $Z_2(\langle G \cup I \rangle)$ is a neutrosophic ring which is not semiprime.

Now we proceed on to find means to give any non trivial class of semiprime neutrosophic group rings. Thus a neutrosophic group ring R $(\langle G \cup I \rangle)$ (ring) is semiprime if and only if R $(\langle G \cup I \rangle)$ contains no nonzero ideal with square zero.

Using this definition we have the following nice result.

**THEOREM 3.2.2:** *Let K be a field of characteristic zero. Then $K(\langle G \cup I \rangle)$, the neutrosophic group ring of the neutrosophic group $(\langle G \cup I \rangle)$ over K is semiprime.*

*Proof:* Suppose the neutrosophic group ring K ($\langle G \cup I \rangle$) is not semiprime; then we have the neutrosophic group ring that contains a neutrosophic nonzero-ideal A with $A^2 = 0$. Let

$$\alpha = \sum_{i=1}^{n} k_i x_i \in A, \alpha \neq 0$$

and let F be the subfield of K generated over rationals by $k_1$, $k_2$, …, $k_n$.

Then F [$\langle G \cup I \rangle$] $\subset$ K [$\langle G \cup I \rangle$] and A $\cap$ F [$\langle G \cup I \rangle$] is a nonzero neutrosophic ideal in F[$\langle G \cup I \rangle$] of square zero. Thus it clearly suffices to assume that K = F or equivalently that K is finitely generated over the rationals. This implies that K is contained in the complex numbers C and we fix an imbedding.

Then K[$\langle G \cup I \rangle$] $\subseteq$ C[$\langle G \cup I \rangle$] and A is a non zero neutrosophic ideal of C[$\langle G \cup I \rangle$] with square zero. This contradicts the Jordan theorem, which says "suppose that K is a subfield of the complex numbers which is closed under complex conjugation, then the neutrosophic group ring K [$\langle G \cup I \rangle$] contains no non zero nil ideal.

Proving the results of *Passman* and *Connell* are direct as in case of group rings.

**THEOREM 3.2.3:** *Let K be a field of characteristic p, p > 0 and let the neutrosophic group $\langle G \cup I \rangle$ have no elements of order p. (if G has no elements of order p clearly the neutrosophic group $\langle G \cup I \rangle$ will not have elements of order p). Then the neutrosophic group ring K [$\langle G \cup I \rangle$] has no non zero neutrosophic nil ideals.*

$$\Delta^I = \Delta(\langle G \cup I \rangle) = \left\{ x \in \langle G \cup I \rangle / \left[ \langle G \cup I \rangle : c_{\langle G \cup I \rangle}^{(X)} \right] \prec \infty \right\}.$$

The following result can be easily proved by the reader analogous to the results of group rings with appropriate modifications.

**THEOREM 3.2.4:** *Let K be a field of characteristic p, p > 0. The following are equivalent:*

1. *K [⟨G ∪I⟩] is semiprime.*
2. *Δ(⟨G ∪I⟩) has no elements of order p.*
3. *⟨G ∪I⟩ has no finite normal subgroup or finite normal neutrosophic subgroup with order divisible by p.*

Let us now define the center of a neutrosophic group ring K [⟨G ∪ I⟩].

**DEFINITION 3.2.10:** *Let K [⟨G ∪I⟩] be a neutrosophic group ring of the neutrosophic group ⟨G ∪I⟩ over K, K any field. The neutrosophic center of the neutrosophic group ring K [⟨G ∪I⟩] denoted by*

$$C^I = \left\{ \alpha \in K[\langle G \cup I \rangle] / \alpha\beta = \beta\alpha \text{ for all } \beta \in K[\langle G \cup I \rangle] \right\}.$$

*Clearly $C^I \neq \phi$ for I and 1 ∈ $C^I$. Further if C is the center of the group ring KG ⊂ K [⟨G ∪I⟩] then C ⊆ $C^I$.*

We have the following interesting result:

**LEMMA 3.2.2:** Let K [⟨G ∪ I⟩] be a semiprime neutrosophic group ring with neutrosophic center $C^I$. If α ∈ $C^I$ then α is a zero divisor in K [⟨G ∪ I⟩] if and only if it is a zero divisor in $C^I$.

*Proof:* Given K[⟨G ∪ I⟩] is a semiprime neutrosophic group ring. Let α be a zero divisor neutrosophic or other wise in $C^I$ (As C ⊂ $C^I$ we may have $C^I$ containing both neutrosophic zero divisors and zero divisors). If α is a zero divisor in $C^I$ then certainly it is a zero divisor in K [⟨G ∪ I⟩].

We have to prove only the converse i.e., if α is a zero divisor (or neutrosophic zero divisor) in K [⟨G ∪ I⟩] to show α is a zero divisor in $C^I$. Let β ∈ K [⟨G ∪ I⟩], β ≠ 0 with αβ = 0. Let ⟨H ∪ I⟩ be the neutrosophic subgroup of ⟨G ∪ I⟩ generated by support of $\alpha$. Since $\alpha \in C^I$, we conclude that $\langle H \cup I \rangle$ is a finitely generated normal neutrosophic subgroup of $\langle G \cup I \rangle$ with $H \subseteq \Delta (\langle G \cup I \rangle)$.

Write $\beta = \sum \beta_i x_i$ with $\beta_i \in K [\langle G \cup I \rangle]$ and with the $x_i$ in cosets of $\langle H \cup I \rangle$ in $\langle G \cup I \rangle$. It follows $\alpha \beta_i = 0$ for all i.

Let $T = \{\gamma \in K [\langle H \cup I \rangle] / \alpha \gamma = 0\}$ since $\beta_i \in T$ for all i and $\beta \neq 0$, we conclude that $T \neq 0$. Clearly T is a neutrosophic right ideal of $K [\langle H \cup I \rangle]$. Moreover since $\alpha \in C^I$, i.e., $\alpha$ is central in $K [\langle G \cup I \rangle]$ and $\langle H \cup I \rangle$ is normal in $\langle G \cup I \rangle$ we have for $x \in \langle G \cup I \rangle$, $x^{-1} T x \subset K [\langle H \cup I \rangle]$ and $\alpha x^{-1} T x = x^{-1} (\alpha T) x = 0$ (if $x = g T$, $x^{-1}$ is then taken as $g^{-1}T$.)

Thus $x^{-1}Tx \subseteq T$. This yields $x^{-1} T x = T$. Thus $\alpha$ is a zero divisor in $C^I$ and the result follows, As neutrosophic group rings form a special class of neutrosophic rings we have introduced them and have analyzed to some extent.

However we have given several examples and many problems in the last chapter for the reader.

Now we proceed on to define other types of neutrosophic rings using semigroups and Smarandache semigroups (S-semigroups). As these all are illustrations of neutrosophic rings we do not wish to dwell deep into them.

### 3.3 Neutrosophic Semigroup Rings and their Generalizations

In this section we proceed on to define two types of neutrosophic rings viz. Neutrosophic semigroup rings and neutrosophic Smarandache semigroup rings and just give some of their properties.

**DEFINITION 3.3.1:** *Let R be any commutative ring with unit or a field, $\langle S \cup I \rangle$ a neutrosophic semigroup. The neutrosophic semigroup ring R $\langle S \cup I\rangle$ of the neutrosophic semigroup $\langle S \cup I\rangle$ over the ring R is defined analogous to neutrosophic group ring* *just defined in the earlier section. R ⟨S ∪ I⟩ consists of all finite formal sums of the form $\alpha = \sum_{i=1}^{n} r_i s_i$ ; $n < \infty$, $r_i \in R$ and $s_i \in ⟨S ∪ I⟩$ satisfying the following conditions.*

1. *Let $\alpha, \beta \in R(⟨S ∪ I⟩)$ where*

$$\alpha = \sum_{i=1}^{n} r_i s_i \ \ and \ \ \beta = \sum_{i=1}^{m} t_i s_i \ ;$$

   *$\alpha = \beta$ if and only if $n = m$ and $r_i = t_i$ for each i.*

2. *Let $\alpha, \beta \in R[⟨S ∪ I⟩]$,*

$$\alpha + \beta = \sum_{i=1}^{K} (r_i + t_i) \, s_i$$

   *where*

$$\alpha = \sum_{i=1}^{n} r_i s_i \ \ and \ \ \beta = \sum_{i=1}^{m} t_i s_i \ .$$

   *Clearly $\alpha + \beta \in R[⟨S ∪ I⟩]$.*

3. *$0 = 0s_1 + \ldots + 0s_n = \sum_{i=1}^{n} 0 s_i$ serves as the additive identity in $R[⟨S ∪ I⟩]$.*

4. *$\alpha + 0 = 0 + \alpha = \alpha$ for all $\alpha \in R[⟨S ∪ I⟩]$.*

5. *For every*

$$\alpha = \sum_{i=1}^{n} r_i s_i$$

   *we have a uniquely defined*

$$\beta = -\alpha = \sum_{i=1}^{n} -r_i s_i$$

   *such that*

$$\alpha + \beta = \alpha + (-\alpha) = \sum_{i=1}^{n} \left[ r_i + (-r_i) \right] s_i = 0.$$

   *$-\alpha$ is called the inverse of $\alpha$ with respect to +. Thus $R[⟨S ∪ I⟩]$ is an additive abelian group.*

6. *For $\alpha$, $\beta \in R\,[\langle S \cup I \rangle]$, where*

$$\alpha = \sum_{i=1}^{n} r_i s_i \ \ and \ \ \beta = \sum_{j=1}^{m} t_j s_j$$

$$\alpha . \beta = \left( \sum_{i=1}^{n} r_i s_i \right) \left( \sum_{j=1}^{m} t_j s_j \right) = \sum_{k=1}^{\delta} t_k s_k$$

*where $s_k = s_i\, s_j$ and $t_k = \displaystyle\sum_{s_i s_j = k} r_i t_j$ .*

*Since R is a distributive structure we see $R\,[\langle G \cup I \rangle]$ is also a distributive structure; i.e., $\alpha\,[\beta + \gamma] = \alpha\beta + \alpha\gamma$ and $[\alpha + \beta]\,\gamma = \alpha\gamma + \beta\gamma$ for all $\alpha$, $\beta$, $\gamma \in R\,[\langle S \cup I \rangle]$.*

Thus $R[\langle S \cup I \rangle]$ is a neutrosophic semigroup under multiplication. So $\{R[\langle S \cup I \rangle], +, .\}$, is a neutrosophic ring called the neutrosophic semigroup ring of the neutrosophic semigroup $\langle S \cup I \rangle$ over the ring R.

Clearly $\langle S \cup I \rangle \subseteq R\,[\langle S \cup I \rangle]$ as $1 \in R$; only if $\langle S \cup I \rangle$ is a neutrosophic monoid we may have $R \subseteq R\,[\langle S \cup I \rangle]$.

We just give an example of a neutrosophic semigroup ring.

**Example 3.3.1:** Let Q be the field of rationals, $[\langle S \cup I \rangle]$ be a neutrosophic semigroup given by the following table:

| 0  | 0 | 1  | h  | k  | I  | hI | kI |
|----|---|----|----|----|----|----|----|
| 0  | 0 | 0  | 0  | 0  | 0  | 0  | 0  |
| 1  | 0 | 1  | h  | k  | I  | hI | kI |
| h  | 0 | h  | 0  | h  | hI | 0  | hI |
| k  | 0 | k  | h  | 1  | kI | hI | I  |
| I  | 0 | I  | hI | kI | I  | hI | kI |
| hI | 0 | hI | 0  | hI | hI | 0  | hI |
| kI | 0 | kI | hI | I  | kI | hI | I  |

$$\langle S \cup I \rangle \quad = \quad \{0, 1, h, k, I, hI, kI\}$$
$$\cong \quad \{0, 1, 2, 3, I, 2I, 3I\}$$
$$= \quad \langle Z_4 \cup I \rangle\,;$$

the neutrosophic semigroup under multiplication modulo 4. Now Q [⟨S ∪ I⟩] is the neutrosophic semigroup ring.

Now we see this neutrosophic semigroup ring has zero divisors and nontrivial idempotents as the neutrosophic semigroup has nontrivial zero divisors and idempotents. For take $\alpha = 3.h$ and $\beta = 8hI$, $\alpha\beta = 0$ as $h^2 = 0$.

As in case of group rings one can study all the properties of neutrosophic semigroup rings. This forms yet another nice class of neutrosophic rings which has a very rich algebraic generalized structure.

Now we can have yet another class of neutrosophic rings which we can define as semigroup neutrosophic rings. Here we take basically any semigroup and a neutrosophic ring and build the semigroup neutrosophic ring. Here also we mainly take only commutative neutrosophic ring with unit and not pseudo neutrosophic rings.

Thus throughout this section we assume or take only commutative neutrosophic rings with unit, which has been defined in the earlier chapter.

**DEFINITION 3.3.2:** *Let S be any semigroup, ⟨R ∪ I⟩ be any commutative neutrosophic ring with unit. ⟨R ∪ I⟩ [S] is defined as the semigroup neutrosophic ring which consists of all finite formal sums of the form $\sum_{i=1}^{n} r_i s_i$ ; $n < \infty$, $r_i \in \langle R \cup I \rangle$ and $s_i \in S$.*

*This semigroup neutrosophic ring is defined analogous to the group ring or semigroup ring.*

***Note:*** Here only the ring is taken as the commutative neutrosophic ring with unit. For in all the earlier cases we have taken only the commutative ring with unit or a field, only now we shift to semigroups over neutrosophic rings.

We just illustrate with some examples.

***Example 3.3.2:*** Let $\langle Z_6 \cup I \rangle$ be the neutrosophic ring S the semigroup S(3). $\langle Z_6 \cup I \rangle$ [S(3)] is a semigroup neutrosophic ring. Clearly $[\langle Z_6 \cup I \rangle]$ [S(3)] is a non commutative neutrosophic ring.

We will see yet another example which is a commutative neutrosophic ring.

***Example 3.3.3:*** Let $\langle Z_2 \cup I \rangle = \{0, 1, I, 1 + I\}$ be the neutrosophic ring of characteristic 2 and S = $\{0, 1, g, h, k, e\}$ under multiplication given by the following table:

|   | 0 | 1 | g | h | k | e |
|---|---|---|---|---|---|---|
| 0 | 0 | 0 | 0 | 0 | 0 | 0 |
| 1 | 0 | 1 | g | h | k | e |
| g | 0 | g | k | 0 | g | k |
| h | 0 | h | 0 | h | 0 | h |
| k | 0 | k | g | 0 | k | g |
| e | 0 | e | k | h | g | 1 |

Clearly S = the semigroup isomorphic with $Z_6 = \{0, 1, 2, …, 5\}$ under multiplication modulo 6.

Further S is a monoid with zero divisors and idempotents. $[\langle Z_2 \cup I \rangle][S]$ is the semigroup neutrosophic ring of finite order which is commutative. Thus by this method we can find several other neutrosophic rings.

$$(1 + h)^2 \quad = \quad 1 + h^2 + 2h$$
$$= \quad 1 + h,$$

is the nontrivial idempotent.

$$(I + hI)^2 \quad = \quad 1 + hI;$$

is a nontrivial neutrosophic idempotent of $[Z_2 \cup I]$ [S].

$$(1 + g + k)\ (1 + g) \quad = \quad (1 + g + k + g + k + g) = 1 + g,$$

is not an idempotent or zero divisor. Take

$$(1 + g + k)\ (1 + g + k) \quad = \quad 1 + g + k + g + k + g + k + g + k$$
$$= \quad 1.$$

Thus $(1 + g + k)$ is a unit.

Now $(1 + g + k)\ (I + gI + kI) = I$ is a neutrosophic unit.

Thus we see this neutrosophic ring has neutrosophic units, units, idempotents, neutrosophic idempotents, zero divisors and neutrosophic zero divisors.

We can have several such interesting examples of neutrosophic rings. On similar lines we can take a commutative neutrosophic ring with unit and any group G and form the group neutrosophic ring which we will first define and then illustrate it with examples.

**DEFINITION 3.3.3:** *Let $\langle R \cup I \rangle$ be a commutative neutrosophic ring with unit and G any group. The group neutrosophic ring of the group G over the neutrosophic ring $\langle R \cup I \rangle$ denoted by $\langle R \cup I \rangle$ [G] consists of all finite formal sums of the form*

$$\sum_{i=1}^{n} r_i g_i, n < \infty; r_i \in \langle R \cup I \rangle$$

*and $g_i \in G$. Take*

$$\alpha = \sum_{i=1}^{n} r_i g_i \ \ and \ \ \beta = \sum_{i=1}^{n} s_i g_i$$

*$\alpha$, $\beta \in \langle R \cup I \rangle$ [G]; $\alpha = \beta$ if and only if $r_i = s_i$ and $m = n$ for each i. For $\alpha$, $\beta \in [\langle R \cup I \rangle]$ [G].*

$$\alpha + \beta = \sum_{i} (r_i + s_i) g_i \ \in \langle R \cup I \rangle \ [G].$$

*$0 = \sum_{i=1}^{n} 0 g_i$ serves as the additive identity, for each*

$$\alpha = \sum_{i=1}^{n} r_i g_i$$

*take*

$$\beta = \sum_{i=1}^{n} -\alpha_i g_i = (-\alpha),$$

*then*

$$\alpha + \beta = \alpha + (-\alpha) = 0.$$

*For $\alpha$, $\beta \in \langle R \cup I \rangle$ [G] define*

$$\alpha . \beta = \left( \sum_{i=1}^{n} r_i g_i \right) \left( \sum_{j=1}^{m} s_j g_j \right)$$

$$= \sum_k t_k g_k, \quad g_k = g_i g_j \quad and \quad t_k = \sum r_i s_j.$$

*Clearly $\alpha . \beta \in \langle R \cup I \rangle$ [G].*

*Thus $\langle R \cup I \rangle$ [G] is the group neutrosophic ring of the group G over the neutrosophic ring $\langle R \cup I \rangle$. Clearly $\langle R \cup I \rangle$ [G] is a neutrosophic ring, $\langle R \cup I \rangle$ [G] is commutative if and only if G is a commutative group.*

First we give some examples of them.

**Example 3.3.4:** Let $\langle Z_{12} \cup I \rangle$ be the neutrosophic ring. G the cyclic group of order 12, i.e., G = $\langle g / g^{12} = 1 \rangle$, $\langle Z_{12} \cup I \rangle$[G] is the group neutrosophic ring of the group G over the neutrosophic ring $\langle Z_{12} \cup I \rangle$. Let $\alpha = 6 + 6g^6$. $\alpha^2 = 36 + 36 + 2.36$ $g^6 = 0$ (mod 12). Thus $\alpha$ is a zero divisor of the group neutrosophic ring. Take $\alpha = 6 + 6g^6$ and $\beta = 6I + 6g$ in $\langle Z_{12} \cup I \rangle$ [G]. Clearly $\alpha\beta = 0$ is a neutrosophic zero divisor of $\langle Z_{12} \cup I \rangle$ [G]. Take $\alpha = 6 + 6g$ and $\beta = 2 + 6I + 4g^2 + 4g^3I$ in $\langle Z_{12} \cup I \rangle$ [G]. $\alpha\beta = 0$ is a neutrosophic zero divisor of $\langle Z_{12} \cup I \rangle$ [G].

We give yet another example.

**Example 3.3.5:** Let $\langle Q \cup I \rangle$ be the neutrosophic ring of integers, $S_3$ be the symmetric group of degree 3. $\langle Q \cup I \rangle$ [$S_3$] is the group neutrosophic ring. Clearly $\langle Q \cup I \rangle$ [$S_3$] is a non commutative group neutrosophic ring. Consider $(1 - p_1)$ $(1 + p_1 + p_2 + p_3 + p_4 + p_5) = 0$ thus the group neutrosophic ring has nontrivial divisors of zero. Take $\alpha = 1 - p_1$ and $\beta = 1 + p_1I + p_2 + p_3 + p_4 + p_5$. Cleary $\alpha\beta = 0$. This is a nontrivial neutrosophic divisor of zero.

Consider $\left[\dfrac{1}{6}\left(1 + p_1 + p_2 + p_3 + p_4 + p_5\right)\right]^2$

$$= \dfrac{1}{6}\left(1 + p_1 + p_2 + p_3 + p_4 + p_5\right),$$

so $\langle Q \cup I \rangle\,[S_3]$ has nontrivial idempotents.

If we take $\alpha = \dfrac{1}{3}\left(I + p_4 I + p_5 I\right)$ in $\langle Q \cup I \rangle\,[S_3]$ we have $\alpha^2 = \alpha$. Thus $\langle Q \cup I \rangle\,[S_3]$ has also non trivial neutrosophic idempotents. Let

$\alpha \quad = \quad 1 + p_4 + p_5,\ 1 - 5\,\alpha = x$ and

$y \quad = \quad 1 - \dfrac{5\alpha}{14},\ x, y \in \langle Q \cup I \rangle[S_3]$.

$x \cdot y \quad = \quad \left(1 - 5\alpha\right)\left(1 - \dfrac{5\alpha}{14}\right)$

$\qquad\quad = \quad 1 - 5\alpha + \dfrac{25 \times 3\alpha}{14} - \dfrac{5\alpha}{14}$

$\qquad\quad = \quad 1 - \dfrac{70\alpha + 75\alpha - 5\alpha}{14} = 1 - 0\ .$

Thus $\langle Q \cup I \rangle\,[S_3]$ has nontrivial units. In fact $\langle Q \cup I \rangle\,[S_3]$ has non trivial neutrosophic units also.

For if $\beta = 1 + p_1$, $x = 1 - 3\beta$ and $y = I - \dfrac{3\beta I}{5}$. Clearly $xy = I$.

Several interesting problems in this direction can be defined and innumerable examples can be given. However as our aim is to introduce neutrosophic ring and show the existence of such new classes of neutrosophic rings we do not go deep into group neutrosophic rings but propose several problems analogous to group rings in the last chapter of this book.

Next we define the notion of Smarandache semigroup neutrosophic rings. Throughout this section when we say Smarandache semigroup S we mean S is a semigroup and has a proper subset A where A is a group under the operations of the semigroup S.

Now we see when we replace the semigroup in the definition of semigroup neutrosophic ring by a Smarandache semigroup (S-semigroup) what are the properties we can enumerate specially and distinctly using the fact the S-semigroup has a proper subset which is a group.

**DEFINITION 3.3.4:** *Let $\langle K \cup I \rangle$ be commutative neutrosophic ring with unit. S be a Smarandache semigroup. The S-semigroup neutrosophic ring $\langle K \cup I \rangle [S]$ is defined analogous to group neutrosophic rings, the groups will be replaced by S semigroups. Thus $\langle K \cup I \rangle [S]$ will be a neutrosophic ring and it will contain all finite formal sums of the form*

$$\sum_{i=1}^{n} \alpha_i s_i \, , \; \alpha_t \in \langle K \cup I \rangle \text{ and } s_i \in S \text{ with } n < \infty.$$

*The S-semigroup neutrosophic ring will be commutative if and only if the S-semigroup is commutative.*

Now we proceed on to prove a nice result which show S-semigroup neutrosophic rings can be realized as the generalization of group neutrosophic rings.

**THEOREM 3.3.1:** *Let S be a S semigroup, $\langle K \cup I \rangle$ be a neutrosophic ring. $\langle K \cup I \rangle [S]$ be the S-semigroup neutrosophic ring of the S-semigroup S over the neutrosophic ring $\langle K \cup I \rangle$. $\langle K \cup I \rangle [S]$ always contains a proper subset which is a group neutrosophic ring.*

*Proof:* Given $\langle K \cup I \rangle [S]$ is a S-semigroup neutrosophic ring where S is a S-semigroup. The fact S is a S-semigroup guarantees of a proper subset A of S which is a group under the operations of S. Thus $\langle K \cup I \rangle [A]$ is a group neutrosophic ring of the group A over the neutrosophic ring $\langle K \cup I \rangle$. Clearly $\langle K \cup I \rangle [A]$ is a proper subset of $\langle K \cup I \rangle [S]$. Hence the claim.

Now we proceed on to give an example of a S-semigroup neutrosophic ring.

***Example 3.3.6:*** Let $\langle Z_7 \cup I \rangle$ be any neutrosophic ring. S = S(3) is the S-semigroup. $\langle Z_7 \cup I \rangle$ [S(3)] is the S-semigroup neutrosophic ring. Take T = $\langle Z_7 \cup I \rangle$ [S$_3$], S$_3 \subseteq$ S(3) and S$_3$ is the symmetric group of degree 3. Clearly T is the group neutrosophic ring of the S-semigroup neutrosophic ring.

***Note:*** A S-semigroup neutrosophic ring can have more than one group neutrosophic ring.

Another thing of interest is that both the group neutrosophic ring as well S-semigroup neutrosophic ring contains a proper subset which is a group ring. Thus it is easy to note that the notions of S-semigroup neutrosophic ring and group neutrosophic rings are still generalization of group rings.

**THEOREM 3.3.2:** *Every S-semigroup neutrosophic ring contains a proper subset which is a group ring.*

*Proof:* Given $\langle K \cup I \rangle$ [S] is a S – semigroup neutrosophic ring of the S-semigroup S over the neutrosophic ring $\langle K \cup I \rangle$. Consider the proper subset A of S, where A is a group. Clearly $\langle K \cup I \rangle$ [A] is a group neutrosophic ring of the group A over the neutrosophic ring $\langle K \cup I \rangle$. Consider K $\subseteq \langle K \cup I \rangle$, K is a commutative ring with unit or a field. Take KA, KA is the group ring of the group A over the ring K. Thus KA $\subseteq \langle K \cup I \rangle$ [A] $\subseteq \langle K \cup I \rangle$ [S]. Hence the claim.

**COROLLARY:** Every group neutrosophic ring $\langle K \cup I \rangle$ [G] contains a proper subset which is a group ring.

*Proof:* Follows from the fact that $\langle K \cup I \rangle$ [S] has a group neutrosophic ring which has a proper subset which is a group ring. Hence the claim.

We just illustrate this situation before we proceed on to prove more result.

***Example 3.3.7:*** Consider the S-semigroup neutrosophic ring $\langle Z_{12} \cup I \rangle$ [S(8)]. Clearly the set $Z_{12}$ S$_8$ is the group ring contained in the S-semigroup neutrosophic ring $\langle Z_{12} \cup I \rangle$ [S(8)] as $Z_{12} \subseteq \langle Z_{12} \cup I \rangle$ and $S_8 \subseteq S(8)$ is the symmetric group of degree 8. Not only this $\langle Z_{12} \cup I \rangle$ [S(8)] contains several such group rings. For take

$$P = \left\{ \begin{pmatrix} 1 & 2 & 3 & ... & 8 \\ 1 & 2 & 3 & ... & 8 \end{pmatrix}, \begin{pmatrix} 1 & 2 & 3 & 4 & 5 & 6 & 7 & 8 \\ 2 & 1 & 3 & 4 & 5 & 6 & 7 & 8 \end{pmatrix} \right\};$$

P is a group of order 2 and $Z_{12}$ P is a group ring. Take T to be the group generated by

$$g = \begin{pmatrix} 1 & 2 & 3 & 4 & 5 & 6 & 7 & 8 \\ 2 & 3 & 4 & 5 & 6 & 7 & 8 & 1 \end{pmatrix}.$$

Clearly $T = \langle g \rangle$ and $Z_{12}$ T is a group ring of the cyclic group of order 8 over the ring of integers $Z_{12}$.

Now we see this S-semigroup neutrosophic ring $\langle Z_{12} \cup I \rangle$ (S(8)) contains nontrivial zero divisors and neutrosophic zero divisors, idempotents, neutrosophic idempotents, units and neutrosophic units.

As our motivation is not to dwell deep into the concept of S-semigroup neutrosophic rings but only show there exists several classes of neutrosophic rings and such generalized rich class happens to be the S-semigroup neutrosophic rings.

We prove just a condition for the S-semigroup neutrosophic ring to have proper zero divisors and non trivial units.

**THEOREM 3.3.3:** Let S be a finite S-semigroup and $\langle K \cup I \rangle$ any neutrosophic ring. Then $\langle K \cup I \rangle$ [S] has proper divisors of zero and neutrosophic divisors of zero. Further if $|K| > 3$, $[\langle K \cup I \rangle]$ [S] has non trivial units and neutrosophic units.

*Proof:* Given S is a finite S-semigroup, so S contains a finite proper subset A which is a group. Let $|A| = n$.
Take $\alpha = \sum_{x \in A} x \in \langle K \cup I \rangle[S]$.

Now $\alpha^2 = n \alpha$ so $\alpha (\alpha - n) = 0$. Thus $\langle K \cup I \rangle [S]$ has proper divisors of zero.

If we take $\alpha = \sum_{x \in A} Ix$ then $\alpha^2 = n I \alpha$ so $\alpha(\alpha - nI) = 0$. Thus $\langle K \cup I \rangle [S]$ has non trivial neutrosophic zero divisors. Consider $d \in K$ with $d \neq 0$ or $1$ once $|K| > 3$, $x = 1 - d\alpha$ is a nontrivial unit in $\langle K \cup I \rangle [S]$ with inverse $y = 1 - \dfrac{d\alpha}{(dn - 1)}$ i.e., $xy = 1$.

Hence the claim.

Now we illustrate this theorem by the following example:

**Example 3.3.8:** Let $\langle Q \cup I \rangle [S(5)]$ be the S-semigroup neutrosophic ring. Let

$$H = \left\{ \left\langle \begin{pmatrix} 1 & 2 & 3 & 4 & 5 \\ 2 & 3 & 4 & 5 & 1 \end{pmatrix} \right\rangle \right\}$$

be the group generated by

$$\begin{pmatrix} 1 & 2 & 3 & 4 & 5 \\ 2 & 3 & 4 & 5 & 1 \end{pmatrix}.$$

Let

$$\begin{pmatrix} 1 & 2 & 3 & 4 & 5 \\ 2 & 3 & 4 & 5 & 1 \end{pmatrix} = g$$

then

$$H = \left\{ e = \begin{pmatrix} 1 & 2 & 3 & 4 & 5 \\ 2 & 3 & 4 & 5 & 1 \end{pmatrix} = 1, g, g^2, g^3, g^4, g^5 = e = 1 \right\}.$$

Take $\alpha = 1 + g + g^2 + g^3 + g^4$, $\alpha^2 = 5\alpha$ so $\alpha (\alpha - 5) = 0$. Clearly $\alpha \neq 0$ and $\alpha \neq 5$. Hence $\alpha$ is a nontrivial divisor of zero.

Take

$$x = 1 - 3\alpha$$

and

$$y = 1 - \dfrac{3\alpha}{14}.$$

Consider

$$xy = (1-3\alpha)\left(1-\frac{3\alpha}{14}\right)$$

$$= 1-3\alpha-\frac{3\alpha}{14}+\frac{45\alpha}{14}$$

$$= 1-\frac{42\alpha}{14}-\frac{3\alpha}{14}+\frac{45\alpha}{14}=1.$$

Thus $\langle Q \cup I\rangle$ [S(5)] has non trivial divisors of zero and nontrivial units. Take

$\alpha_1 = I. \ 1 + Ig + Ig^2 + Ig^3 + Ig^5,$
$\alpha^2_1 = 5\,\alpha_1$ and $\alpha_1\,(\alpha_1 - 5) = 0.$

Thus $\alpha_1$ is a neutrosophic zero divisor.

Several such results can be proved. We have given for the interested reader many problems in the last chapter of this book. With the study of these group neutrosophic rings, neutrosophic group rings, semigroup neutrosophic rings and neutrosophic semigroup rings, one natural question that would rise in the minds of the readers in what will be the algebraic structure when we take both the group and the ring to be neutrosophic. The answer is we can define and we will call them as only strong or generalized structures for they will also continue to be a neutrosophic rings. To this end we make the following definition.

**DEFINITION 3.3.5:** *Let $\langle G \cup I\rangle$ be a neutrosophic group of $\langle K \cup I\rangle$ be a neutrosophic field or a commutative neutrosophic ring $\langle K \cup I\rangle$ [$\langle G \cup I\rangle$] is defined to be a neutrosophic group neutrosophic ring which consists of all finite formal sums of the form $\sum_{i=1}^{n}\alpha_i g_i; n < \infty \ \ \alpha_i \in \langle K \cup I\rangle$ and $g_i \in \langle G \cup I\rangle.$*

*This algebraic structure has a ring structure, for this neutrosophic group neutrosophic ring is defined by replacing a group by a neutrosophic group and the ring by a neutrosophic ring.*

It is clear that neutrosophic group neutrosophic ring is a neutrosophic ring with unit. Neutrosophic group neutrosophic ring will be commutative if and only if the neutrosophic group is commutative.

We illustrate this by the following examples:

***Example 3.3.9:*** Let $\langle Z \cup I \rangle$ be the neutrosophic ring of integers and $\langle G \cup I \rangle = \{1,\ g,\ g^2,\ I,\ gI,\ g^2I/\ g^3 = 1$ and $I^2 = I\}$ be the neutrosophic group. Then $\langle Z \cup I \rangle\ [\langle G \cup I \rangle]$ is the neutrosophic group neutrosophic ring. This is a most generalized algebraic structure. This contains non trivial zero divisors and units. For $(1 - g)\ (1 + g + g^2) = 0$ is a zero divisor. $(I - Ig)\ (1 + g + g^2) = 0$ is the neutrosophic zero divisor.

This is clearly an infinite commutative neutrosophic ring of characteristic zero but is not a neutrosophic integral domain.

Now we give yet another example of a commutative finite neutrosophic ring built using a neutrosophic group and a neutrosophic ring.

***Example 3.3.10:*** Consider the neutrosophic ring $\langle Z_2 \cup I \rangle$ and the neutrosophic group $\langle G \cup I \rangle = \{1, g, g^2, g^3, g^4, g^5, g^6, g^7, I, gI, g^2I, g^3I, g^4I, g^5I, g^6I, g^7I/\ g^8 = 1$ and $I^2 = I\}$.

Clearly the neutrosophic group neutrosophic ring $\langle Z_2 \cup I \rangle$ ($\langle G \cup I \rangle$)) is a commutative finite neutrosophic ring having non trivial divisors of zero, both neutrosophic and ordinary; for take

$x = (1 + g^2 + g^4 + g^6)$,
$y = (1 + g^2)$ in $\langle Z_2 \cup I \rangle\ [\langle G \cup I \rangle]$.

Clearly x . y = 0. So x, y is a zero divisor.
Take
$\alpha = 1 + g^4$ and
$\beta = I + g^2 + Ig^4 + g^6$ in $\langle Z_2 \cup I \rangle\ [\langle G \cup I \rangle]$.

$$\alpha\beta = (1 + g^4)(1 + g^2 + g^4I + g^6)$$
$$= I + g^2 + g^4I + g^6 + Ig^4 + g^6 + I + g^2 = 0.$$

Thus $\alpha\beta$ is a neutrosophic zero divisor of $\langle Z_2 \cup I \rangle [\langle G \cup I \rangle]$.

***Example 3.3.11:*** Let $\langle Z_2 \cup I \rangle$ be the neutrosophic ring and $\langle G \cup I \rangle = \{1, g, g^2, g^3, g^4, g^5, I, gI, g^2I, g^3I, g^4I, g^5I/ g^6 = 1$ and $I^2 = I\}$ be the neutrosophic group. $\langle Z_2 \cup I \rangle [\langle G \cup I \rangle]$ is the neutrosophic group neutrosophic ring.

Consider $\alpha = (1 + g^2 + g^4)$ in $\langle Z_2 \cup I \rangle [\langle G \cup I \rangle]$. $\alpha^2 = \alpha$ thus $\alpha$ is an idempotent of the neutrosophic group neutrosophic ring. If we take $\alpha_1 = I + g^2I + g^4I$ in $\langle Z_2 \cup I \rangle [\langle G \cup I \rangle]$. Clearly $\alpha^2_1 = \alpha_1$ and $\alpha_1$ is a neutrosophic idempotent of $\langle Z_2 \cup I \rangle [\langle G \cup I \rangle]$. Take $\beta = 1 + g^3$.

Clearly $\beta^2 = 0$ so is a zero divisor. If we take $\alpha = 1 + g^3$ and $\gamma = I + g + g^2 + g^3I + g^4 + g^5$ in $\langle Z_2 \cup I \rangle [\langle G \cup I \rangle]$ then $\alpha\gamma = 0$, $\alpha$, $\gamma$ is a neutrosophic zero divisor.

Now we proceed on to give an example of a non commutative neutrosophic group neutrosophic ring.

***Example 3.3.12:*** Take $\langle Z_3 \cup I \rangle$ to be a neutrosophic ring. Let $\langle G \cup I \rangle = \langle D_{26} \cup I \rangle = \{1, a, b, b^2, b^3, b^4, b^5, ab, ab^2, ab^3, ab^4, ab^5, I, aI, bI, b^2I, b^3I, b^4I, b^5I, abI, ab^2I, ab^3I, ab^4I, ab^5I / a^2 = 1, b^6 = 1$ and $I^2 = I$ with $bab = a\}$ be a non-commutative neutrosophic group. Then $\langle Z_3 \cup I \rangle [\langle D_{26} \cup I \rangle]$ is a non commutative finite neutrosophic group neutrosophic ring.

It is left as an exercise for the reader to determine both neutrosophic and ordinary zero divisor, idempotents and units in $\langle Z_3 \cup I \rangle [\langle D_{26} \cup I \rangle]$.

Now we have got the following interesting results which can be made use of while studying the properties of these algebraic structure.

**THEOREM 3.3.4:** *Every neutrosophic group neutrosophic ring contains a proper subset which is a neutrosophic group ring.*

*Proof:* Given $\langle R \cup I \rangle$ $[\langle G \cup I \rangle]$ is a neutrosophic group neutrosophic ring. Consider the proper subset T = R $[\langle G \cup I \rangle]$ in $\langle R \cup I \rangle$ $[\langle G \cup I \rangle]$; Clearly T is a neutrosophic group ring.

For if we consider the example 3.3.12, $Z_3[\langle D_{26} \cup I \rangle]$ is the neutrosophic group ring contained in the neutrosophic group neutrosophic ring $\langle Z_3 \cup I \rangle$ $[\langle D_{26} \cup I \rangle]$.

It is seen that neutrosophic group neutrosophic ring is the most generalized form. For we see it contains group ring and group neutrosophic ring as proper subsets.

The proofs of these can be got in an analogous way. The interested reader can develop an algebraic structure analogous to the group ring. The problems given in the last chapter will be helpful to researchers for constructing these structures.

We define the ideals in these neutrosophic group neutrosophic ring as strong neutrosophic ideals. However we can have all types of subrings but only strong neutrosophic ideals. This is a marked difference between the ideals and subrings in neutrosophic group neutrosophic rings. One can define the notion of strong neutrosophic maximal ideals, strong neutrosophic minimal ideals, strong neutrosophic prime ideals and so an.

Now we proceed on to define some more special structures now only in case of these group neutrosophic rings, semigroup neutrosophic rings, neutrosophic group rings, neutrosophic semigroup rings and S-semigroup neutrosophic rings.

**DEFINITION 3.3.6:** *Let $\langle K \cup I \rangle$ [G] be a group neutrosophic ring. If KG the group ring is semisimple we call the group neutrosophic ring to be pseudo semisimple.*

Now we define when a semigroup neutrosophic ring is pseudo semisimple.

**DEFINITION 3.3.7:** *Let $\langle K \cup I \rangle[S]$ be a semigroup neutrosophic ring of a semigroup S over the neutrosophic ring $\langle K \cup I \rangle$. If the semigroup ring KS is semisimple then we call $\langle K \cup I \rangle$ [S] to be pseudo semisimple.*

Now we define when are neutrosophic group rings and neutrosophic semigroup rings pseudo semisimple.

**DEFINITION 3.3.8:** *Let K[⟨G ∪ I⟩] be a neutrosophic group ring. We say K [⟨G ∪I⟩] is pseudo semisimple if KG the group ring contained in K [⟨G ∪I⟩] is semisimple.*

Now we say when the neutrosophic semigroup ring K [⟨S ∪ I⟩] is pseudo semisimple.

**DEFINITION 3.3.9:** *Let K[⟨S ∪I⟩] be a neutrosophic semigroup ring of the neutrosophic semigroup ⟨S ∪I⟩ over the field K or a commutative ring with unit. We say K[⟨S ∪ I⟩] to be pseudo semisimple if KS the semigroup ring of the semigroup S over the ring K is semisimple.*

Interested reader can characterize pseudo semisimple rings.

Now we proceed on to define Smarandache semisimple in case of S-semigroup neutrosophic rings.

**DEFINITION 3.3.10:** *Let ⟨K ∪ I⟩ [S] be a S-semigroup neutrosophic ring. Let A be a proper subset of S which is a group. If the group ring KA is semisimple then we say ⟨K ∪I⟩ [S] is a Smarandache semisimple (S-semisimple).*

Now we proceed on to define the notion of pseudo domain or pseudo division ring.

**DEFINITION 3.3.11:** *Let ⟨K ∪I⟩ [G] be any group neutrosophic ring of the group G over the neutrosophic ring ⟨K ∪I⟩. We say ⟨K ∪I⟩ [G] to be a pseudo domain or a pseudo division ring if the group ring KG ⊆ ⟨K ∪I⟩ [G] is embeddable in a domain or a division ring respectively.*

***Note:*** In case of semigroup neutrosophic ring ⟨K ∪ I⟩ [S] where S is a semigroup and ⟨K ∪ I⟩ any neutrosophic ring, if the

semigroup ring KS ⊆ ⟨K ∪ I⟩ [S] is embeddable in a domain or a division ring then we say the semigroup neutrosophic ring is a pseudo domain or pseudo division ring respectively.

Further when we have a neutrosophic group ring K[⟨G ∪ I⟩] where K is a field or a commutative ring with unit and ⟨G ∪ I⟩ a neutrosophic group, we say K [⟨G ∪ I⟩] to be a pseudo domain or a pseudo division ring if the group ring KG ⊆ K[⟨G ∪ I⟩] is embeddable in a domain or a division ring respectively.

Suppose K[⟨S ∪ I⟩] be a neutrosophic semigroup ring of the neutrosophic semigroup ⟨S ∪ I⟩ over the ring K, we say K[⟨S ∪ I⟩] is a pseudo domain or a pseudo division ring if the semigroup ring KS ⊆ K[⟨S ∪ I⟩] is embeddable in a domain or a division ring respectively. Several interesting properties in this direction can be derived.

Now on similar lines we can define the notion of pseudo prime and pseudo semiprime in case of the 4 types of neutrosophic rings viz. group neutrosophic ring, semigroup neutrosophic ring, neutrosophic group ring and neutrosophic semigroup ring. Just we define it in case of group neutrosophic ring and neutrosophic group ring, the definition for the other two types follows in an analogous way.

**DEFINITION 3.3.12:** *Let ⟨K ∪ I⟩ [G] be a group neutrosophic ring of the group G over the neutrosophic ring ⟨K ∪I⟩. We say ⟨K ∪I⟩ [G] is pseudo prime if KG ⊆ ⟨K ∪ I⟩ [G] is prime. We say the group neutrosophic ring ⟨K ∪ I⟩ [G] is pseudo semiprime if KG ⊆ ⟨K ∪I⟩ [G] is semiprime.*

*If K [⟨G ∪ I⟩] be the neutrosophic group ring, of the neutrosophic group ⟨G ∪I⟩ over the field K (or a commutative ring with unit K). We say the neutrosophic group ring K (⟨G ∪ I⟩) is pseudo prime if the group ring KG ⊆K[⟨G ∪I⟩] is prime. Similarly K[⟨G ∪I⟩] is pseudo semiprime if KG ⊆K[⟨G ∪I⟩] is semiprime.*

On similar lines the reader is requested to define the notion of pseudo prime and pseudo semiprime in case of neutrosophic semigroup rings and semigroup neutrosophic rings.

Now when we speak of neutrosophic group neutrosophic ring we say [⟨R ∪ I⟩] [⟨G ∪ I⟩] is weakly pseudo semiprime if one of the subsets which is a group neutrosophic ring ⟨R ∪ I⟩ [G] or the neutrosophic group ring R[⟨G ∪ I⟩] is pseudo semiprime similarly we say the neutrosophic group neutrosophic ring [⟨K ∪ I⟩] [⟨G ∪ I⟩] is weakly pseudo prime if the neutrosophic group ring K [⟨G ∪ I⟩] or group neutrosophic ring ⟨K ∪ I⟩ [G] is pseudo prime.

On similar lines we define weak pseudo domain or weak pseudo division ring in the neutrosophic group neutrosophic ring if the group neutrosophic ring [⟨K ∪ I⟩][G] or neutrosophic group ring K [⟨G ∪ I⟩] is a pseudo domain or a pseudo division ring.

Like wise we define the notion of weak pseudo semisimple.

Now we define the notion of neutrosophic group neutrosophic ring to satisfy a pseudo polynomial identity.

**DEFINITION 3.3.13:** *Let ⟨K ∪ I⟩ [⟨G ∪ I⟩] be a neutrosophic group neutrosophic ring we say ⟨K ∪ I⟩ [⟨G ∪ I⟩] satisfies a pseudo polynomial identity of degree n if K[⟨G ∪ I⟩] or [⟨K ∪ I⟩] [G] satisfies a polynomial identity of degree n.*

Now we can define neutrosophic semigroup neutrosophic ring.

**DEFINITION 3.3.14:** *Let ⟨K ∪ I⟩ be a neutrosophic ring and ⟨S ∪ I⟩ neutrosophic semigroup. ⟨K ∪ I⟩ [⟨S ∪ I⟩] is defined as the neutrosophic semigroup neutrosophic ring analogous to neutrosophic group neutrosophic ring where the neutrosophic group is replaced by the neutrosophic semigroup. Thus ⟨K ∪ I⟩ [⟨S ∪ I⟩], the neutrosophic semigroup neutrosophic ring consists of all finite formal sums of the form $\sum_{i=1}^{n} k_i s_i$ where $n < \infty$ and $k_i \in ⟨K ∪ I⟩$ and $s_i \in ⟨S ∪ I⟩$.*

It is a matter of routine to verify ⟨K ∪ I⟩ [⟨S ∪ I⟩] is a neutrosophic ring. It is still interesting to note that the neutrosophic semigroup neutrosophic ring contains the following nice substructure viz. semigroup neutrosophic ring ⟨K ∪ I⟩ [S] and neutrosophic semigroup ring K⟨S ∪ I⟩.

All concepts defined for neutrosophic group neutrosophic rings can be defined for this class of neutrosophic ring keeping in mind that we have a semigroup instead of a group. Now we still proceed on to define a still generalized class of neutrosophic rings.

We call a neutrosophic semigroup ⟨S ∪ I⟩ to be a Smarandache neutrosophic semigroup if ⟨S ∪ I⟩ contains a proper subset P, such that P is a group under the operations of ⟨S ∪ I⟩; we just illustrate by an example a S- neutrosophic semigroup.

***Example 3.3.13:*** Let ⟨Z ∪ I⟩ be a neutrosophic semigroup under multiplication ×.P = {1, −1, x} is a semigroup. Thus (⟨Z ∪ I⟩, ×) is a S- neutrosophic semigroup.

For more about S- neutrosophic semigroup refer [142-3].

Now we define the notion of S- neutrosophic semigroup ring and S- neutrosophic semigroup neutrosophic rings.

**DEFINITION 3.3.15:** *Let K be any ring ⟨S ∪ I⟩ be a S-neutrosophic semigroup. K[⟨S ∪ I⟩] is the S- neutrosophic semigroup ring of the S- neutrosophic semigroup [⟨S ∪I⟩] over the ring K. K [⟨S ∪I⟩] consists of all finite formal sums of the form $\sum_{i=1}^{n} \alpha_i s_i$ where $\alpha_i \in K$ and $s_i \in ⟨S ∪ I⟩$ and $n < \infty$. As in case of S semigroup neutrosophic rings it can be easily verified that K[⟨S ∪I⟩] satisfies a ring structure. In fact K [⟨S ∪ I⟩] is a neutrosophic ring of a special type.*

Further K[⟨S ∪ I⟩] is a generalized structure for it contains all subsets which are neutrosophic semigroup rings. S-semigroup rings.

First we illustrate this structure by the following example:

***Example 3.3.14:*** Let $\langle S \cup I \rangle$ be the S- neutrosophic semigroup given by the following table:

| * | 0 | $g_1$ | $g_2$ | $g_3$ | $g_1I = I$ | $g_2I$ | $g_3I$ |
|---|---|-------|-------|-------|-----------|--------|--------|
| 0 | 0 | 0 | 0 | 0 | 0 | 0 | 0 |
| $g_1$ | 0 | $g_1$ | $g_2$ | $g_3$ | I | $g_2I$ | $g_3I$ |
| $g_2$ | 0 | $g_2$ | 0 | $g_2$ | $g_2I$ | 0 | $g_2I$ |
| $g_3$ | 0 | $g_3$ | $g_2$ | $g_1$ | $g_3I$ | $g_2I$ | I |
| $I = g_1I$ | 0 | I | $g_2I$ | $g_3I$ | I | $g_2I$ | $g_3I$ |
| $g_2I$ | 0 | $g_2I$ | 0 | $g_2I$ | $g_2I$ | 0 | $g_2I$ |
| $g_3I$ | 0 | $g_3I$ | $g_2I$ | I | $g_3I$ | $g_2I$ | I |

P = {$g_1$, $g_3$} is a group under * contained in $\langle S \cup I \rangle$ so $\langle S \cup I \rangle$ is a S- neutrosophic semigroup.

Take Q to be the field of rationals Q[$\langle S \cup I \rangle$] is S-neutrosophic semigroup ring.

Now we define still a generalized structure viz. the S-neutrosophic semigroup neutrosophic ring, $\langle K \cup I \rangle$ [$\langle S \cup I \rangle$] where $\langle K \cup I \rangle$ is a neutrosophic ring and $\langle S \cup I \rangle$ is a S-neutrosophic semigroup. $\langle K \cup I \rangle$ [$\langle S \cup I \rangle$], the S- neutrosophic semigroup neutrosophic ring is a neutrosophic ring and it is defined analogous to neutrosophic group neutrosophic ring.

It consists as usual finite formal sums of the form $\sum_{i=1}^{n} \alpha_i s_i; n < \infty$, $\alpha_i \in \langle K \cup I \rangle$ and $s_i \in \langle S \cup I \rangle$.

This S neutrosophic semigroup neutrosophic ring is commutative if and only if $\langle S \cup I \rangle$ is commutative. $\langle K \cup I \rangle$ [$\langle S \cup I \rangle$] is Smarandache commutative if and only if every proper subset of $\langle S \cup I \rangle$ which is a group under the operations of $\langle S \cup I \rangle$ is commutative. Smarandache weakly commutative if and only if $\langle S \cup I \rangle$ has atleast a proper subset A, which is a commutative group under the operation of $\langle S \cup I \rangle$.

It is still interesting to note that S- neutrosophic semigroup neutrosophic ring $\langle K \cup I \rangle$ [$\langle S \cup I \rangle$] is the most generalized structure among the neutrosophic rings defined for it contains a proper semigroup ring KS. It contains a proper S-semigroup ring KS, it contains a S-semigroup neutrosophic ring $\langle K \cup I \rangle$ [S], it contains a neutrosophic semigroup ring K[$\langle S \cup I \rangle$] and a S- neutrosophic semigroup ring K [$\langle S \cup I \rangle$].

Thus this is the presently known rich structure among neutrosophic rings for it contains also group rings, group neutrosophic rings and neutrosophic group ring and neutrosophic group neutrosophic ring given by the sets KA, (A $\subseteq \langle S \cup I \rangle$ is a group), $\langle K \cup I \rangle$ [A], K$\langle A \cup I \rangle$ and $\langle K \cup I \rangle$ [$\langle A \cup I \rangle$] respectively.

Now we have so far discussed only about the pseudo neutrosophic rings and pseudo neutrosophic ideals in neutrosophic rings.

First we have found it very difficult get pseudo neutrosophic rings which contain the unit. Only when we can find a non trivial class of pseudo neutrosophic rings with unit which are commutative we would be in a position to define the notion of group pseudo neutrosophic rings, semigroup pseudo neutrosophic rings and so on. We have given several examples of pseudo neutrosophic rings but they do not contain unit.

# SUGGESTED PROBLEMS

In this chapter we have suggested 246 number of problems for the researcher. We do this mainly to attract students towards this new notion of neutrosophic rings. Further it is pertinent to mention that some problems are simple, some little difficult and some problems can be taken as open problems for further research.

1. Find all the ideals in the neutrosophic group ring, $Z_2 \langle G \cup I \rangle$ where $G = \langle 1, g, g^2, \ldots, g^{11}, I, gI, g^2I, \ldots, g^{11}I / g^{12} = 1$ and $I^2 = I \rangle$.

2. Does the neutrosophic group ring $Z_2 \langle G \cup I \rangle$ defined in problem 1 have pseudo neutrosophic ideals?

3. Does the neutrosophic group ring $Z_2 \langle G \cup I \rangle$ defined in problem 1 have neutrosophic subrings which are not neutrosophic ideals?

4. Let $Z \langle G \cup I \rangle$ be the neutrosophic group ring of the group $\langle G \cup I \rangle = \{1, g, g^2, \ldots, g^6, I, gI, g^2I, \ldots, g^6I / g^7 = 1$ and $I^2 = I\}$. Does $Z\langle G \cup I \rangle$ have a minimal neutrosophic ideal?

5.  Can Z⟨G ∪ I⟩ given in problem 4 have a maximal neutrosophic ideal? How many maximal neutrosophic ideals does Z ⟨G ∪ I⟩ have?

6.  Is all neutrosophic ideals (of the group rings Z⟨G ∪ I⟩ given in problem 4) in Z⟨G ∪ I⟩ principal?

7.  Does neutrosophic group ring Z⟨G ∪ I⟩ given in problem 4 have prime neutrosophic ideals?

8.  Does the neutrosophic group ring Z⟨G ∪ I⟩ given in problem 4 satisfy ACC condition or DCC condition on ideals?

9.  Does the neutrosophic group ring Z(⟨G ∪ I⟩) given in problem 4 have any pseudo neutrosophic ideal?

10. Let $⟨G ∪ I⟩ = \left\{ \begin{pmatrix} a & b \\ c & d \end{pmatrix} \right\}$, a, b, c, d, ∈ ⟨$Z_3$ ∪ I⟩; the neutrosophic ring of integers modulo 3 and ad − bc ≠ 0}. Clearly ⟨G ∪ I⟩ is a non commutative finite group. Z (⟨G ∪ I⟩) be the neutrosophic group ring of ⟨G ∪ I⟩ over Z.

    a.  Does Z(⟨G ∪ I⟩) have neutrosophic right ideals? If so find them.

    b.  Does Z(⟨G ∪ I⟩) have neutrosophic left ideals?

    c.  Find all the two sided neutrosophic ideals of Z(⟨G ∪ I⟩).

    d.  Does Z⟨G ∪ I⟩ have any pseudo neutrosophic ideals?

    e.  Does Z⟨G ∪ I⟩ have any maximal neutrosophic ideal, if any find them?

f. Can $Z\langle G \cup I \rangle$ have any minimal neutrosophic ideal? Justify your claim!

g. Can $Z\langle G \cup I \rangle$ have any principal neutrosophic ideals?

h. Does $Z\langle G \cup I \rangle$ have any prime neutrosophic ideals?

i. Can in a neutrosophic group ring $Z\langle G \cup I \rangle$ we have a neutrosophic ideal which is both maximal and minimal. Justify your answer.

11. Give an example of a neutrosophic group ring which has only neutrosophic ideals but no pseudo neutrosophic ideals.

12. Does the neutrosophic group ring $Z_3 (\langle G \cup I \rangle)$ where $Z_3 = \{0, 1, 2\}$ and $G = \langle 1, g, I, gI / g^2 = 1$ and $I^2 = I\}$ have neutrosophic ideals? or pseudo neutrosophic ideals? Answer your claim!

13. Does there exist a neutrosophic group ring $R\langle G \cup I \rangle$ which has only pseudo neutrosophic ideals and no neutrosophic ideals?

14. Let $Z\langle G \cup I \rangle$ be a neutrosophic group ring of the neutrosophic group $\langle G \cup I \rangle = \{1, g, g^2, \ldots, g^{p-1}, I, gI, g^2I, \ldots, g^{p-1}I / g^p = 1$ where $p$ is an odd prime and $I^2 = I\}$ and $Z$ the ring of integers. Let $P = \langle 0, n (1 + g + \ldots + g^{p-1}), n (I + gI + \ldots + g^{p-1} I) / n \in Z \rangle \subseteq Z (\langle G \cup I \rangle)$. Is $P$ a pseudo neutrosophic ideal of $Z\langle G \cup I \rangle$? Is $P$ a neutrosophic ideal of $Z\langle G \cup I \rangle$? Does $Z\langle G \cup I \rangle$ have infinite number of neutrosophic ideals? Justify your claim.

15. Characterize those neutrosophic group rings, which are loyal?

16. Give an example of a loyal neutrosophic group ring.

17. Give an example of a loyal neutrosophic ideal of a neutrosophic group ring.

18. Let $R\langle G \cup I \rangle$ be given by $R = Z_p = \{0, 1, 2, \ldots, p - 1\}$; p a prime and $G = \{1, g, g^2, \ldots, g^{18}, I, gI, \ldots, g^{18}I / g^{19} = 1$ and $I^2 = I\}$. Can $R\langle G \cup I \rangle$ have any quasi neutrosophic ideals?

19. Given $R\langle G \cup I \rangle$ is a neutrosophic group ring where $R = Z_{12}$ and $G \cup I = \langle \{g / g^{16} = 1\} \cup I \rangle = \{1, g, \ldots, g^{15}, I, gI, \ldots, g^{15}I / g^{16} = 1$ and $I^2 = I\}$. Does $R\langle G \cup I \rangle$ have any neutrosophic ideals? (If so find them explicitly).

20. Does $R\langle G \cup I \rangle$ given in problem 19 have any quasi neutrosophic ideals? (If so find them explicitly).

21. Give an example of bonded quasi neutrosophic ideals in a neutrosophic group ring.

22. Give an example of a loyal quasi neutrosophic ideal in a neutrosophic ring.

23. Can the neutrosophic group ring $Z\langle G \cup I \rangle$ where $\langle G \cup I \rangle = \{1, g, g^2, g^3, g^4, g^5, I, gI, g^2I, g^3I, g^4I, g^5I / g^6 = 1$ and $I^2 = I\}$ have any
    a. Bonded quasi neutrosophic ideals?
    b. Loyal neutrosophic ideals?
    Justify your claim!

24. Does their exists a bonded quasi neutrosophic ideal which is loyal?

25. Obtain any interesting relation between the two types of quasi neutrosophic ideals, bonded and quasi!

26. Illustrate by an example a bonded quasi neutrosophic right ideal which is not a bonded quasi neutrosophic left ideal and vice versa.

27. Give an example of a loyal quasi neutrosophic left ideal of a neutrosophic group ring which is not a loyal quasi neutrosophic right ideal and (vice versa).

28. Suppose S is a quasi neutrosophic ideal say over two subrings $P_1$ and $P_2$. Can we say any relation exists between $P_1$ and $P_2$?

29. Does the neutrosophic group ring $Z_{15}$ $(\langle G \cup I \rangle)$ where $Z_{15}$ is the ring of integers modulo 15 and $\langle G \cup I \rangle = \langle \{g, g^{30} = 1\} \cup I \rangle$ have
    a. Quasi neutrosophic ideals?
    b. Bonded quasi neutrosophic ideals?
    c. Loyal quasi neutrosophic ideals?
    d. Pseudo quasi neutrosophic ideals?
    e. Pseudo quasi strong neutrosophic ideals?

30. Does the neutrosophic group ring $R\langle G \cup I \rangle$ where $R = Z$, the ring of integers and $\langle G \cup I \rangle = \langle \{g \mid g^{16} = 1\} \cup I \rangle = \{1, g, g^2, ..., g^{15}I \mid g^{16} = 1 \text{ and } I^2 = I\}$ contain
    a. a bonded pseudo quasi neutrosophic ideal?
    b. a strongly bonded pseudo quasi neutrosophic ideal?
    Justify your claim!

31. Give an example of a neutrosophic group ring, $R\langle G \cup I \rangle$ which has only pseudo quasi neutrosophic ideals and no strong pseudo quasi neutrosophic ideals.

32. Does their exist a neutrosophic group ring which has only strong pseudo quasi neutrosophic ideals?

33. Give an example of a neutrosophic group ring having only strong pseudo bonded neutrosophic ideals?

34. Illustrate by an example a neutrosophic group ring which has only pseudo bonded neutrosophic ideals and not strong pseudo bonded neutrosophic ideals.

35. Does their exists a neutrosophic group ring which has no neutrosophic idempotents and idempotents?

36. Find all the zero divisors and neutrosophic zero divisors in the neutrosophic group ring $Z_{12} \langle G \cup I \rangle$ where $Z_{12}$ is the ring of integers modulo 12 and $\langle G \cup I \rangle = \{1, g, g^2, g^3, g^4, g^5, I, gI, g^2I, g^3I, g^4I, g^5I / g^6 = 1$ and $I^2 = I\}$.

37. Find all the zero divisors and neutrosophic zero divisors in $Z_7 \langle G \cup I \rangle$ where $Z_7$ is the ring of integers modulo 7 and $\langle G \cup I \rangle = \{1, g, g^2, g^3, g^4, g^5, I, gI, g^2I, g^3I, g^4I, g^5I / g^6 = 1$ and $I^2 = I\}$.

38. Let $Z_p = \{0, 1, 2, \ldots, p-1\}$ be the ring of integers modulo p, p a prime. $\langle G \cup I \rangle$ be a neutrosophic group of order n. Let $Z_p (\langle G \cup I \rangle)$ be the neutrosophic group ring of the neutrosophic group $\langle G \cup I \rangle$ over $Z_p$.
    a. Find the number of zero divisors and neutrosophic zero divisors of $Z_p (\langle G \cup I \rangle)$ if p/n?
    b. Does $Z_p (\langle G \cup I \rangle)$ have zero divisors and neutrosophic zero divisors if (p, n) = 1.
    c. Find all units and neutrosophic units in $Z_p \langle G \cup I \rangle$ when p/n and (p, n) = 1.

39. Given $Z_{10} = \{0, 1, 2, 3, 4, 5, \ldots, 9\}$, the ring of integers modulo 10. $\langle G \cup I \rangle = \{1, g, g^2, g^3, g^4, I, gI, g^2I, g^3I g^4I / g^5 = 1$ and $I^2 = I\}$ the neutrosophic group. Find all units and neutrosophic units of $Z_{10} (\langle G \cup I \rangle)$.

40. Does the neutrosophic group ring $Z_2 (\langle G \cup I \rangle)$ where $\langle G \cup I \rangle = \{1, g, g^2, g^3, g^4, I, gI, g^2I, g^3I / g^4 = 1$ and $I^2 = I\}$ have nontrivial units and non trivial neutrosophic units?

41. Find all idempotents and neutrosophic idempotents in the neutrosophic group ring $Z_2 (\langle G \cup I \rangle)$ given in the problem 40.

42. Can a torsion free neutrosophic group over any field have neutrosophic zero divisors and zero-divisors?

*Note:* This problem is analogous to the problem of zero divisor conjecture "Does a group ring have zero divisors?" if G is a torsion free group over any field F.

43. Give an example of a neutrosophic unique product group, which is not a strongly neutrosophic unique product group.

44. Give an example of a strongly neutrosophic two unique product.

45. Can a neutrosophic unique product group ring FG over a field F have neutrosophic divisors of zero?

46. Does a strongly neutrosophic two unique product group ring FG over a field F have neutrosophic divisors of zero?

47. Give an example of an infinite neutrosophic group which is not a neutrosophic unique product group.

48. Define the notion of full neutrosophic ring of matrices over a neutrosophic ring and illustrate it with examples.

49. Illustrate by an example the notion of pseudo neutrosophic submodule.

50. Give an example of a neutrosophic submodule which is not a pseudo neutrosophic submodule.

51. Does their exist any pseudo neutrosophic submodule which is a neutrosophic submodulo and vice versa?

52. Define the notion of irreducible neutrosophic module and illustrate it by an example.

53. Is every neutrosophic ideal of a neutrosophic group ring $R(\langle G \cup I \rangle)$, a neutrosophic R-module? Justify your claim.

54. Define free neutrosophic module and illustrate it with an example!

55. Obtain conditions under which $R$ $(\langle G \cup I \rangle) \cong S$ $(\langle H \cup I \rangle)$ where R and S commutative rings with unit or fields and $\langle G \cup I \rangle$ and $\langle H \cup I \rangle$ are neutrosophic groups?

Remarks: This is analogous to the famous isomorphism problem for group rings, i.e., "Is $RG \cong SH$? R and S fields and G and H groups."

This conjecture has remained open for over 3 to 4 decades. Only partial solutions are available.

56. Let $Z_{12}$ $(\langle G \cup I \rangle)$ and $Z_{10}$ $(\langle H \cup I \rangle)$ be two neutrosophic group rings where $\langle G \cup I \rangle = \{1, g, g^2, g^3, g^4, I, gI, g^2I, g^3I, g^4I \ / \ g^5 = 1$ and $I^2 = I\}$ and $\langle H \cup I \rangle = \{1, h, h^2, h^3, h^4, h^5, I, hI, h^2I, h^3I, h^4I, h^5I \ / \ h^6 = 1$ and $I^2 = I\}$.
Find a neutrosophic homomorphism from $Z_{12}$ $(\langle G \cup I \rangle)$ to $Z_{10}$ $(\langle H \cup I \rangle)$ so that ker $\phi$ is a subring of $Z_{12}$ $(\langle G \cup I \rangle)$ having some neutrosophic elements. Justify your claim.
(Hint: It may happen $\phi$ $(\alpha g \ I) = \phi$ $(\alpha g)$. (I) where $\phi(\alpha g) = 0$ so that $\phi$ $(\alpha g I) = 0$).

Now we propose a open problem for the reader which needs more researching and number theoretic techniques.

57. Let $Z_p$ $(\langle G \cup I \rangle)$ and $Z_q$ $(\langle H \cup I \rangle)$ be two neutrosophic group rings, where p and q are primes, $p \neq q$ and $o\langle G \cup I \rangle = 2n$ and $o\langle H \cup I \rangle = 2m$. Find a necessary and sufficient condition on p, q, m and n so that $Z_p$ $(\langle G \cup I \rangle) \cong Z_q$ $(\langle H \cup I \rangle)$.
    a. Is it possible or not with p and q distinct primes?
    b. If p and q are non primes will it be possible to find conditions on m and n so that $Z_p$ $(\langle G \cup I \rangle) \cong Z_q(\langle H \cup I \rangle)$? (Assume both $\langle G \cup I \rangle$ and $\langle H \cup I \rangle$ are abelian).

58. Let $Z_2$ $(\langle G \cup I \rangle)$ and $Z_3$ $(\langle H \cup I \rangle)$ be two neutrosophic group rings where $\langle G \cup I \rangle = \{1, g, g^2, I, gI, g^2I \ / \ g^3 = 1$ and $I^2 = I\}$ and $\langle H \cup I \rangle = \{1, h, I, hI \ / \ h^2 = 1$ and $I^2 = I\}$. Construct a neutrosophic group ring homomorphism with a non trivial kernel.

59. Let R be the field of reals. Is R ($\langle G \cup I \rangle$) $\cong$ R ($\langle H \cup I \rangle$) imply $\langle G \cup I \rangle \cong \langle H \cup I \rangle$ as neutrosophic groups?

60. Give an example of isomorphic neutrosophic group rings where $\langle G \cup I \rangle \neq \langle H \cup I \rangle$.

61. Give an example of a semiprime neutrosophic group ring.

62. Is the neutrosophic group ring $Z_7$ ($\langle G \cup I \rangle$) where $Z_7 = \{0, 1, 2, \ldots, 6\}$ ring of integers modulo 7 and $\langle G \cup I \rangle = \{1, g, g^2, g^3, g^4, g^5, g^6, I, gI, g^2I, g^3I, g^4I, g^5I, g^6I / g^7 = 1$ and $I^2 = I\}$ semiprime? Justify your claim.

63. Is the group ring K ($\langle G \cup I \rangle$) where K is a field of characteristic zero and $\langle G \cup I \rangle$ a neutrosophic group of finite order semiprime? Prove your claim.

64. Make some innovative definitions about neutrosophic group rings satisfying polynomial identity.

65. Suppose $\langle G \cup I \rangle$ is a neutrosophic group having a abelian neutrosophic subgroup $\langle A \cup I \rangle$ with $[\langle G \cup I \rangle : \langle A \cup I \rangle] = n < \infty$, will the neutrosophic group ring K $[\langle G \cup I \rangle]$ satisfy the standard polynomial identity of degree 2n? Justify your claim.

66. Suppose K $[\langle G \cup I \rangle]$ satisfies a nontrivial polynomial identity of degree n. Is $[\langle G \cup I \rangle : \Delta^I] \leq \underline{n}$?

67. Let K $[\langle G \cup I \rangle]$ be a prime neutrosophic group ring. Suppose K $[\langle G \cup I \rangle]$ satisfies a polynomial identity of degree n.
    a. Will $\Delta [\langle G \cup I \rangle]$ be torsion free and abelian?
    b. Is $[\langle G \cup I \rangle : \Delta^I] \leq n/2$? Substantiate your claim.

68. Suppose $\langle G \cup I \rangle$ is a finitely generated group and K $[\langle G \cup I \rangle]$ satisfies a polynomial identity of degree n. Is $[\langle G \cup I \rangle : \Delta^I] \leq n/2$? Justify your answer.

69. Let $K[\langle G \cup I \rangle]$ be a semiprime neutrosophic group ring with a or the neutrosophic center $C^I$. $\alpha \in C^I$ is a zero divisor in $K[\langle G \cup I \rangle]$ if and only if it is a zero divisor in $C^I$. Illustrate this by a concrete example.

70. Let $K[\langle G \cup I \rangle]$ be a semiprime group ring and let $\alpha \in K[\langle G \cup I \rangle]$. If $\alpha$ is not a right divisor of zero, then there exists $\gamma \in K[\langle G \cup I \rangle]$ such that $\theta(\alpha\gamma)$ is central and not a zero divisor in $K[\langle G \cup I \rangle]$. $\theta : K[\langle G \cup I \rangle] \rightarrow K[\Delta^I]$ ($\theta$ is clearly a K-linear map but not a ring homomorphism defined by if
$$\alpha = \sum_{x \in \langle G \cup I \rangle} k_x x \, , \; \theta(\alpha) = \sum_{x \in \Delta^I = \Delta(\langle G \cup I \rangle)} k_x x \, ).$$
Prove the result. Illustrate this result by a example.

71. Is the neutrosophic group ring $K[\langle G \cup I \rangle]$ where $K = Q$ (ring of rationals) and $\langle G \cup I \rangle = \{1, g, g^2, g^3, \ldots, g^{14}, I, gI, \ldots, g^{14}I \; / \; I^2 = I$ and $g^{15} = 1\}$ semiprime? Find all neutrosophic divisors of zero and divisors of zero. Clearly $C^I = K[\langle G \cup I \rangle]$.

72. Define for a neutrosophic group ring $K[\langle G \cup I \rangle]$, the notion of semisimple. Hence or otherwise define the Jacobson radical of $K[\langle G \cup I \rangle]$.

73. Find a necessary and sufficient condition for a neutrosophic group ring $K[\langle G \cup I \rangle]$ to be semisimple. If KG is semisimple does it imply $K[\langle G \cup I \rangle]$ is semisimple? If $K[\langle G \cup I \rangle]$ is semisimple does it imply $K[G]$ is semisimple?

74. Suppose the neutrosophic group $\langle G \cup I \rangle$ has no elements of order p and K is a field of characteristic p. Is the neutrosophic group ring $K[\langle G \cup I \rangle]$ semisimple?

75. Let $Z_{19}[\langle G \cup I \rangle]$ be a neutrosophic group ring where $Z_{19} = \{0, 1, 2, \ldots 18\}$ the prime field of characteristic 19 and $\langle G \cup I \rangle = \{1, g, g^2, g^3, g^4, \ldots, g^{11}, I, gI, g^2I, g^3I, g^4I, \ldots, g^{11}I \; / \; g^{12} = 1$ and $I^2 = I\}$ . Is $Z_{19}[\langle G \cup I \rangle]$ semisimple?

76. Define the notion of solvability in case of neutrosophic groups. Illustrate a solvable neutrosophic group by an example.

77. Describe JK [⟨G ∪ I⟩]. Is it always a neutrosophic nil ideal?

78. Find a necessary and sufficient condition for a neutrosophic group ring K [⟨G ∪ I⟩] to satisfy a polynomial identity.

79. Derive some interesting properties about group neutrosophic rings.

80. If G is a unique product group and ⟨R ∪ I⟩ the neutrosophic ring of reals. Can the group neutrosophic ring ⟨R ∪ I⟩ [G] have non trivial divisors of zero. Establish your claim!

81. Let G be a two unique product group and ⟨Z ∪ I⟩ the neutrosophic ring of integers, can the group neutrosophic ring ⟨Z ∪ I⟩ [G] have non trivial neutrosophic zero divisors? Can ⟨Z ∪ I⟩ [G] have non trivial units and non trivial neutrosophic units?

82. Find the necessary and sufficient condition for the group neutrosophic ring ⟨K ∪ I⟩ [G] to be semisimple.

83. Give an example of a group neutrosophic ring which is not semisimple.

84. Give an example of a group neutrosophic ring which is semisimple.

85. Let $S_\infty$ denote the infinite symmetric group consisting of all permutations moving only finitely many points. Is ⟨K ∪ I⟩ [$S_\infty$] semisimple for all neutrosophic fields?

86. Suppose G is a group having no elements of order p, in case ⟨K ∪ I⟩, the neutrosophic ring has characteristic p. Is ⟨K ∪ I⟩ [G] semisimple?

87. Let $G = \langle g/g^{15} = 1 \rangle$ and $\langle Z_5 \cup I \rangle$ be the neutrosophic field of characteristic 5. Is the group neutrosophic ring $\langle Z_5 \cup I \rangle$ [G] semisimple?

88. Let $G = \langle g / g^{16} = 1 \rangle$ and $\langle Z_2 \cup I \rangle$ be the neutrosophic field of characteristic 2. Is the group neutrosophic ring $\langle Z_2 \cup I \rangle$[G] semisimple?

89. Find necessary and sufficient condition for $\langle K \cup I \rangle$ [G], the group neutrosophic ring to satisfy a polynomial identity.

90. Let $\langle K \cup I \rangle$ be a neutrosophic field of characteristic p and suppose $\langle K \cup I \rangle$ [G] satisfies a polynomial identity of degree n. Does G necessarily have a normal p-abelian subgroup A of finite index? If this is the case, are [G:A] and |A| necessarily bounded as functions of n?

91. Let $\langle K \cup I \rangle$ be a neutrosophic field of characteristic p and let $\langle K \cup I \rangle$[G] satisfy a polynomial identity of degree n. Is S(G) solvable of bounded derived length? (S(G) the subgroup of G generated by all solvable normal subgroups of G).

92. Let $\langle K \cup I \rangle$ [G] satisfy a polynomial identity of degree n;
    a. Is $[G : \Delta(G)] \leq n/2$ in all cases?
    b. What are the best possible bounds for [G: S(G)] and for [G: A] in case $[\langle K \cup I \rangle]$ [G] is semiprime?

93. Let $\langle K \cup I \rangle$ be a neutrosophic field of characteristic 0. Is it necessarily true that $D^n (\langle K \cup I \rangle$ [G]) = $\gamma^n(G)$?

94. Let $\langle K \cup I \rangle$ be a neutrosophic field of characteristic p. Will the group neutrosophic ring $\langle K \cup I \rangle$[G] satisfy a polynomial identity if G has a p abelian subgroup of finite index?
    a. Is the converse of the above statement true?
    b. Illustrate this situation by an example when p = 3.

95. Suppose $\langle K \cup I \rangle$ [G] satisfies a polynomial identity of degree n will $[G : \Delta(G) ] \leq n/2$ and $| \Delta(G) | < \infty$?

96. Let $\alpha$, $\beta \in \langle K \cup I \rangle$ [G] with $\alpha\beta = 1$. Must $\beta\alpha$ necessarily also equal 1 in a group neutrosophic ring?

97. Let $\langle K \cup I \rangle$ be a neutrosophic field of characteristic zero. Let $e \in \langle K \cup I \rangle$ [G] be an idempotent. Is the trace of e necessarily rational? If so find an algebraic proof of this fact.

98. Let G be a torsion free group and $\langle K \cup I \rangle$ be any neutrosophic field. Is it true that $\langle K \cup I \rangle$ [G] has no proper divisors of zero. Are all units in $\langle K \cup I \rangle$ [G] necessarily trivial?

99. What will be the situation if in the above problem (i.e., problem 98)G is replaced by a nilpotent group?

100. Find the necessary and sufficient condition for $\langle K \cup I \rangle$ [G] to be embeddable in a neutrosophic division ring.

101. Find a necessary and sufficient condition for the group neutrosophic ring $\langle K \cup I \rangle$ [G], which is embeddable in a neutrosophic division ring to satisfy a polynomial identity.

102. Let $\langle K \cup I \rangle$ be a neutrosophic field of characteristic p and let $e \in \langle K \cup I \rangle$ [G] be an idempotent or a neutrosophic idempotent.
     a.  Is the trace of e necessarily rational?
     b.  Does trace of $e \in GF(p)$?

103. Find a necessary and sufficient condition for the group neutrosophic ring $\langle K \cup I \rangle$ [G] to be
     a.  Primitive.
     b.  Semiperfect.
     c.  Algebraic.

104. Is it true $\langle K \cup I \rangle$ [G] is noetherian if and only if G has a series of subgroup $G = G_0 \supseteq G_1 \ldots \supseteq G_n = \langle 1 \rangle$, such that $G_{i+1}$ is normal in $G_i$ and $G_i / G_{i+1}$ is either a finite group or an infinite cyclic group.

105. Find all units, neutrosophic units, idempotents, neutrosophic idempotents of the group neutrosophic ring $\langle Z_{12} \cup I \rangle$ $[A_4]$? ($A_4$, the alternating subgroup of $S_4$).

106. Find all zero divisors, neutrosophic zero divisors and units and neutrosophic units of the group neutrosophic ring $\langle Z_{12} \cup I \rangle$ $[D_{2.6}]$ where $D_{2.6} = \{a, b \mid a^2 = b^6 = 1, bab = a\}$ the dihedral group of order 12.

107. When will the two group neutrosophic rings, $\langle R \cup I \rangle$ $[G] \cong \langle S \cup I \rangle$ $[H]$?

108. Will $\langle Z_p \cup I \rangle$ $[G]$ where $G = \langle g \, / \, g^n = 1 \rangle$ have zero divisors? $(n, p) = 1$ or $p/n$.

109. Find all zero divisors and neutrosophic zero divisors of $\langle Z_p \cup I \rangle$ $(S_p)$, where p is a prime and $S_p$ is the symmetric group of degree p.

110. When will the group neutrosophic ring $\langle Z_p \cup I \rangle$ $[G]$ isomorphic with $\langle Z_q \cup I \rangle$ $[G']$? Find conditions on p and q.

111. Let A and B be neutrosophic ideals in $\langle K \cup I \rangle$ $[G]$;
    a. Is $\theta(A)$ an ideal in $\langle K \cup I \rangle$ $[A]$?
    b. Is $A \neq 0$ if and only if $\theta(A) \neq 0$?
    c. $AB = 0$, does it imply that $\theta(A) \theta(B) = 0$?

112. Prove if $\langle K \cup I \rangle$ is a neutrosophic field of characteristic 0 & G any group, $\langle K \cup I \rangle [G]$ contains no non-zero-nil-ideals.

113. Can we say if $\langle K \cup I \rangle$ is a neutrosophic field of characteristic 0 where the field K is not algebraic over the rational field Q and G is a group, will the neutrosophic group ring $\langle K \cup I \rangle$ $[G]$ be semisimple?

114. Let $\langle K \cup I \rangle$ be a neutrosophic field. G a group and let F be a family of subgroups of G. Suppose that if $H \in F$ then will $\langle K \cup I \rangle$ $[H]$ be semisimple?

115. If S is a finite subset of G then will there exist a subgroup H in F with S ⊆ H. Is ⟨K ∪ I⟩ [G] semisimple?

116. Let N be a normal subgroup of G with G/N = B, B abelian suppose ⟨K ∪ I⟩ [N] is semisimple and G has no elements of order p in case ⟨K ∪ I⟩ has characteristic p. Will ⟨K ∪ I⟩ [G] be semisimple?

117. Let G be a finite non abelian group and let ⟨K ∪ I⟩ [G] satisfy a polynomial identity of degree n. What can one say about elements x ∈ G \ (Z(G)) is in the centre of G?

118. Let S(4) be the symmetric semigroup. Is the semigroup neutrosophic ring ⟨$Z_4$ ∪ I⟩ ($S_4$) semisimple? Justify your claim.

119. Find all neutrosophic zero divisors of the semigroup neutrosophic ring ⟨$Z_p$ ∪ I⟩ S(n); where (n, p) = 1.

120. Can the semigroup neutrosophic ring ⟨$Z_p$ ∪ I⟩ S(p) have zero divisors and neutrosophic zero divisors?

121. Find a nontrivial neutrosophic right ideal of the semigroup neutrosophic ring ⟨$Z_2$ ∪ I⟩ (S(3)) which is not a neutrosophic left ideal. Also find a neutrosophic right ideal of ⟨$Z_2$ ∪ I⟩ [S (3)].

122. Find condition on the neutrosophic ring ⟨K ∪ I⟩ and on the semigroup S so that the semigroup neutrosophic ring ⟨K ∪ S⟩S is semisimple?

123. Find under what conditions the semigroup neutrosophic ring ⟨K ∪ I⟩ [S] satisfies a polynomial identity.

124. Find a necessary and sufficient condition for the semigroup neutrosophic ring ⟨K ∪ I⟩[S] to be embeddable in a neutrosophic division ring.

125. Find a necessary and sufficient condition for the semigroup neutrosophic ring $\langle K \cup I \rangle$ [G] to be
   a. Primitive.
   b. Semiperfect.
   c. Algebraic.

126. Let $\alpha, \beta \in \langle K \cup I \rangle$ [S], a semigroup neutrosophic ring with $\alpha\beta = 1$. Must $\beta\alpha$ necessarily also be equal to 1? Justify your claim.

127. Can we find left units in the semigroup neutrosophic ring $\langle Z_3 \cup I \rangle$ [S₃]?

128. Does their exists non trivial units and neutrosophic units in the group neutrosophic ring $\langle Z_2 \cup I \rangle$ (S₆)?

129. Can a semigroup neutrosophic ring have a principal neutrosophic ideal? Justify your claim.

130. Find a necessary and sufficient condition for the semigroup neutrosophic ring $\langle K \cup I \rangle$ [S] to be semisimple.

131. Is $\langle Z \cup I \rangle$ [S₄] semisimple?

132. Can $\langle Q \cup I \rangle$ [S(n) ] be semisimple?

133. Will $\langle Z_p \cup I \rangle$ [S(p)] be semisimple?

134. When will a semigroup neutrosophic ring $\langle K \cup I \rangle$ [S] have nil ideals?

135. Find a necessary and sufficient condition for the semigroup neutrosophic ring $\langle K \cup I \rangle$[S] to be prime.

136. Find a necessary and sufficient condition for the semigroup neutrosophic ring $\langle K \cup I \rangle$ [S] to be semiprime.

137. Is the semigroup neutrosophic ring $\langle Q \cup I \rangle$ [S(4)] semiprime?

138. Find the center of the semigroup neutrosophic ring $\langle Z_2 \cup I \rangle$ [S(3)].

139. Let $\langle K \cup I \rangle$ [S] be a S-semigroup neutrosophic ring. Find conditions on S and $\langle K \cup I \rangle$ so that $\langle K \cup I \rangle$ [S] contains a non trivial group ring which is semiprime.

140. Does the S-semigroup neutrosophic ring $\langle Z_{11} \cup I \rangle$ [S(5)] contain a semiprime group ring?

141. Find interesting properties about the semigroup neutrosophic rings $\langle R \cup I \rangle$S.

142. Find a necessary and sufficient condition for a S-semigroup neutrosophic ring $\langle K \cup I \rangle$ [S] to have non trivial prime group rings.

143. Can a S-semigroup neutrosophic ring be prime? If so give an example of one such S-semigroup neutrosophic ring.

144. Does there exist a S-semigroup neutrosophic ring which satisfies a polynomial identity of degree n?

145. Find conditions on the S-semigroup so that the S-semigroup neutrosophic ring $\langle K \cup I \rangle$ [S] contains a group ring which satisfies the polynomial identity of degree n.

146. Can the S-semigroup neutrosophic ring $\langle Q \cup I \rangle$ [S(n)] satisfy the polynomial identity of degree n?

147. When will the S-semigroup neutrosophic ring $\langle K \cup I \rangle$ [S] be semiprime?

148. Find conditions on the S-semigroup S so that the S-semigroup neutrosophic ring $\langle K \cup I \rangle$[S] has a proper subset which is a group ring say K[G], which is semiprime.

149. Find a necessary and sufficient condition for the S-semigroup neutrosophic ring to be semisimple.

150. Is the S-semigroup neutrosophic ring $\langle R \cup I \rangle$ [S(5)] semisimple where R is the field of reals.

151. Find conditions on the S-semigroup neutrosophic ring $\langle K \cup I \rangle$ [S] to have a group ring which is semisimple.

152. Find the necessary and sufficient conditions for the S-semigroup ring $\langle K \cup I \rangle$ [S] to be
    a. Primitive.
    b. Semiperfect.
    c. Algebraic.

153. Find a necessary and sufficient condition for the S-semigroup neutrosophic ring $\langle K \cup I \rangle$ [S] to be embeddable in a neutrosophic division ring.

154. Can we say in a S-semigroup neutrosophic ring $\langle K \cup I \rangle$ [S] for $\alpha, \beta \in \langle K \cup I \rangle$ [S] with $\alpha\beta = 1$, then must it imply that $\beta\alpha = 1$?

155. Find all units in the S-semigroup neutrosophic ring $\langle Z_{10} \cup I \rangle$ [S(4)].

156. Find a maximal neutrosophic ideal of the S-semigroup ring $\langle Z_2 \cup I \rangle$ [S(5)].

157. Find a right neutrosophic ideal of $\langle Z_3 \cup I \rangle$ [S(3)] which is not a left neutrosophic ideal of $\langle Z_3 \cup I \rangle$ [S(3)]. Does $\langle Z_3 \cup I \rangle$ [S (3)] contain any

    a. Minimal neutrosophic ideal?
    b. Principal neutrosophic ideal?
    c. Prime neutrosophic ideal?

158. Find all units and neutrosophic units of the neutrosophic group neutrosophic ring $[\langle Q \cup I \rangle]$ $[\langle G \cup I \rangle]$ where $\langle G \cup I \rangle = \{1, g, g^2, g^3, g^4, g^5, I, gI, g^2I, g^3I, g^4I, g^5I \, / \, I^2 = I$ and $g^6 = 1\}$. Does $\langle Q \cup I \rangle$ $[\langle G \cup I \rangle]$ have non trivial zero divisors?

159. Does the neutrosophic group neutrosophic ring $\langle Z_p \cup I \rangle$ $[\langle G \cup I \rangle]$ where $\langle G \cup I \rangle = \{1, g, g^2, \ldots, g^{p-1}, I, gI, g^2I, \ldots, g^{p-1}I$ / $I^2 = I$ and $g^p = 1\}$ p a prime have non trivial neutrosophic zero divisors, units and idempotents?

160. Find strong neutrosophic ideal of the neutrosophic group neutrosophic ring $\langle Z_6 \cup I \rangle$ $[\langle G \cup I \rangle]$ where $\langle G \cup I \rangle = \{1, g, g^2, g^3, \ldots g^{11}, I, gI, \ldots g^{11}I$ / $g^{12} = 1$ and $I^2 = I\}$. Does this neutrosophic group neutrosophic ring have both maximal and minimal neutrosophic strong ideals?

161. Find a neutrosophic strong right ideal of the neutrosophic group neutrosophic ring $\langle Z_{10} \cup I \rangle$ $[D_{2.4} \cup I]$ where $\langle D_{2.4} \cup I \rangle$ = $\{1, a, b, b^2, b^3, ab, ab^2, ab^3, I, aI, bI, b^2I, b^3I, abI, ab^2I, ab^3I$ / $I^2 = I$, $a^2 = b^4 = 1$, $bab = a\}$. Does it have any neutrosophic strong prime ideal? Find a neutrosophic left ideal in $\langle Z_{10} \cup I \rangle$ $[\langle D_{2.4} \cup I \rangle]$.

162. Find some interesting properties about neutrosophic group neutrosophic rings.

163. Find a necessary and sufficient condition for a neutrosophic group neutrosophic ring $\langle K \cup I \rangle$ $[\langle G \cup I \rangle]$ to satisfy a polynomial identity.

164. Find a necessary and sufficient condition for the neutrosophic group neutrosophic ring to have non trivial non zero nil ideals.

165. Find a necessary and sufficient condition for the neutrosophic group neutrosophic ring to be embeddable in a neutrosophic division ring.

166. When is the neutrosophic group neutrosophic ring $\langle K \cup I \rangle$ $[\langle G \cup I \rangle]$ semisimple?

167. Is the neutrosophic group neutrosophic ring $\langle Q \cup I \rangle$ $[\langle G \cup I \rangle]$ where $\langle G \cup I \rangle = = \{1, g, g^2, \ldots, g^{p-1}, I, gI, g^2I, \ldots, g^{p-1}I$ / $g^p = 1$ and $I^2 = I\}$ semiprime? Justify your claim.

168. Find the necessary and sufficient condition for the neutrosophic group neutrosophic ring $\langle K \cup I \rangle [\langle G \cup I \rangle]$ to be
     a. Primitive.
     b. Semiperfect.
     c. Algebraic.

169. Suppose $\langle K \cup I \rangle [G \cup I]$ satisfies a polynomial identity of degree n. Is $[\langle G \cup I \rangle : \Delta (\langle G \cup I \rangle)] \leq n/2$ in all cases?

170. Suppose $\langle K \cup I \rangle [\langle G \cup I \rangle]$ is the neutrosophic group neutrosophic ring where $\langle K \cup I \rangle$ is a neutrosophic field of characteristic p, p > 0.

     a. Is $\langle K \cup I \rangle [\langle G \cup I \rangle]$ semiprime?

     b. Does $\langle K \cup I \rangle [\langle G \cup I \rangle]$ semiprime imply $\Delta (\langle G \cup I \rangle)$ has no elements of order p?

     c. If $\Delta (\langle G \cup I \rangle)$ has no elements of order p does it imply $\langle K \cup I \rangle [\langle G \cup I \rangle]$ is semiprime?

     d. Suppose $\langle G \cup I \rangle$ has no finite neutrosophic normal subgroup with order divisible by p does it not imply $\langle K \cup I \rangle [\langle G \cup I \rangle]$ is semiprime?

171. Let $\langle K \cup I \rangle$ be a neutrosophic field of characteristic p, p > 0 and let $\langle G \cup I \rangle$ have no elements of order p. Does $\langle K \cup I \rangle [\langle G \cup I \rangle]$ have non zero nil ideal?

172. $\langle K \cup I \rangle [\langle G \cup I \rangle]$ is a neutrosophic group neutrosophic ring where $\langle K \cup I \rangle$ is a neutrosophic field of characteristic 0. Is $\langle K \cup I \rangle [\langle G \cup I \rangle]$ semiprime? Justify your claim.

173. $\langle Q \cup I \rangle [\langle G \cup I \rangle]$ be a neutrosophic group neutrosophic ring. Is $\langle Q \cup I \rangle [\langle G \cup I \rangle]$ semiprime?

174. Let $\langle K \cup I \rangle$ $[\langle G \cup I \rangle]$ be a prime neutrosophic group neutrosophic ring. Suppose $\langle K \cup I \rangle$ $[\langle G \cup I \rangle]$ satisfies a polynomial identity of degree n. Is $\Delta$ $(\langle G \cup I \rangle)$ torsion free abelian? Will $[\langle G \cup I \rangle : \Delta (\langle G \cup I \rangle)] \leq n / 2$.

175. Prove or disprove if $[\langle K \cup I \rangle]$ $[\langle G \cup I \rangle]$ satisfy a polynomial identity of degree n, then $[\langle G \cup I \rangle : \Delta (\langle G \cup I \rangle)] \leq n$ !

176. Prove if $\langle G \cup I \rangle$ is an infinite torsion free neutrosophic group, then $\langle K \cup I \rangle$ $[\langle G \cup I \rangle]$ has no divisor of zero. If K is of characteristic p, p an odd prime will $\langle K \cup I \rangle$ $[\langle G \cup I \rangle]$ have non trivial units?

177. When will any two neutrosophic group neutrosophic rings be isomorphic?

    Note: The famous isomorphism problem for group rings has been proposed some 40 years ago.

178. Find a necessary and sufficient condition for the group neutrosophic ring $\langle R \cup I \rangle$ [G] to be pseudo semisimple.

179. Let $\langle R \cup I \rangle$ [G] be a group neutrosophic ring where $\langle R \cup I \rangle$ is the neutrosophic field of reals G a t-u-p group. Is $\langle R \cup I \rangle$ [G] pseudo semisimple?

180. Let $\langle Z_n \cup I \rangle$ [G] be a group neutrosophic ring, n not a prime. G a u-p group. Is $\langle Z_n \cup I \rangle$ [G] pseudo semisimple?

181. Give an example of a group neutrosophic ring which is not pseudo semisimple.

182. Find a necessary and sufficient condition so that the semigroup neutrosophic ring $\langle K \cup I \rangle$ [G] is pseudo semisimple.

183. Let $\langle Q \cup I \rangle$ [S(n)] be the semigroup neutrosophic ring where $\langle Q \cup I \rangle$ is the neutrosophic field of rationals. Is $\langle Q \cup I \rangle$ [S(n)] pseudo semisimple?

184. Give an example of a semigroup neutrosophic ring which is not pseudo semisimple.

185. Give a class of semigroup neutrosophic rings which are pseudo semisimple.

186. Find a necessary and sufficient condition for a neutrosophic group ring to be pseudo semisimple.

187. Find a necessary and sufficient condition for a neutrosophic semigroup ring to be pseudo semisimple.

188. Find a necessary and sufficient condition for a S-semigroup neutrosophic ring to be Smarandache semisimple.

189. Give a class of S -semigroup neutrosophic rings which are Smarandache semisimple.

190. Is the S-semigroup neutrosophic ring $\langle Q \cup I \rangle$ [S (n)] Smarandache semisimple? Is the S-semigroup neutrosophic ring $\langle Z_n \cup I \rangle$ [S(n)] Smarandache semisimple?

191. Give a class of neutrosophic group ring which is pseudo semisimple?

192. Is Q [$\langle D_{2n} \cup I \rangle$] pseudo semisimple?

193. Is the neutrosophic group ring $Z_n$ ($\langle D_{2m} \cup I \rangle$) pseudo semisimple? [n / m, (m, n) = 1 or m / n].

194. Let p > 2 be a prime, $Z_p$ the prime field of characteristic p. $Z_p$ [$\langle G \cup I \rangle$] be the neutrosophic group, $\langle G \cup I \rangle$ has no elements of order p. Is $Z_p$ [$\langle G \cup I \rangle$] pseudo semisimple? Justify your claim!

195. Let Q [$\langle S \cup I \rangle$] be the neutrosophic semigroup ring where $\langle S \cup I \rangle \cong \langle Z_n \cup I \rangle$; $\langle Z_n \cup I \rangle$ is a semigroup generated under multiplication modulo n and Q the field of rationals. Is Q[$\langle S \cup I \rangle$] pseudo semisimple?

196. Obtain a necessary and sufficient condition for the neutrosophic semigroup ring K [⟨S ∪ I⟩] to be a pseudo domain or pseudo division ring.

197. Find a class of neutrosophic semigroup rings K [⟨S ∪ I⟩] which is a pseudo domain or a pseudo division ring.

198. Is the neutrosophic semigroup ring Q [⟨S ∪ I⟩] (where S ≅ $Z_n$, the semigroup and ⟨S ∪ I⟩ is a neutrosophic semigroup under multiplication modulo n; (n a prime)) a pseudo domain or a pseudo division ring?

199. Find a necessary and sufficient condition for the neutrosophic group ring K [⟨G ∪ I⟩] where ⟨G ∪ I⟩ is a neutrosophic group to be a pseudo domain or a pseudo division ring.

200. Give a class of neutrosophic group rings K [⟨G ∪ I⟩] to be a pseudo domain or a pseudo division ring.

201. Is the neutrosophic group ring K [⟨G ∪ I⟩] where G is a t-u-p group or t-p group and K any field a pseudo domain or a pseudo division ring? Justify your answer.

202. Is [⟨K ∪ I⟩][G] ≅ K [⟨G ∪ I⟩]? (G a group ⟨G ∪ I⟩ a neutrosophic group) i.e., Is a group neutrosophic ring isomorphic to a neutrosophic group ring?

203. Is [K ∪ I][S] ≅ K[⟨G ∪ I⟩]? i.e., is a semigroup neutrosophic ring isomorphic to a neutrosophic semigroup ring?

204. Find a necessary and sufficient condition for a neutrosophic group ring to be
    a. Pseudo prime.
    b. Pseudo semiprime.

205. Obtain a class of pseudo prime (pseudo semiprime) neutrosophic group rings.

206. Obtain a necessary and sufficient condition for a group neutrosophic ring to be
    a. Pseudo prime.
    b. Pseudo semiprime.

207. Find a class of group neutrosophic rings which are pseudo prime.

208. Is the group neutrosophic ring $\langle Q \cup I \rangle$ $[S_n]$ pseudo prime?

209. Is the semigroup neutrosophic ring $\langle Q \cup I \rangle$ $[S(n)]$ semiprime?

210. Is the group neutrosophic ring $\langle R \cup I \rangle$ $[G]$ where $[R \cup I]$ is the neutrosophic real field and G a torsion free abelian group pseudo semiprime?

211. Can the group neutrosophic ring $\langle Z_{10} \cup I \rangle$ $[G]$ where G is a cyclic group of order 20 be pseudo semiprime or pseudo prime?

212. Can the semigroup neutrosophic ring $[Z_{10} \cup I]$ $[S(10)]$ be pseudo semiprime? Justify your claim.

213. Can the neutrosophic group ring $Z[\langle D_{2.5} \cup I \rangle]$ be pseudo semiprime? Justify your answer.

214. Q $[\langle D_{2.5} \cup I \rangle]$ be the neutrosophic group ring. Is Q $[\langle D_{2.5} \cup I \rangle]$ pseudo prime? Justify your claim.

215. Characterize those neutrosophic group neutrosophic ring $[\langle K \cup I \rangle]$ $[\langle G \cup I \rangle]$ which are pseudo semiprime.
    a. Weakly semiprime,
    b. Weakly pseudo prime,
    c. Weakly pseudo semisimple,
    d. Weakly pseudo domain or weakly pseudo division ring.

216. Give a class of neutrosophic group neutrosophic rings [⟨K ∪ I⟩] [⟨G ∪ I⟩] which are weakly pseudo semisimple.

217. Is the neutrosophic group ring [⟨K ∪ I⟩] [⟨G ∪ I⟩] where ⟨K ∪ I⟩ = ⟨Q ∪ I⟩ the neutrosophic field of rationals, ⟨G ∪ I⟩ the neutrosophic group, where G is the infinite cyclic group a weak pseudo domain? Justify your claim.

218. Find a necessary and sufficient condition for the neutrosophic group neutrosophic ring to satisfy a pseudo polynomial identity of degree n.

219. Find a class of neutrosophic group neutrosophic ring which satisfies the pseudo polynomial identity of degree n.

220. Let [⟨K ∪ I⟩] [⟨G ∪ I⟩] be the neutrosophic group neutrosophic ring where [⟨K ∪ I⟩] = [⟨Q ∪ I⟩] the neutrosophic ring of rationals and ⟨G ∪ I⟩ = {1, g, … $g^{n-1}$, I, gI, … $g^{n-1}I$ / $g^n$ = 1 and $I^2$ = I} be the neutrosophic group. Does the neutrosophic group neutrosophic ring satisfy any pseudo polynomial identity?

221. Find some interesting properties about neutrosophic semigroup neutrosophic rings.

222. Find a necessary and sufficient condition for the neutrosophic semigroup neutrosophic ring [⟨K ∪ I⟩] [⟨S ∪ I⟩] to be semisimple.

223. Find conditions on the neutrosophic ring ⟨K ∪ I⟩ and the neutrosophic semigroup, ⟨S ∪ I⟩ so that the neutrosophic semigroup neutrosophic ring [⟨K ∪ I⟩] [⟨S ∪ I⟩] satisfies a polynomial identity of degree n.

224. Find a necessary and sufficient condition for the neutrosophic semigroup neutrosophic ring to be embeddable in a neutrosophic division ring.

225. Find a necessary and sufficient conditions for the neutrosophic semigroup neutrosophic ring [⟨K ∪ I⟩] [⟨S ∪ I⟩] to be
    a. Primitive.
    b. Prime.
    c. Semiprime.
    d. Semiperfect.
    e. Algebraic.

226. Find a necessary and sufficient condition on [⟨K ∪ I⟩] [⟨S ∪ I⟩] to be
    a. Weakly pseudo prime.
    b. Weakly pseudo semisimple.

227. Obtain some interesting properties about S-neutrosophic semigroup rings K [⟨S ∪ I⟩].

228. Find a necessary and sufficient condition for the S-neutrosophic semigroup ring to be semisimple.

229. Is the S-neutrosophic semigroup ring Q [⟨S ∪ I⟩] (where S ≅ {$Z_{12}$, semigroup under multiplication modulo 12 and ⟨S ∪ I⟩ is the S-neutrosophic semigroup} semisimple?

230. Find the necessary and sufficient condition for the S-neutrosophic semigroup ring to be
    a. Prime.
    b. Semiprime.
    c. Semi perfect and
    d. Embeddable in a neutrosophic division ring.

231. Obtain some interesting properties about the S-neutrosophic semigroup neutrosophic rings [⟨K ∪ I⟩] [⟨S ∪ I⟩].

232. Find a necessary and sufficient condition for the S-neutrosophic semigroup neutrosophic ring [⟨K ∪ I⟩] [⟨S ∪ I⟩] to be semisimple.

233. When will the S-neutrosophic semigroup neutrosophic ring $[\langle K \cup I \rangle] [\langle S \cup I \rangle]$ be pseudo semisimple?

234. Find necessary and sufficient condition for the S-neutrosophic semigroup neutrosophic ring $[\langle K \cup I \rangle] [\langle S \cup I \rangle]$ to be
    a. Prime.
    b. Semiprime.
    c. Perfect.
    d. Algebraic.

235. Find a necessary and sufficient condition for the S-neutrosophic semigroup neutrosophic ring $[\langle K \cup I \rangle] [\langle S \cup I \rangle]$ to satisfy the polynomial identity of degree n.

236. Find conditions on the neutrosophic ring $\langle K \cup I \rangle$ and on the S-neutrosophic semigroup $\langle S \cup I \rangle$, so that the S-neutrosophic semigroup neutrosophic ring is
    a. pseudo prime
    b. pseudo semiprime
    c. pseudo embeddable in a neutrosophic domain or in a neutrosophic division ring.
    d. pseudo semisimple.

237. Let $\langle Z_{12} \cup I \rangle = \langle S \cup I \rangle$ be a S-neutrosophic semigroup and $\langle K \cup I \rangle$ the neutrosophic rational field of characteristic zero. Does the S-neutrosophic semigroup neutrosophic ring $[\langle K \cup I \rangle] [\langle S \cup I \rangle]$ have maximal neutrosophic ideals? Find zero divisors, units and idempotents in $[\langle K \cup I \rangle] [\langle S \cup I \rangle]$. Show in general every neutrosophic subring of $[\langle K \cup I \rangle] [\langle S \cup I \rangle]$ is not a neutrosophic ideal of $[\langle K \cup I \rangle] [\langle S \cup I \rangle]$.

    Is this S-neutrosophic semigroup neutrosophic ring semiperfect? Justify your claim.

238. Does their exists a S-neutrosophic semigroup neutrosophic ring $[\langle K \cup I \rangle] [\langle S \cup I \rangle]$ which is a neutrosophic domain or a neutrosophic division ring?

239. Find condition on $\langle K \cup I \rangle$ and $\langle S \cup I \rangle$ so that $[\langle K \cup I \rangle] [\langle S \cup I \rangle]$ contains a commutative domain.

240. Does their exist a class of pseudo neutrosophic rings which contains 1?

241. Characterize those pseudo neutrosophic ring with 1 if they exist.

242. Does the group neutrosophic ring $\langle K \cup I \rangle [S_n]$ have any pseudo neutrosophic subring?

243. Can the group neutrosophic ring $\langle Z_{12} \cup I \rangle [S_6]$ have any nontrivial pseudo neutrosophic subring? Justify your claim.

244. Does the semigroup neutrosophic ring $\langle Z_6 \cup I \rangle [S(3)]$ have any pseudo neutrosophic subring?

245. Find a necessary and sufficient condition for the group neutrosophic ring $[\langle K \cup I \rangle][G]$ to have pseudo neutrosophic subrings and pseudo neutrosophic ideals.

246. Find all pseudo neutrosophic subrings of $\langle Z_9 \cup I \rangle [D_{2.9}]$.

# REFERENCES

An elaborate references is given mainly for the reader to study the notion of neutrosophic ring in those innovative directions carried out by several researchers.

1.    Abian Alexander and McWorter William, *On the structure of pre p-rings*, Amer. Math. Monthly, Vol. 71, 155-157, (1969).

2.    Adaoula Bensaid and Robert. W. Vander Waal, *Non-solvable finite groups whose subgroups of equal order are conjugate*, Indagationes Math., New series, No.1(4), 397-408, (1990).

3.    Allan Hayes, *A characterization of f-ring without non-zero nilpotents*, J. of London Math. Soc., Vol. 39, 706-707, (1969).

4.    Allevi.E, *Rings satisfying a condition on subsemigroups*, Proc. Royal Irish Acad., Vol. 88, 49-55, (1988).

5.    Andruszkiewicz.R, *On filial rings*, Portug. Math., Vol. 45, 136-149, (1988).

6.    Atiyah.M.F and MacDonald.I.G, *Introduction to Commutative Algebra*, Addison Wesley, (1969).

7.    Aubert.K.E and Beck.L, *Chinese rings*, J. Pure and Appl. Algebra, Vol.24, 221-226, (1982).

8.    Bovdi Victor and Rosa A.L., *On the order of unitary subgroup of a modular group algebra*, Comm. in Algebra, Vol. 28, 897-1905, (2000).

9. Chen Huanyin, *On generalized stable rings*, Comm. in Algebra, Vol. 28, 1907-1917, (2000).

10. Chen Huanyin, *Exchange rings having stable range one,* Inst. J. Math. Sci., Vol. 25, 763-770, (2001).

11. Chen Huanyin, *Regular rings with finite stable range*, Comm. in Algebra, Vol. 29, 157-166, (2001).

12. Chen Jain Long and Zhao Ying Gan, *A note on F-rings*, J. Math. Res. and Expo. Vol. 9, 317-318, (1989).

13. Connel.I.G, *On the group ring*, Canada. J. Math. Vol. 15, 650-685, (1963).

14. Corso Alberto and Glaz Sarah, *Guassian ideals and the Dedekind- Merlens lemma*, Lecture notes in Pure and Appl. Math., No. 217, 113-143, Dekker, New York, (2001).

15. Cresp.J and Sullivan.R.P, *Semigroup in rings*, J. of Aust. Math. Soc., Vol. 20, 172-177, (1975).

16. Dubrovin.N.I, *Chain rings*, Rad. Math., Vol. 5, 97-106, (1989).

17. Dutta. T.K, *On generalized semi-ideal of a ring*, Bull. Calcutta Math. Soc., Vol. 74, 135-141, (1982).

18. Erdogdu.V, *Regular multiplication rings*, J. Pure and Appl. Algebra, Vol. 31 , 55-59 (1989).

19. Gray.M, *A Radical Approach to Algebra*, Addison Wesley, (1970).

20. Gupta.V, *A generalization of strongly regular rings*, Acta. Math. Hungar. , Vol. 43, 57-61, (1984).

21.    Han Juncheol and Nicholson W.K., *Extensions of clean rings*, Comm. in Algebra, Vol. 29, 2589-2595, (2001).

22.    Han Yang, *Strictly wild algebras with radical square zero*, Arch-Math. (Basel), Vol. 76, 95-99, (2001).

23.    Herstein.I.N, *Topics in Algebra*, John Wiley and Sons, (1964).

24.    Herstein.I.N, *Topics in Ring theory*, Univ. of Chicago Press, (1969).

25.    Higman.D.G, *Modules with a group of operators*, Duke Math. J., Vol. 21, 369-376, (1954).

26.    Higman.G, *The units of group rings*, Proc. of the London Math. Soc., Vol. 46, 231-248, (1940).

27.    Hirano Yasuyuki and Suenago Takashi, *Generalizations of von Neuman regular rings and n-like rings*, Comment Math. Univ. St. Paul., Vol. 37, 145-149, (1988).

28.    Hirano Yasuyuki, *On $\pi$-regular rings with no infinite trivial subring*, Math. Scand., Vol. 63, 212-214, (1988).

29.    Jacobson.N, *Theory of rings*, American Mathematical Society, (1943).

30.    Jacobson.N, *Structure of ring*, American Mathematical Society, (1956).

31.    Jin Xiang Xin, *Nil generalized Hamiltonian rings*, Heilongiang Daxue Ziran Kexue Xuebao, No. 4, 21-23, (1986).

32.    Johnson P.L, *The modular group ring of a finite p-group*, Proc. Amer. Math. Soc., Vol. 68, 19-22, (1978).

33.    Katsuda.R, *On Marot Rings*, Proc. Japan Acad., Vol. 60, 134-138, (1984).

34.    Kim Nam Kyin and Loe Yang, *On right quasi duo-rings which are π-regular,* Bull. Korean Math. Soc., Vol. 37, 217-227, (2000).

35.    Kishimoto Kazuo and Nagahara Takas, *On G-extension of a semi-connected ring*, Math. J. of Okayama Univ., Vol. 32, 25-42, (1990).

36.    Krempa.J, *On semigroup rings*, Bull. Acad. Poland Sci. Ser. Math. Astron. Phy., Vol. 25, 225-231, (1977).

37.    Lah Jiang, *On the structure of pre J-rings*, Hung-Chong Chow, 65[th] anniversary volume, Math. Res. Centre, Nat. Taiwan Univ., 47-52, (1962).

38.    Lang.S, *Algebra*, Addison Wesley, (1984).

39.    Ligh.S and Utumi.Y, *Direct sum of strongly regular rings and zero rings*, Proc. Japan Acad., Vol. 50, 589-592, (1974).

40.    Lin Jer Shyong, *The structure of a certain class of rings*, Bull. Inst. Math. Acad. Sinica, Vol. 19, 219-227, (1991).

41.    Louis Dale, *Monic and Monic free ideals in a polynomial semiring*, Proc. Amer. Math. Soc., Vol.56, 45-50, (1976).

42.    Louis Halle Rowen, *Ring theory*, Academic Press, (1991).

43.    Mason.G, *Reflexive ideals*, Comm. in Algebra, Vol. 17, 1709-1724, (1981).

44.    Northcott.D.G, *Ideal theory*, Cambridge Univ. Press, (1953).

45.    Padilla Raul, Smarandache Algebraic structures, Bull. of Pure and Appl. Sci., Vol. 17E, 119-121, (1998).

46.    Passman.D.S, *Infinite Group Rings*, Pure and Appl. Math., Marcel Dekker, (1971).

47.    Passman.D.S, *The Algebraic Structure of Group Rings*, Inter-science Wiley, (1977).

48.    Peric Veselin, *Commutativity of rings inherited by the location of Herstein's condition*, Rad. Math., Vol. 3, 65-76, (1987).

49.    Pimenov K.L. and Yakovlev. A. V., *Artinian Modules over a matrix ring, Infinite length modules*, Trends Math. Birkhauser Basel, Bie. 98, 101-105, (2000).

50.    Putcha Mohan.S and Yaqub Adil, *Rings satisfying a certain idempotency condition*, Portugal Math. No.3, 325-328, (1985).

51.    Raftery.J.G, *On some special classes of prime rings*, Quaestiones Math., Vol. 10, 257-263, (1987).

52.    Ramamurthi.V.S, *Weakly regular rings*, Canada Math. Bull., Vol. 166, 317-321, (1973).

53.    Richard.P.Stanley, *Zero square rings,* Pacific. J. of Math., Vol. 30, 8-11, (1969).

54.    Searcold.O.Michael, *A Structure theorem for generalized J rings*, Proc. Royal Irish Acad., Vol. 87, 117-120, (1987).

55.    Shu Hao Sun, *On the least multiplicative nucleus of a ring*, J. of Pure and Appl. Algebra, Vol. 78, 311-318, (1992).

56.    Smarandache, Florentin, (editor), Proceedings of the First International Conference on Neutrosophy, Neutrosophic

Set, Neutrosophic Probability and Statistics, University of New Mexico, (2001).

57.     Smarandache, Florentin, *A Unifying Field in Logics: Neutrosophic Logic,* Preface by Charles Le, American Research Press, Rehoboth, 1999, 2000. Second edition of the Proceedings of the First International Conference on *Neutrosophy, Neutrosophic Logic, Neutrosophic Set, Neutrosophic Probability and Statistics,* University of New Mexico, Gallup, (2001).

58.     Smarandache, Florentin, *Special Algebraic Structures*, in Collected Papers, Abaddaba, Oradea, Vol.3, 78-81 (2000).

59.     Smarandache Florentin, Multi structures and Multi spaces, (1969) www.gallup.unm.edu/~smarandache/transdis.txt

60.     Smarandache, Florentin, *Definitions Derived from Neutrosophics*, In Proceedings of the First International Conference on Neutrosophy, Neutrosophic Logic, Neutrosophic Set, Neutrosophic Probability and Statistics, University of New Mexico, Gallup, 1-3 December (2001).

61.     Smarandache, Florentin, *Neutrosophic Logic— Generalization of the Intuitionistic Fuzzy Logic,* Special Session on Intuitionistic Fuzzy Sets and Related Concepts, International EUSFLAT Conference, Zittau, Germany, 10-12 September 2003.

62.     Van Rooyen.G.W.S, *On subcommutative rings*, Proc. of the Japan Acad., Vol. 63, 268-271, (1987).

63.     Vasantha Kandasamy W.B., *On zero divisors in reduced group rings over ordered groups*, Proc. of the Japan Acad., Vol. 60, 333-334, (1984).

64.     Vasantha Kandasamy W.B., *On semi idempotents in group rings*, Proc. of the Japan Acad., Vol. 61, 107-108, (1985).

65.    Vasantha Kandasamy W.B., *A note on the modular group ring of finite p-group*, Kyungpook Math. J., Vol. 25, 163-166, (1986).

66.    Vasantha Kandasamy W.B., *Zero Square group rings*, Bull. Calcutta Math. Soc., Vol. 80, 105-106, (1988).

67.    Vasantha Kandasamy W.B., *On group rings which are p-rings*, Ganita, Vol. 40, 1-2, (1989).

68.    Vasantha Kandasamy W.B., *Semi idempotents in semi group rings*, J. of Guizhou Inst. of Tech., Vol. 18, 73-74, (1989).

69.    Vasantha Kandasamy W.B., *Semigroup rings which are zero square ring*, News Bull. Calcutta Math. Soc., Vol. 12, 8-10, (1989).

70.    Vasantha Kandasamy W.B., *A note on the modular group ring of the symmetric group $S_n$*, J. of Nat. and Phy. Sci., Vol. 4, 121-124, (1990).

71.    Vasantha Kandasamy W.B., *Idempotents in the group ring of a cyclic group*, Vikram Math. Journal, Vol. X, 59-73, (1990).

72.    Vasantha Kandasamy W.B., *Regularly periodic elements of a ring*, J. of Bihar Math. Soc., Vol. 13, 12-17, (1990).

73.    Vasantha Kandasamy W.B., *Semi group rings of ordered semigroups which are reduced rings*, J. of Math. Res. and Expo., Vol. 10, 494-493, (1990).

74.    Vasantha Kandasamy W.B., *Semigroup rings which are p-rings*, Bull. Calcutta Math. Soc., Vol. 82, 191-192, (1990).

75.    Vasantha Kandasamy W.B., *A note on pre J-group rings*, Qatar Univ. Sci. J., Vol. 11, 27-31, (1991).

76.    Vasantha Kandasamy W.B., *A note on semigroup rings which are Boolean rings*, Ultra Sci. of Phys. Sci., Vol. 3, 67-68, (1991).

77.    Vasantha Kandasamy W.B., *A note on the mod p-envelope of a cyclic group*, The Math. Student, Vol.59, 84-86, (1991).

78.    Vasantha Kandasamy W.B., *A note on units and semi idempotents elements in commutative group rings*, Ganita, Vol. 42, 33-34, (1991).

79.    Vasantha Kandasamy W.B., *Inner Zero Square ring*, News Bull. Calcutta Math. Soc., Vol. 14, 9-10, (1991).

80.    Vasantha Kandasamy W.B., *On E-rings*, J. of Guizhou. Inst. of Tech., Vol. 20, 42-44, (1991).

81.    Vasantha Kandasamy W.B., *On semigroup rings which are Marot rings*, Revista Integracion, Vol.9, 59-62, (1991).

82.    Vasantha Kandasamy W.B., *Semi idempotents in the group ring of a cyclic group over the field of rationals*, Kyungpook Math. J., Vol. 31, 243-251, (1991).

83.    Vasantha Kandasamy W.B., *A note on semi idempotents in group rings*, Ultra Sci. of Phy. Sci., Vol. 4, 77-78, (1992).

84.    Vasantha Kandasamy W.B., *Filial semigroups and semigroup rings*, Libertas Mathematica, Vol.12, 35-37, (1992).

85.    Vasantha Kandasamy W.B., *n-ideal rings*, J. of Southeast Univ., Vol. 8, 109-111, (1992).

86.    Vasantha Kandasamy W.B., *On generalized semi-ideals of a groupring*, J. of Qufu Normal Univ., Vol. 18, 25-27, (1992).

87.    Vasantha Kandasamy W.B., *On subsemi ideal rings*, Chinese Quat. J. of Math., Vol. 7, 107-108, (1992).

88.    Vasantha Kandasamy W.B., *On the ring $Z_2 S_3$*, The Math. Student, Vol. 61, 246-248, (1992).

89.    Vasantha Kandasamy W.B., *Semi group rings that are pre-Boolean rings*, J. of Fuzhou Univ., Vol. 20, 6-8, (1992).

90.    Vasantha Kandasamy W.B., *Group rings which are a direct sum of subrings*, Revista Investigacion Operacional, Vol. 14, 85-87, (1993).

91.    Vasantha Kandasamy W.B., *On strongly sub commutative group ring*, Revista Ciencias Matematicas, Vol. 14, 92-94, (1993).

92.    Vasantha Kandasamy W.B., *Semigroup rings which are Chinese ring*, J. of Math. Res. and Expo., Vol.13, 375-376, (1993).

93.    Vasantha Kandasamy W.B., *Strong right S-rings*, J. of Fuzhou Univ., Vol. 21, 6-8, (1993).

94.    Vasantha Kandasamy W.B., *s-weakly regular group rings*, Archivum Mathematicum, Tomus. 29, 39-41, (1993).

95.    Vasantha Kandasamy W.B., *A note on semigroup rings which are pre p-rings*, Kyungpook Math. J., Vol.34, 223-225, (1994).

96.    Vasantha Kandasamy W.B., *A note on the modular semigroup ring of a finite idempotent semigroup*, J. of Nat. and Phy. Sci., Vol. 8, 91-94, (1994).

97.    Vasantha Kandasamy W.B., *Coloring of group rings*, J. Inst. of Math. and Comp. Sci., Vol. 7, 35-37, (1994).

98.   Vasantha Kandasamy W.B., *f-semigroup rings*, The Math. Edu., Vol. XXVIII, 162-164, (1994).

99.   Vasantha Kandasamy W.B., *J-semigroups and J-semigroup rings*, The Math. Edu., Vol. XXVIII, 84-85, (1994).

100.  Vasantha Kandasamy W.B., *On a new type of group rings and its zero divisor*, Ult. Sci. of Phy. Sci., Vol. 6, 136-137, (1994).

101.  Vasantha Kandasamy W.B., *On a new type of product rings*, Ult. Sci. of Phy. Sci., Vol.6, 270-271, (1994).

102.  Vasantha Kandasamy W.B., *On a problem of the group ring $Z_p S_n$* , Ult. Sci. of Phy. Sci., Vol.6, 147, (1994).

103.  Vasantha Kandasamy W.B., *On pseudo commutative elements in a ring*, Ganita Sandesh, Vol. 8, 19-21, (1994).

104.  Vasantha Kandasamy W.B., *On rings satisfying $A^r = b^s = (ab)^t$*, Proc. Pakistan Acad. Sci., Vol. 31, 289-292, (1994).

105.  Vasantha Kandasamy W.B., *On strictly right chain group rings*, Hunan. Annele Math., Vol. 14, 47-49, (1994).

106.  Vasantha Kandasamy W.B., *On strong ideal and subring of a ring*, J. Inst. Math. and Comp. Sci., Vol.7, 197-199, (1994).

107.  Vasantha Kandasamy W.B., *On weakly Boolean group rings*, Libertas Mathematica, Vol. XIV, 111-113, (1994).

108.  Vasantha Kandasamy W.B., *Regularly periodic elements of group ring*, J. of Nat. and Phy. Sci., Vol. 8, 47-50, (1994).

109.  Vasantha Kandasamy W.B., *Weakly Regular group rings*, Acta Ciencia Indica.,  Vol. XX, 57-58, (1994).

110.    Vasantha Kandasamy W.B., *Group rings which satisfy super ore condition*, Vikram Math. J., Vol. XV, 67-69, (1995).

111.    Vasantha Kandasamy W.B., *Obedient ideals in a finite ring*, J. Inst. Math. and Comp. Sci., Vol. 8, 217-219, (1995).

112.    Vasantha Kandasamy W.B., *On group semi group rings*, Octogon, Vol. 3, 44-46, (1995).

113.    Vasantha Kandasamy W.B., *On Lin group rings*, Zesztyty Naukowe Poli. Rzes., Vol. 129, 23-26, (1995).

114.    Vasantha Kandasamy W.B., *On Quasi-commutative rings*, Caribb. J. Math. Comp. Sci. Vol.5, 22-24, (1995).

115.    Vasantha Kandasamy W.B., *On semigroup rings in which $(xy)^n = xy$*, J. of Bihar Math. Soc., Vol. 16, 47-50, (1995).

116.    Vasantha Kandasamy W.B., *On the mod p-envelope of $S_n$*, The Math. Edu., Vol. XXIX, 171-173, (1995).

117.    Vasantha Kandasamy W.B., *Orthogonal sets in group rings*, J. of Inst. Math. and Comp. Sci., Vol.8, 87-89, (1995).

118.    Vasantha Kandasamy W.B., *Right multiplication ideals in rings*, Opuscula Math., Vol.15, 115-117, (1995).

119.    Vasantha Kandasamy W.B., *A note on group rings which are F-rings*, Acta Ciencia Indica, Vol. XXII, 251-252, (1996).

120.    Vasantha Kandasamy W.B., *Finite rings which has isomorphic quotient rings formed by non-maximal ideals*, The Math. Edu., Vol. XXX, 110-112, (1996).

121.    Vasantha Kandasamy W.B., *$I^*$-rings*, Chinese Quat. J. of Math., Vol. 11, 11-12, (1996).

122.    Vasantha Kandasamy W.B., *On ideally strong group rings*, The Math. Edu., Vol. XXX, 71-72, (1996).

123.    Vasantha Kandasamy W.B., *Gaussian Polynomial rings*, Octogon, Vol.5, 58-59, (1997).

124.    Vasantha Kandasamy W.B., *On semi nilpotent elements of a ring*, Punjab Univ. J. of Math. , Vol. XXX, 143-147, (1997).

125.    Vasantha Kandasamy W.B., *On tripotent elements of a ring*, J. of Inst. of Math. and Comp. Sci., Vol. 10, 73-74, (1997).

126.    Vasantha Kandasamy W.B., *A note on f-group rings without non-zero nilpotents*, Acta Ciencia Indica, Vol. XXIV, 15-17, (1998).

127.    Vasantha Kandasamy W.B., *Inner associative rings*, J. of Math. Res. and Expo., Vol. 18, 217-218, (1998).

128.    Vasantha Kandasamy W.B., *On a quasi subset theoretic relation in a ring*, Acta Ciencia Indica, Vol. XXIV, 9-10, (1998).

129.    Vasantha Kandasamy W.B., *On SS-rings*, The Math. Edu., Vol. XXXII, 68-69, (1998).

130.    Vasantha Kandasamy W.B., *Group rings which have trivial subrings*, The Math. Edu., Vol. XXXIII, 180-181, (1999).

131.    Vasantha Kandasamy W.B., *On E-semi group rings*, Caribbean J. of Math. and Comp. Sci., Vol. 9, 52-54, (1999).

132.    Vasantha Kandasamy W.B., *On demi- modules over rings*, J. of Wuhan Automotive Politechnic Univ., Vol. 22, 123-125, (2000).

133.    Vasantha Kandasamy W.B., *On finite quaternion rings and skew fields*, Acta Ciencia Indica, Vol. XXIV, 133-135, (2000).

134.    Vasantha Kandasamy W.B., *On group rings which are $\gamma_n$ rings*, The Math. Edu., Vol. XXXIV, 61, (2000).

135.    Vasantha Kandasamy W.B., *CN rings*, Octogon, Vol.9, 343-344, (2001).

136.    Vasantha Kandasamy W.B., *On locally semi unitary rings*, Octogon, Vol.9, 260-262, (2001).

137.    Vasantha Kandasamy W.B., *Tight rings and group rings,* Acta Ciencia Indica, Vol. XXVII, 87-88, (2001).

138.    Vasantha Kandasamy W.B., *Smarandache Semigroups*, American Research Press, Rehoboth, NM, (2002).

139.    Vasantha Kandasamy W.B., *On Smarandache pseudo ideals in rings*, (2002).
        http://www.gallup.unm.edu/~smaranandache/pseudoideals.pdf

140.    Vasantha Kandasamy W.B., *Smarandache Zero divisors*, (2002).
        http://www.gallup.unm.edu/~smarandache/ZeroDivisor.pdf

141.    Vasantha Kandasamy W.B., *Finite zeros and finite zero-divisors*, Varahmihir J. of Math. Sci., Vol. 2, (To appear), 2002.

142.    Vasantha Kandasamy W.B. and Smarandache, F., *Basic Neutrosophic and their Applications to Fuzzy and Neutrosophic models*, Hexis, Church Rock, 2004.

143.    Vasantha Kandasamy W.B. and Smarandache, F., *Some Neutrosophic Algebraic Structures and Neutrosophic N-Algebraic Structures*, Hexis, Church Rock, 2006.

144.    Victoria Powers, *Higher level orders on non-commutative rings*, J. Pure and Appl. Algebra, Vol. 67, 285-298, (1990).

145.    Vougiouklis Thomas, *On rings with zero divisors strong V-groups*, Comment Math. Univ. Carolin J., Vol. 31, 431-433, (1990).

146.    Wilson John. S., *A note on additive subgroups of finite rings*, J. Algebra, No. 234, 362-366, (2000).

147.    Yakub Adil, *Structure of weakly periodic rings with potent extended commutators*, Int. J. of Math. Sci., Vol. 25, 299-304, (2001).

148.    Zariski.O and Samuel.P, *Commutative Algebra*, Van Nostrand Reinhold, (1958).

149.    Zhang-Chang quan, *Inner Commutative Rings*, Sictiuan Daxue Xuebao, Vol. 26, 95-97, (1989).

# INDEX

**B**

**P**

# ABOUT THE AUTHORS

**Dr.W.B.Vasantha Kandasamy** is an Associate Professor in the Department of Mathematics, Indian Institute of Technology Madras, Chennai, where she lives with her husband Dr.K.Kandasamy and daughters Meena and Kama. Her current interests include Smarandache algebraic structures, fuzzy theory, coding/ communication theory. In the past decade she has guided 11 Ph.D. scholars in the different fields of non-associative algebras, algebraic coding theory, transportation theory, fuzzy groups, and applications of fuzzy theory of the problems faced in chemical industries and cement industries. Currently, four Ph.D. scholars are working under her guidance.

She has to her credit 633 research papers of which 211 are individually authored. Apart from this, she and her students have presented around 328 papers in national and international conferences. She teaches both undergraduate and post-graduate students and has guided over 51 M.Sc. and M.Tech. projects. She has worked in collaboration projects with the Indian Space Research Organization and with the Tamil Nadu State AIDS Control Society. This is her 27th book.

She can be contacted at vasantha@iitm.ac.in
You can visit her work on the web at: http://mat.iitm.ac.in/~wbv
www.vasantha.net

---

**Dr. Dr.Florentin Smarandache** is an Associate Professor of Mathematics at the University of New Mexico in USA. He published over 75 books and 100 articles and notes in mathematics, physics, philosophy, psychology, literature, rebus. In mathematics his research is in number theory, non-Euclidean geometry, synthetic geometry, algebraic structures, statistics, neutrosophic logic and set (generalizations of fuzzy logic and set respectively), neutrosophic probability (generalization of classical and imprecise probability). Also, small contributions to nuclear and particle physics, information fusion, neutrosophy (a generalization of dialectics), law of sensations and stimuli, etc.

He can be contacted at smarand@unm.edu